\documentclass[opre,nonblindrev]{informs3_modified1}

\DoubleSpacedXI 

\usepackage{endnotes}
\let\footnote=\endnote
\let\enotesize=\normalsize
\def\notesname{Endnotes}%
\def\makeenmark{\hbox to1.275em{\theenmark.\enskip\hss}}
\def\enoteformat{\rightskip0pt\leftskip0pt\parindent=1.275em
  \leavevmode\llap{\makeenmark}}

\usepackage{natbib}
 \bibpunct[, ]{(}{)}{,}{a}{}{,}%
 \def\bibfont{\small}%

\usepackage{amsmath,amsfonts,amssymb,graphicx}
\usepackage{bm}
\usepackage{bbm}
\usepackage{algorithm}
\usepackage{algorithmic}
\usepackage{subcaption}
\usepackage{amsfonts}
\usepackage{amsbsy}
\usepackage{enumitem}
\usepackage{amsmath}
\usepackage{multirow}

\TheoremsNumberedThrough     
\ECRepeatTheorems

\EquationsNumberedThrough    

\newlength\myindent 
\newcommand\bindent{%
  \begingroup 
  \setlength{\itemindent}{\myindent} 
  \addtolength{\algorithmicindent}{\myindent} 
}
\newcommand\eindent{\endgroup} 
\begin{document}


\RUNAUTHOR{Lam and Zhang}

\RUNTITLE{Distributionally Constrained Black-Box Stochastic Gradient Estimation and Optimization}

\TITLE{Distributionally Constrained Black-Box Stochastic Gradient Estimation and Optimization}

\ARTICLEAUTHORS{%
\AUTHOR{Henry Lam}
\AFF{Department of Industrial Engineering and Operations Research, Columbia University, New York, NY 10027, \EMAIL{henry.lam@columbia.edu}} 
\AUTHOR{Junhui Zhang}
\AFF{Department of Applied Physics and Applied Mathematics, Columbia University, New York, NY 10027, \EMAIL{jz2903@columbia.edu}} 
} 

\ABSTRACT{%
We consider stochastic gradient estimation using only black-box function evaluations, where the function argument lies within a probability simplex. This problem is motivated from gradient-descent optimization procedures in multiple applications in distributionally robust analysis and inverse model calibration  involving decision variables that are probability distributions. We are especially interested in obtaining gradient estimators where one or few sample observations or simulation runs apply simultaneously to all directions. 
Conventional zeroth-order gradient schemes such as simultaneous perturbation face challenges as the required moment conditions that allow the ``canceling" of higher-order biases cannot be satisfied without violating the simplex constraints. We investigate a new set of required conditions on the random perturbation generator, which leads us to a class of implementable gradient estimators using Dirichlet mixtures. We study the statistical properties of these estimators and their utility in constrained stochastic approximation, including both Frank-Wolfe and mirror descent update schemes. We demonstrate the effectiveness of our procedures and compare with benchmarks via several numerical examples.
}%



\maketitle

%


\section{Introduction}
\label{sec:intro}
We study stochastic gradient estimation with noisy black-box function evaluation, where the function argument lies within a probability simplex. More precisely, given a real-valued function $Z(\cdot)$ defined on an $n$-dimensional set of probability weights, we want to compute its (properly defined) gradient $\nabla Z$, by using only noisy observations of $Z$ at suitably chosen points. Our particular focus is on schemes that can simultaneously estimate all components of the gradient with only a single or a pair of observations on $Z$. 

Gradient estimation with only noisy function evaluation arises ubiquitously in stochastic or simulation-based optimization, when one uses stochastic approximation (SA) or gradient descent and the model is only accessible as a black box (e.g., \citealt{glasserman2013monte,asmussen2007stochastic,fu2006gradient,l1991overview}). In derivative-free optimization, this is also known as the zeroth-order gradient oracle \citep{ghadimi2013stochastic}. The classical Kiefer-Wolfowitz (KW) scheme \citep{kushner2003stochastic} is applicable precisely in the above optimization setting. At each iteration, KW uses a pair of noisy evaluation to compute a finite difference that approximates the partial derivative at each dimension, and convergence is provably attained via suitable decays of the perturbation size and iteration step size. The number of evaluations at each iteration required by this scheme, however, grows with the problem dimension. To address this, schemes that allow a simultaneous estimation of derivatives at all dimensions have been investigated. These schemes randomly generate a perturbation vector according to some well-chosen distribution, and evaluate the function at the realized perturbation. Then, through multiplying by a factor related to the generating distribution (and possibly taking differences of these evaluations), an estimate of the entire gradient is obtained. Prominent examples of such approaches include the simultaneous perturbation stochastic approximation (SPSA) \citep{spall1992multivariate,spall1997one} that uses distributions with enough masses at zero, Gaussian smoothing \citep{nesterov2017random,maggiar2018derivative,ghadimi2013stochastic} that uses multivariate Gaussian distributions, and uniform sampling \citep{flaxman2004online} to generate the perturbation vector. Recent works further propose the generation of dependent perturbation vectors, such as orthogonal sampling \citep{choromanski2018structured,rowland2018geometrically}, which can be viewed as a middle ground between classical finite difference with perturbation along each dimension and randomized perturbation, and are shown to improve estimation accuracy. Comprehensive comparisons among these approaches in terms of optimization efficiency can be found in \cite{berahas2019linear,berahas2019theoretical}.

In this paper we are interested in designing random-perturbation gradient estimators described above when the function input, or the decision variable in the corresponding optimization, is a probability distribution. Such type of optimization appears in two recent applications. One is \emph{distributionally robust optimization (DRO)}, which comprises the computation of robust bounds or decisions for stochastic problems under uncertainty. More specifically, this method advocates performance evaluation under the worst-case scenario among all input probability distributions within a so-called uncertainty set or ambiguity set \citep{wiesemann2014distributionally,delage2010distributionally,ben2013robust,goh2010distributionally,bertsimas2018robust} that is postulated to contain the truth with high likelihood. This idea roots in robust optimization \citep{ben2009robust,bertsimas2011theory} to handle optimization with uncertain parameter, in this case the underlying probability distribution. DRO has gained substantial popularity in recent years, with applications in control \citep{petersen2000minimax}, economics \citep{hansen2008robustness}, finance \citep{glasserman2014robust}, dynamic pricing \citep{lim2007relative} and queueing \citep{jain2010optimality}, using various uncertainty sets such as neighborhood balls with $\phi$-divergence \citep{ben2013robust,bayraksan2015data,jiang2016data,lam2016robust,lam2018sensitivity,gotoh2018robust,duchi2021}, Renyi divergence \citep{atar2015robust,blanchet2020distributionally} and Wasserstein distance \citep{blanchet2019quantifying,gao2016distributionally,chen2018robust, esfahani2018data}, or moment \citep{ghaoui2003worst,delage2010distributionally,goh2010distributionally,wiesemann2014distributionally,hanasusanto2015distributionally}, shape  \citep{popescu2005semidefinite,van2016generalized,li2016ambiguous,lam2017tail} and marginal \citep{chen2018distributionally} constraints. 

In the context of stochastic simulation, DRO entails an approach to quantify the so-called input uncertainty. To explain the latter, simulation modeling consists of feeding random variates drawn from input distributions which, through the system logic, generates outputs for subsequent performance analyses. When the input distributions are corrupted or calibrated from data, the statistical error can propagate to the outputs, a problem known as input uncertainty (see the surveys  \citealt{barton2012tutorial,song2014advanced,henderson2003input,chick2006bayesian,lam2016advanced,corlu2020stochastic}). Under such situations, DRO can be used to construct performance bounds that are especially competitive in high-uncertainty, nonparametric situations where the modeler may only have minimal information on the input distributions \citep{hu2012robust,glasserman2014robust,ghosh2019robust,ghosh2015mirror,lam2016robust,hu2015robust}, or when the modeler wants to address uncertainty in directions beyond typical assumptions, such as dependency in the input processes \citep{lam2018sensitivity}. This approach also can be viewed as an optimization-based alternative to obtain statistically accurate output confidence intervals via its connection to the empirical likelihood \citep{lam2017optimization}. These usages are in contrast to traditional statistical methods such as bootstrap resampling \citep{barton1993uniform,barton2001resampling,cheng1997sensitivity,barton2013quantifying,song2015quickly,lam2018subsampling}, the delta method \citep{cheng1998two,lin2015single}, or Bayesian posterior sampling \citep{chick2001input,zouaoui2003accounting,zouaoui2004accounting,xie2014bayesian,zhu2020risk}. In order to compute the worst-case values, one needs to run gradient descents on the probability simplex for which our estimator will be useful. 

The second application is \emph{inverse model calibration}. This refers to the calibration of stochastic input models from output-level data, a task faced by various operational settings where direct input data are unavailable due to administrative or data collection reasons. It is closely related to so-called model validation that aims to assess the match between simulation and real-world data \citep{sargent2010verification,kleijnen1995verification}. In scientific applications, it is also known as simulation-based or likelihood-free inference \citep{cranmer2020frontier} as the likelihood of the target quantity is only accessible via black-box simulation models, and is viewed as an inverse problem \citep{robert2016approximate,tarantola2005inverse}. One approach is to minimize a distance between the postulated model class and real-world data, where the distance metric leads to criteria such as regularized least square \citep{tarantola2005inverse}, maximum entropy \citep{donoho1992maximum,csiszar1991least,goeva2014reconstructing} and Jensen-Shannon divergence minimization \citep{louppe2019adversarial}. A more recent approach utilizes the distributionally robust framework described above, where the input-output match is incorporated as constraints in the uncertainty set of a suitable DRO problem \citep{OR2019,bai2020distributionally}. When the input parameters are probability distributions, the involved distance minimization or robust optimization leads to a distributional optimization problem in which our gradient descent schemes are again relevant.

Besides the two main applications above, optimization over probability distributions is potentially useful in other contexts. Examples include matrix games that optimize randomized strategies \citep{nemirovski2009robust}. Another potential context is reinforcement learning, when using policy gradients embedded in nonlinear optimization schemes \citep{schulman2015trust,choromanski2018structured,fazel2018global}. Procedures such as the evolutionary strategy \citep{salimans2017evolution} advocate the use of derivative-free optimization to overcome the high variance in long-horizon problems exhibited by conventional score-function-based policy gradients. Our stochastic gradient estimator applies as a derivative-free approach when the randomized policies are represented via discrete probability distributions on the action space.

In all the above applications, we stochastically optimize over decision variables that are probability distributions. Moreover, the objective function involves a black-box simulation model where we can only observe noisy model outputs, thus susceptible only to biased gradient estimators \citep{zazanis1993convergence,fox1989replication,frolov1963calculation} in contrast to unbiased alternatives (such as the infinitesimal perturbation analysis \citep{heidelberger1988convergence,ho1983infinitesimal,PG1991,l1990unified}, the likelihood ratio or the score function method \citep{glynn1990likelihood,rubinstein1986score,reiman1989sensitivity}, measure-valued differentiation \citep{heidergott2008measure,heidergott2010gradient}, or other variants \citep{fu1992extensions,rubinstein1992sensitivity,hong2009estimating,peng2018new}. This motivates us to investigate black-box gradient estimators evaluated over probability simplices. Note that, as the decision variable is the input distribution, the objective is typically undefined or incomputable when the decision variable is outside this feasible region. Thus, in order to successfully run the zeroth-order gradient-based iteration, we need to obtain a gradient estimator on the constrained space, which is precisely our focus.



A challenge in obtaining gradient estimator via random perturbation, compared to previously studied unconstrained settings, lies in the satisfaction of the distributional and moment conditions required for the random perturbation vector. These conditions ensure that the higher-order error terms in the function evaluation or finite difference under the perturbation can be sufficiently ``canceled out", so that the resulting gradient estimator enjoys good bias properties. For instance, the standard version of SPSA requires the perturbation vector to have independent, mean-zero, symmetric components with finite reciprocal moments \citep{spall1992multivariate}. Because of the constraints in our problem, however, the perturbation vector we generate must have correlated and unsymmetric components, and with restrictive moment behaviors, thus making these previous conditions difficult to satisfy. Similar issues arise one way or the other in other similar randomized schemes such as Gaussian smoothing and uniform sampling.  Moreover, challenges also arise in using the idea of central finite-differencing in these schemes, as the probability constraint may, depending on the evaluation point of interest, only allow movement in either the forward or backward direction. 

To overcome all these challenges, in this paper we present a new set of conditions to ensure that the perturbation vector within the probability simplex leads to higher-order cancellation. We then propose a specific set of estimators using the Dirichlet mixture that satisfies these new conditions. We show how these estimators are readily implementable by offering precise guidelines. We also demonstrate their effectiveness in both estimating gradients and in applying to SA. In particular, we consider both the Frank-Wolfe (FW) \citep{reddi2016stochastic,lafond2015online} and the mirror descent (MD) \citep{beck2003mirror} counterparts in addressing the constraints, and show how correctly tuning our gradient estimators leads to provable convergences in both schemes. We also compare our schemes with simple finite differences that perturb along randomly chosen dimensions, and discuss the similarities and distinctions with gradient estimators under unconstrained domains \citep{berahas2019theoretical}.


The rest of the paper is as follows. Section \ref{sec:problem} introduces our problem setting and proposed framework. Section \ref{sec:analysis} analyzes the bias and variance properties of our gradient estimators. Section \ref{sec:Dir} constructs explicit configurations of our estimators using Dirichlet mixtures. Section \ref{sec:comparisons} compares among our constructions, with naive finite differences, and unconstrained gradient estimators. Section \ref{sec:SA} studies FWSA and MDSA using our estimators. Section \ref{sec:numerics} demonstrates numerical results on both the qualities of gradient estimation and optimization procedures. All technical proofs are deferred to the Appendix. Finally, a preliminary conference version of this work appears in \cite{baedistributionally}, which only contains the basic idea of the proposed gradient estimators and simple numerical illustrations. All bias-variance analyses on the estimators, algorithmic details and guarantees of optimization procedures, and the elaborate numerical studies in this work are beyond \cite{baedistributionally}.


\section{Problem Setting and Estimation Framework}\label{sec:problem}
Consider the function $Z(\cdot):\mathcal P\to\mathbb R$, where $$\mathcal P=\{\mathbf p=(p_1,\ldots,p_n)\in\mathbb R^n:p_1+\cdots+p_n=1,p_i\geq0\ \forall i=1,\ldots,n\}$$ 
is an $n$-dimensional probability simplex. We are interested in estimating $\nabla Z(\mathbf p)$, the gradient of $Z$. To be more precise, note that, because of the probability-simplex constraint, an arbitrary perturbation of $\mathbf p$ needed to define a gradient may shoot outside the domain of $Z$, and thus, when speaking of gradient, we mean a directional gradient $\nabla Z(\mathbf p)\in\mathbb R^n$ defined such that, for any perturbation of $\mathbf p\in\mathcal P$ to $\mathbf q\in\mathcal P$, we have
\begin{equation}
Z(\mathbf q)-Z(\mathbf p)=\nabla Z(\mathbf p)'(\mathbf q-\mathbf p)+o(\|\mathbf q-\mathbf p\|)\label{grad def}
\end{equation}

Here and below, we use $\|\cdot\|$ to denote the $L_2$ norm of a vector and the spectral norm of a matrix. 

Equivalently, we can define each component of $\nabla Z(\mathbf p)=(Z_1(\mathbf p),\ldots,Z_n(\mathbf p))'$ as 
\begin{equation}
Z_i(\mathbf p)=\lim_{\epsilon\to0_+}\frac{Z((1-\epsilon)\mathbf p+\epsilon\mathbf e_i)-Z(\mathbf p)}{\epsilon}\label{grad def1}
\end{equation}
where $\mathbf e_i$ is a vector that takes value 1 at component $i$ and 0 otherwise. That is, $Z_i(\mathbf p)$ is defined by mixing $\mathbf p$ with a point mass at component $i$ and evaluating the increment as the mixture parameter shrink to 0. It is straightforward to see that \eqref{grad def} and \eqref{grad def1} coincide, up to a translation of $c\mathbbm 1$ for any constant $c$ and $\mathbbm 1$ being the vector of 1's (note that \eqref{grad def} as defined is non-unique up to this translation since $c\mathbbm 1'(\mathbf q-\mathbf p)=0$ for any probability vectors $\mathbf p$ and $\mathbf q$, whereas \eqref{grad def1} is unique). In robust statistics, definition \eqref{grad def1} represents the so-called influence function \citep{hampel2011robust} that measures the sensitivity with respect to changes in the data distribution, when viewing $Z(\mathbf p)$ as a statistical functional on the data distribution $\mathbf p$.

Given this definition of the gradient, we can also define higher-order derivatives recursively via
\begin{equation*}
Z_{i_1,\ldots,i_k,i_{k+1}}(\mathbf p)=\lim_{\epsilon\to0_+}\frac{Z_{i_1,\ldots,i_k}((1-\epsilon)\mathbf p+\epsilon\mathbf e_{i_{k+1}})-Z_{i_1,\ldots,i_k}(\mathbf p)}{\epsilon}
\end{equation*}
In fact, using these definitions, we can get the corresponding Taylor expansion formula:
\begin{lemma}
Suppose $Z(\cdot):\mathcal P\to\mathbb R$ is $\mathcal C^{k}$ and $Z_{i_1,\ldots,i_{k+1}}$ exists for all $(i_1,\ldots,i_{k+1}) \in [n]^{k+1} $ for every $\mathbf p \in \mathcal P$. Let $\mathbf p_1 = (p_{1,1},\ldots,p_{1,n})$ and $\mathbf p_2  = (p_{2,1},\ldots,p_{2,n})\in \mathcal P$. Then there exists $\xi \in [0,1]$ such that
\begin{align*}
        Z(\mathbf p_2) &= Z(\mathbf p_1) + \sum_{l=1}^{k}\frac{1}{l!}\sum_{(i_1,\ldots,i_l)\in [n]^{l}}Z_{i_1,\ldots i_l}(\mathbf p_{1})\prod_{j=1}^l ( p_{2,i_j}- p_{1,i_j})\\ 
        &+\frac{1}{(k+1)!}\sum_{(i_1,\ldots,i_{k+1})\in [n]^{k+1}}Z_{i_1,\ldots i_{k+1}}((1-\xi)\mathbf p_{1} +\xi \mathbf p_2)\prod_{j=1}^{k+1} ( p_{2,i_j}- p_{1,i_j})
\end{align*}
\label{lemma:taylor}
\end{lemma}
where $[N],N\in \mathbbm N_{+}$ denotes $\{1,2,\ldots,N\}$. 

This result will be useful later in the analysis of the bias of the gradient estimators.

Our goal is to estimate $\nabla Z(\mathbf p)$ when given only the capability to output noisy function evaluations. That is, given any $\mathbf q$, we can output $\hat Z(\mathbf q)$ such that $E[\hat Z(\mathbf q)]=Z(\mathbf q)$. Note that, in many cases, evaluating $\hat Z(\cdot)$ can be costly, and the naive finite-difference method that approximates the gradient via the right hand side of \eqref{grad def1} with a small $\epsilon>0$ leads to $n+1$ number of function evaluations, which can be computationally costly in high dimensions. To overcome this, we consider estimating $\nabla Z(\mathbf p)$ by randomly generating a probability vector, say $\boldsymbol\delta$, and mixing it with the original probability vector $\mathbf p$. This can be thought as using a perturbation vector $\boldsymbol\delta-\mathbf p$. Our scheme then follows by evaluating $\hat Z$ at the perturbed probability vector and ``projecting" properly to get derivative estimates simultaneously for all dimensions. We consider three variants of gradient estimators, and will discuss the faced statistical challenges and their properties in the next section.
\\



\noindent\emph{Single Function Evaluation (SFE). }Our most basic scheme is to use
\begin{equation}
\frac{\hat Z((1-c)\mathbf p+c\boldsymbol\delta)}{c}S(\mathbf p,\boldsymbol\delta)\label{est 1}
\end{equation}
where $c>0$ and $S(\mathbf p,\boldsymbol\delta)\in\mathbb R^n$ are chosen by the user. When $R$ simulation budget is available, we would randomly draw $\boldsymbol\delta^j,j=1,\ldots,R$, evaluate \eqref{est 1} at each drawn $\boldsymbol\delta^j$ with an independent evaluation of $\hat Z$, and at the end output their average, i.e.,
\begin{equation}
\hat\psi_{SFE}=\frac{1}{R}\sum_{j=1}^R\frac{\hat Z^j((1-c)\mathbf p+c\boldsymbol\delta^j)}{c}S(\mathbf p,\boldsymbol\delta^j)\label{est 1 all}
\end{equation}
where $\hat Z^j(\cdot),j=1,\ldots,R$ denote $R$ independent function evaluations. We call \eqref{est 1 all} the SFE estimator since each replication that contributes to the gradient estimate only needs to evaluate $\hat Z$ once.

In \eqref{est 1} and \eqref{est 1 all}, $c>0$ is a perturbation size that could depend on the dimension of the problem $n$. The vector $S(\mathbf p,\boldsymbol\delta)$ represents a factor or ``score function" that aims to cancel out higher-order bias terms in the associated Taylor series.
\\

\noindent\emph{Forward Function Evaluation (FFE).}
Next we consider the use of two function evaluations, one at $(1-c)\mathbf p+c\boldsymbol\delta$ and one at $\mathbf p$ itself. We output
\begin{equation}
\frac{\hat Z^1((1-c)\mathbf p+c\boldsymbol\delta)-\hat Z^2(\mathbf p)}{c}S(\mathbf p,\boldsymbol\delta)\label{est 2}
\end{equation}
where $\hat Z^1(\cdot)$ and $\hat Z^2(\cdot)$ represent two independent evaluations, and $c>0$ and $S(\mathbf p,\boldsymbol\delta)\in\mathbb R^n$ are chosen by the user like before. Analogously, when $2R$ simulation budget is available, we would use
\begin{equation}
\hat\psi_{FFE}=\frac{1}{R}\sum_{j=1}^R\frac{\hat Z^{2j-1}((1-c)\mathbf p+c\boldsymbol\delta^j)-\hat Z^{2j}(\mathbf p)}{c}S(\mathbf p,\boldsymbol\delta^j)\label{est 2 all}
\end{equation}
where $\hat Z^{2j-1}(\cdot),\hat Z^{2j}(\cdot),j=1,\ldots,R$ are $2R$ independent function evaluations. 
\\

\noindent\emph{Central Function Evaluation (CFE)}. 
Lastly, we also consider a variant similar to central finite-differencing. Consider using two function evaluations, one at $(1-c)\mathbf p+c\boldsymbol\delta$ and one at $(1+c)\mathbf p-c\boldsymbol\delta$. We output
\begin{equation}
\frac{\hat Z^1((1-c)\mathbf p+c\boldsymbol\delta)-\hat Z^2((1+c)\mathbf p-c\boldsymbol\delta)}{2c}S(\mathbf p,\boldsymbol\delta)\label{est 3}
\end{equation}
where $\hat Z^1(\cdot)$ and $\hat Z^2(\cdot)$ again represent two independent evaluations. Like FFE, when $2R$ simulation budget is available, we would use
\begin{equation}
\hat\psi_{CFE}=\frac{1}{R}\sum_{j=1}^R\frac{\hat Z^{2j-1}((1-c)\mathbf p+c\boldsymbol\delta^j)-\hat Z^{2j}((1+c)\mathbf p-c\boldsymbol\delta^j)}{2c}S(\mathbf p,\boldsymbol\delta^j)\label{est 3 all}
\end{equation}
where $\hat Z^{2j-1}(\cdot),\hat Z^{2j}(\cdot),j=1,\ldots,R$ are $2R$ independent function evaluations. 

\section{Challenges and Statistical Properties}\label{sec:analysis}
We investigate and compare the statistical errors of the above estimators in estimating the gradient. In doing so we will also highlight the involved challenges. Typically, all the estimators above have variances that scale in order $1/(Rc^2)$, as the $c$ appears in the denominators of these estimators. To understand the bias, supposing $Z$ is sufficiently smooth, Lemma \ref{lemma:taylor} gives the Taylor series
\begin{equation*}
Z((1-c)\mathbf p+c\boldsymbol\delta)=Z(\mathbf p)+c\nabla Z(\mathbf p)'(\boldsymbol\delta-\mathbf p)+\frac{c^2}{2}(\boldsymbol\delta-\mathbf p)'\nabla^2Z(\mathbf p)(\boldsymbol\delta-\mathbf p)+O(c^3)
\end{equation*}
where $\nabla^2Z(\mathbf p)$ is a properly defined Hessian. Now, since $\boldsymbol\delta$ and the function evaluations are independent, we can write the expectation of SFE as
\begin{eqnarray}
&&E\left[\frac{Z((1-c)\mathbf p+c\boldsymbol\delta)}{c}S(\mathbf p,\boldsymbol\delta)\right]\notag\\
&=&E\left[\frac{1}{c}\left(Z(\mathbf p)+c\nabla Z(\mathbf p)'(\boldsymbol\delta-\mathbf p)+\frac{c^2}{2}(\boldsymbol\delta-\mathbf p)'\nabla^2Z(\mathbf p)(\boldsymbol\delta-\mathbf p)+O(c^3)\right)S(\mathbf p,\boldsymbol\delta)\right]\notag\\
&=&Z(\mathbf p)\frac{E[S(\mathbf p,\boldsymbol\delta)]}{c}+\nabla Z(\mathbf p)'E[(\boldsymbol\delta-\mathbf p)S(\mathbf p,\boldsymbol\delta)]+\frac{c}{2}E[tr(\nabla^2Z(\mathbf p)(\boldsymbol\delta-\mathbf p)(\boldsymbol\delta-\mathbf p)')S(\mathbf p,\boldsymbol\delta)]+O(c^2)\label{interim}
\end{eqnarray}
If we can ensure the following three moment conditions
\begin{align}
E[S(\mathbf p,\boldsymbol\delta)]&=\mathbf 0\tag{MC 1}\label{first order}\\
E[(\boldsymbol\delta-\mathbf p)S(\mathbf p,\boldsymbol\delta)']&=\mathbf I-\lambda \mathbbm 1\mathbbm 1'\tag{MC 2}\label{second order}\\
E[(\boldsymbol\delta - \mathbf p)_i(\boldsymbol\delta -\mathbf p)_jS(\mathbf p,\boldsymbol\delta)_k]&=\mu,~\forall~i,j,k\in[n]\tag{MC 3}
\label{third order}
\end{align}
where $\lambda,\mu\in\mathbb R$, $\mathbf I$ is the identity matrix, $\mathbbm 1$ is a vector of ones, and $\mathbf 0$ is a zero matrix each in the suitable dimension, then SFE gives rise to a gradient estimator that has bias $O(c^2)$. As stated later in Theorem \ref{thm_necessary_conditions}, these conditions are in fact necessary for the $O(c^2)$ bias. Here the term $-\lambda \mathbbm 1\mathbbm 1'$ for any constant $\lambda$ appears in \eqref{second order}, and constant $\mu$ appears in \eqref{third order}. When setting both $\lambda$ and $\mu$ to 0, it is clear that the right hand side of \eqref{interim} becomes $\nabla Z(\mathbf p)+O(c^2)$. The non-zero choices of $\lambda$ and $\mu$ are allowed as they also determine $\nabla Z(\mathbf p)$ up to a translational constant with the same error order. Importantly, the non-zero choices of $\lambda$ and $\mu$ provide more flexibility in our estimator tuning to attain a small-order error. If, however, \eqref{third order} cannot be satisfied, then the bias could be $O(c)$, which is less desirable. There is also a possibility that \eqref{third order} is not exact equality but very close, which we will see to be a useful compromise. 

We make our discussions above concrete by deriving the bias and variance of $\hat \psi_{SFE}$. Here we choose $S(\mathbf p,\boldsymbol\delta) = \gamma (\boldsymbol\delta - \mathbf p)$, for some $\gamma\in\mathbbm R$, as it turns out to be a simple yet good choice for bias cancellation. We assume the following function class:
\begin{assumption}[Smoothness of objective function]
Assume that $Z(\cdot)$ is $\mathcal C^2$. In addition, assume that $Z(\cdot)$ is $L_0$-Lipschitz (i.e.,$\lvert  Z(\mathbf p_1)-Z(\mathbf p_2)\rvert \leq L_0\lVert \mathbf p_1 - \mathbf p_2 \rVert$), its gradient is $L_1$-Lipschitz (i.e., $\lVert \nabla Z(\mathbf p_1)-\nabla Z(\mathbf p_2)\rVert \leq L_1\lVert \mathbf p_1 - \mathbf p_2 \rVert$), and its Hessian is $L_2$-Lipschitz (i.e., $\lVert \nabla^2 Z(\mathbf p_1)-\nabla^2 Z(\mathbf p_2)\rVert \leq L_2\lVert \mathbf p_1 - \mathbf p_2 \rVert$). \label{assumption:smooth}
\end{assumption}

Under Assumption \ref{assumption:smooth}, we further define constants $M_0, M_1, M_2$ that depend on $L_0, L_1, L_2$ as $M_0 =\sup_{\mathbf p\in \mathcal P} \lvert Z(\mathbf p) \rvert $, $M_1 = \sup_{\mathbf p\in \mathcal P} \lVert \nabla Z(\mathbf p) \rVert$, and $M_2 = \sup_{\mathbf p\in \mathcal P} \lVert \nabla^2 Z(\mathbf p) \rVert$. The sup exists since $Z$, $\nabla Z$, and $\nabla^2 Z$ are Lipschitz and $\mathcal P$ is bounded.

Since $L_0 = \sup_{\mathbf p\neq \mathbf q\in \mathcal P} \frac{\lvert Z(\mathbf p)-Z(\mathbf p)\rvert}{\lVert\mathbf p-\mathbf q\rVert}=\sup_{\mathbf p\neq \mathbf q\in \mathcal P} \frac{\lvert \nabla Z(\eta \mathbf p+(1-\eta)\mathbf q)'(\mathbf p-\mathbf q)\rvert}{\lVert\mathbf p-\mathbf q\rVert}\leq \sup_{\mathbf p\neq \mathbf q\in \mathcal P} \frac{\lVert \nabla Z(\eta \mathbf p+(1-\eta)\mathbf q)\rVert\cdot \lVert \mathbf p-\mathbf q\rVert }{\lVert\mathbf p-\mathbf q\rVert}=\sup_{\mathbf p,\mathbf q\in \mathcal P} \lVert \nabla Z(\eta \mathbf p+(1-\eta)\mathbf q)\rVert\leq M_1$, where $\eta\in[0,1]$ depends on $\mathbf p,\mathbf q$. Similarly, we can get the relation that $L_1\leq M_2$. 

\begin{assumption}[Bias and variance of the function evaluation]
Assume that the function evaluation is unbiased, i.e., $E[\hat{Z}(\mathbf p)] = Z(\mathbf p)$, and its variance $\sigma^2(\mathbf p):= Var[\hat{Z}(\mathbf p)]$ is bounded by $\sigma^2>0$.
 \label{assumption:smooth_noise}
\end{assumption}

Denote $\psi_{SFE}:=E[\hat \psi_{SFE}]$. We have the following:
\begin{theorem}[Bias and variance of SFE]
Under Assumptions \ref{assumption:smooth} and \ref{assumption:smooth_noise}, we choose $\boldsymbol\delta$ and $S(\mathbf p,\boldsymbol\delta) = \gamma (\boldsymbol\delta - \mathbf p)$ that satisfy \eqref{first order} and \eqref{second order}. Then $\hat \psi_{SFE}$ has a bias given by
\begin{equation}
\psi_{SFE}-\nabla Z(\mathbf p)-\epsilon_0 \mathbbm 1=\mathbf \epsilon \label{SFE bias}
\end{equation}
with $ \epsilon_0 = -\frac{\nabla Z(\mathbf p)'\mathbbm 1}{n}$ and $\lVert \mathbf \epsilon \rVert = O(ncM_2+nc^2L_2)$, and variance given by
\begin{equation*}
E[\lVert \hat \psi_{SFE}-\psi_{SFE}\rVert^2] =O\left(\frac{\gamma n }{Rc^2}[\sigma^2+M_0^2]\right)
\end{equation*}
Moreover, if \eqref{third order} is also satisfied, then
\begin{equation*}
\psi_{SFE}-\nabla Z(\mathbf p)-\epsilon_0 \mathbbm 1=\mathbf \epsilon
\end{equation*}
with $\epsilon_0 = \frac{c\mu}{2}\mathbbm 1'\nabla^2 Z(\mathbf p) \mathbbm 1-\frac{\nabla Z(\mathbf p)'\mathbbm 1}{n}$ and $\lVert \mathbf \epsilon \rVert = O(nc^2L_2)$.
Here $M_0, M_2>0$ are constants that depend on $L_0, L_1, L_2$. In the above, $\nabla Z(\mathbf p)$ refers to the unique version in \eqref{grad def1}.
\label{thm:SFE}
\end{theorem}
Note that \eqref{SFE bias} states that $E[\hat\psi_{SFE}]$ determines the canonical version of $\nabla Z(\mathbf p)$ in \eqref{grad def1} up to the constant $-\frac{\nabla Z(\mathbf p)'\mathbbm 1}{n}$ or $\frac{c\mu}{2}\mathbbm 1'\nabla Z(\mathbf p)\mathbbm 1-\frac{\mathbbm 1'\nabla Z(\mathbf p)}{n}$, depending on how many of \eqref{first order} to \eqref{third order} are satisfied. 
If the configuration is chosen such that only \eqref{first order} and \eqref{second order} are satisfied, then the bias with respect to this canonical $\nabla Z(\mathbf p)$ is of order $O(c)$, whereas if \eqref{third order} is also satisfied, then the bias is improved to order $O(c^2)$.


One might argue that the three conditions \eqref{first order}, \eqref{second order} and \eqref{third order} can be challenging to be satisfied simultaneously, and wonder if all of them are necessary for the bias to be of order $O(c^2)$ in estimating $\nabla Z(\mathbf p)$ up to some translation. We show next that these conditions are indeed necessary.
\begin{theorem}[Necessary conditions for bias cancellation]
If we choose $S(\mathbf p,\boldsymbol\delta) = \gamma ( \boldsymbol\delta-\mathbf p)$, conditions \eqref{first order}, \eqref{second order}, and \eqref{third order} are necessary for $\psi_{SFE}(\mathbf p)-\nabla Z(\mathbf p) -\epsilon_0(\mathbf p) \mathbbm 1= O(c^2)$ for some $\epsilon_0(\mathbf p) \in \mathbbm R$ that depends on $\mathbf p$.
\label{thm_necessary_conditions}
\end{theorem}

Next we consider FFE. It follows a similar, though slightly different, behavior than SFE. Its expectation is
\begin{eqnarray*}
&&E\left[\frac{Z((1-c)\mathbf p+c\boldsymbol\delta)- Z(\mathbf p)}{c}S(\mathbf p,\boldsymbol\delta)\right]\\
&=&E[S(\mathbf p,\boldsymbol\delta)(\boldsymbol\delta-\mathbf p)']\nabla Z(\mathbf p)+\frac{c}{2}E[(\boldsymbol\delta-\mathbf p)'\nabla^2Z(\mathbf p)(\boldsymbol\delta-\mathbf p)S(\mathbf p,\boldsymbol\delta)]+O(c^2)
\end{eqnarray*}
Compared with \eqref{interim}, the term $Z(\mathbf p)E[S(\mathbf p,\boldsymbol\delta)]/c$ disappears. Thus, if we can ensure \eqref{second order} and \eqref{third order} hold, then we have bias $O(c^2)$. Similarly, if only \eqref{second order} is enforced, then bias is $O(c)$. Denote $\psi_{FFE}:=E[\hat \psi_{FFE}]$. We summarize these as:
\begin{theorem}[Bias and variance of FFE]
Under Assumptions \ref{assumption:smooth} and \ref{assumption:smooth_noise}, we choose $\boldsymbol\delta$ and $S(\mathbf p,\boldsymbol\delta) = \gamma (\boldsymbol\delta - \mathbf p)$ that satisfy \eqref{first order} and \eqref{second order}. Then $\psi_{FFE}$ has a bias given by
\begin{equation*}
\psi_{FFE}-\nabla Z(\mathbf p)-\epsilon_0 \mathbbm 1=\mathbf \epsilon 
\end{equation*}
with $ \epsilon_0 = -\frac{\nabla Z(\mathbf p)'\mathbbm 1}{n}$ and $\lVert \mathbf \epsilon \rVert = O(ncM_2+nc^2L_2)$, and variance given by
\begin{equation*}
E[\lVert \hat \psi_{FFE}-\psi_{FFE}\rVert^2] =O(\frac{\gamma n }{Rc^2}[\sigma^2+c^2M_1^2+c^4L_1^2])
\end{equation*}
Moreover, if \eqref{third order} is also satisfied, then
\begin{equation*}
\psi_{FFE}-\nabla Z(\mathbf p)-\epsilon_0 \mathbbm 1=\mathbf \epsilon
\end{equation*}
with $\epsilon_0 = \frac{c\mu}{2}\mathbbm 1'\nabla^2 Z(\mathbf p) \mathbbm 1-\frac{\nabla Z(\mathbf p)'\mathbbm 1}{n}$ and $\lVert \mathbf \epsilon \rVert = O(nc^2L_2)$.
Here $M_2>0$ is a constant that depends on $L_0, L_1, L_2$. In the above, $\nabla Z(\mathbf p)$ refers to the unique version in \eqref{grad def1}.
\label{thm:FFE}
\end{theorem}

Since \eqref{first order} is easy to satisfy (by selecting any mean-zero $S(\mathbf p,\boldsymbol\delta)$), it may appear that SFE is preferable as it only requires one function evaluation to get one ``gradient sample". However, the first term in the Taylor expansion of $\hat\psi_{SFE}$ (not its mean $\psi_{SFE}$) is $Z(\mathbf p)S(\mathbf p,\boldsymbol\delta)/c$, which has a $c$ in the denominator, so with a finite simulation run this term could add significantly to the variance. This phenomenon is intuitive as SFE does not resemble a finite-difference scheme, the latter capturing the differencing needed in evaluating a derivative. In other words, by evaluating the objective function twice in FFE, we may substantially reduce the estimation variance. 

Next we discuss CFE. Thanks to the central finite-differencing, the first and third-order terms in the Taylor series of $Z((1-c)\mathbf p+c\boldsymbol\delta)$ and $Z((1+c)\mathbf p-c\boldsymbol\delta)$ cancel out, so that
$$E\left[\frac{Z((1-c)\mathbf p+c\boldsymbol\delta)- Z((1+c)\mathbf p-c\boldsymbol\delta)}{2c}S(\mathbf p,\boldsymbol\delta)\right]=E[S(\mathbf p,\boldsymbol\delta)(\boldsymbol\delta-\mathbf p)']\nabla Z(\mathbf p)+O(c^2)$$
Thus, to ensure a bias $O(c^2)$, we only require \eqref{second order}. Denote $\psi_{CFE}:=E[\hat \psi_{CFE}]$. We have the following result.

\begin{theorem}[Bias and variance of CFE]
Under Assumptions \ref{assumption:smooth} and \ref{assumption:smooth_noise}, we choose $\boldsymbol\delta$ and $S(\mathbf p,\boldsymbol\delta) = \gamma (\boldsymbol\delta - \mathbf p)$ that satisfy \eqref{first order} and \eqref{second order}. Then $\psi_{CFE}$ has a bias given by
\begin{equation*}
\psi_{CFE}-\nabla Z(\mathbf p)-\epsilon_0 \mathbbm 1=\mathbf \epsilon 
\end{equation*}
with $ \epsilon_0 = -\frac{\nabla Z(\mathbf p)'\mathbbm 1}{n}$ and $\lVert \mathbf \epsilon \rVert = O(nc^2L_2)$, and variance given by
\begin{equation*}
E[\lVert \hat \psi_{CFE}-\psi_{CFE}\rVert^2] =O(\frac{\gamma n }{Rc^2}[\sigma^2+c^2M_1^2+c^4L_1^2])
\end{equation*}
In the above, $\nabla Z(\mathbf p)$ refers to the unique version in \eqref{grad def1}.
\label{thm:CFE}
\end{theorem}
The fact that we only require \eqref{second order} in CFE to ensure bias $O(c^2)$ may appear to be desirable. However, note that $(1+c)\mathbf p-c\boldsymbol\delta$ is not a proper convex mixture and may not correspond to a probability distribution.

Thus, taking into account the variance and the implementability, FFE appears the best one among the three candidates. We will see that \eqref{first order} and \eqref{second order} are quite easy to enforce. In this sense, CFE can result in $O(c^2)$ bias, but it is not implementable generally due to the probability simplex constraint. SFE and FFE, on the other hand, can ensure bias $O(c)$ easily, and FFE can have a significantly better variance with a price of needing twice as many function evaluations. The main challenge, however, is that without the central finite-differencing trick, it is not obvious how to enforce \eqref{second order} and \eqref{third order} to achieve a bias $O(c^2)$ or to significantly diminish the $O(c)$ term. This would be the main investigation in the next section.

\section{Implementable Gradient Estimators Using Dirichlet Mixtures}\label{sec:Dir}
In view of the discussion in Section \ref{sec:problem}, we aim to get an estimator for $\nabla Z(\mathbf p)$ up to a translation, by using a random perturbation and associated score function that satisfy \eqref{first order}, \eqref{second order} and \eqref{third order}. In particular, we will base our choice on a mixture of Dirichlet distributions, i.e.,
\begin{equation*}
\boldsymbol\delta=\sum_{l=1}^K \theta_l\boldsymbol\delta^l,~~~\theta_l\geq 0,~\boldsymbol\delta^l\sim Dir(\boldsymbol\alpha^l),~\boldsymbol\alpha^l=(\alpha_{1}^l,\alpha_{2}^l,...,\alpha_{n}^l)
\end{equation*}
where $(\theta_1,\ldots,\theta_K)\in\mathbb R^K$ are the mixture probabilities. This $\boldsymbol\delta$ is guaranteed to lie in the probability simplex, and recall that we set the score function $S(\mathbf p,\boldsymbol\delta)=\gamma(\boldsymbol\delta-\mathbf p)$ for some $\gamma\in\mathbb R$. Without loss of generality, we assume from now on that $\mathbf p =(p_1,\ldots,p_n)$ such that $p_1\leq \ldots \leq p_n$.

\subsection{Preliminary Properties of the Dirichlet Distribution}
To facilitate our analyses, we list several useful facts for the Dirichlet distribution. Suppose that $\boldsymbol\delta^l \sim Dir(\boldsymbol\alpha^l)$ where $\boldsymbol\alpha^l=(\alpha^l_1,\ldots,\alpha^l_n)$ and  $\alpha^l_0=\sum_{i=1}^n \alpha^l_i$. Then 
\begin{equation}
E[\boldsymbol\delta^l]=\frac{\boldsymbol\alpha^l}{\alpha^l_0}
\label{eqn:delta_moment1}
\end{equation}
\begin{equation}
E[(\boldsymbol\delta^l-E[\boldsymbol\delta^l])(\boldsymbol\delta^l-E[\boldsymbol\delta^l])']=\frac{diag(\boldsymbol\alpha^l/\alpha^l_0)-(\boldsymbol\alpha^l/\alpha^l_0)(\boldsymbol\alpha^l/\alpha^l_0)'}{\alpha^l_0+1}
\label{eqn:delta_moment2}
\end{equation}
and for higher moments,
\begin{equation}
E\left[\prod_{i=1}^n(\boldsymbol\delta^l_i)^{\beta_i}\right]=\frac{\Gamma(\alpha^l_0)}{\Gamma(\alpha^l_0+\sum_{j=1}^n\beta_j)}\prod_{i=1}^n\frac{\Gamma(\alpha^l_i+\beta_i)}{\Gamma(\alpha^l_i)}
\label{eq:higher_moments}
\end{equation}
Using the formula for higher moments above, we get the following result on the third moment:
\begin{lemma}
For $\boldsymbol\delta^l \sim Dir(\boldsymbol\alpha^l)$ where $\boldsymbol\alpha^l=(\alpha^l_1,\ldots,\alpha^l_n)$ and $\alpha^l_0=\sum_{i=1}^n \alpha^l_i$, we have
\begin{align*}
&E[(\boldsymbol\delta^l-E[\boldsymbol\delta^l])_i(\boldsymbol\delta^l-E[\boldsymbol\delta^l])_j(\boldsymbol\delta^l-E[\boldsymbol\delta^l])_k]\nonumber \\
&=\frac{1}{(\alpha^l_0+1)(\alpha^l_0+2)}\cdot \begin{cases} 4(\frac{\alpha^l_i}{\alpha^l_0})^3-6(\frac{\alpha^l_i}{\alpha^l_0})^2+2(\frac{\alpha_i^l}{\alpha^l_0})~~~i=j=k\\
4(\frac{\alpha^l_i}{\alpha^l_0})^2(\frac{\alpha^l_k}{\alpha^l_0})-2(\frac{\alpha^l_i}{\alpha^l_0})(\frac{\alpha^l_k}{\alpha^l_0})~~~i=j\neq k\\
4(\frac{\alpha^l_i}{\alpha^l_0})(\frac{\alpha^l_j}{\alpha^l_0})(\frac{\alpha^l_k}{\alpha^l_0})~~~i\neq j\neq k
\end{cases}
\end{align*}
\label{lemma_third_order}
\end{lemma}

With these properties, we are now ready to construct two implementable estimators from Dirichlet mixtures.
\subsection{An Eligible Dirichlet Mixture}\label{sec:1estimator}
We use a mixture of $K=n(n-1)/2+1$ Dirichlet distributions. Consider
\begin{equation*}
\boldsymbol\delta^* =\sum_{l=1}^{n-1}\sum_{i=l+1}^{n}\theta^{l,i}\boldsymbol\delta^{l,i}+\theta^n\boldsymbol\delta^n
\end{equation*}
where $\boldsymbol\delta^{l,i}\sim Dir(\boldsymbol\alpha^{l,i})$ and $\boldsymbol\delta^{n}\sim Dir(\boldsymbol\alpha^{n})$. We set
\begin{align*}
&\theta^{1}=\theta^{1,i}=\frac{2p_1}{n-1},~i=2,3,\ldots,n\nonumber \\
&\theta^{2}=\theta^{2,i}=\frac{2p_2-\theta^{1}}{n-2},~i=3,4,\ldots,n \nonumber\\
&\ldots\nonumber\\
&\theta^{l}=\theta^{l,i}=\frac{2p_{l}-\sum_{j=1}^{l-1}\theta^{j}}{n-l},~i=l+1,\ldots,n\\
&\ldots\nonumber\\
&\theta^{n-1}=\theta^{n-1,n}=2p_{n-1}-\sum_{j=1}^{n-2}\theta^{j}\nonumber \\
&\theta^{n}=p_n-\frac{1}{2}\sum_{j=1}^{n-1}\theta^{j}
\end{align*}
and
\begin{align*}
&\frac{\boldsymbol\alpha^{l,i}}{\alpha^{l,i}_0}=\frac{1}{2}(\mathbf{e}_i+\mathbf{e}_l)\\
&\alpha^{l,i}_0=\mathbbm{1}'\boldsymbol\alpha^{l,i}=C(\theta^{l,i})^2-1\\
&\boldsymbol\alpha^{n}=\mathbf{e}_n
\end{align*}
where $C$ is a large constant such that $C(\theta^{l,i})^2>1$ for all $l=1,\ldots,n-1$ and $i=l+1,\ldots,n$. Notice that some entries of $\boldsymbol\alpha^{l,i}$ are 0, and by this we mean that the corresponding entries of $\boldsymbol \delta^{l,i}$ are set to 0 and the rest of the components follow a Dirichlet distribution. With this notion, the moment properties \eqref{eqn:delta_moment1}, \eqref{eqn:delta_moment2}, and \eqref{eq:higher_moments} still hold. 


We explain the intuition behind the above design. For each $\boldsymbol\delta^{l,i}\sim Dir(\boldsymbol\alpha^{l,i})$, we have its variance matrix
\begin{align}\label{eq:cov_mat}
    E[(\boldsymbol\delta^{l,i}- E[\boldsymbol\delta^{l,i}])(\boldsymbol\delta^{l,i}- E[\boldsymbol\delta^{l,i}])'] &= \frac{diag(\boldsymbol\alpha^{l,i}/\alpha^{l,i}_0)-(\boldsymbol\alpha^{l,i}/\alpha^{l,i}_0)(\boldsymbol\alpha^{l,i}/\alpha^{l,i}_0)'}{\alpha^{l,i}_0+1}\nonumber\\
    &=\frac{diag((\mathbf{e}_i+\mathbf{e}_l)/2)-(\mathbf{e}_i+\mathbf{e}_l)(\mathbf{e}_i+\mathbf{e}_l)'/4}{\alpha^{l,i}_0+1}\nonumber\\
    &\propto diag(\mathbf{e}_i+\mathbf{e}_l)-(\mathbf{e}_i+\mathbf{e}_l)(\mathbf{e}_i+\mathbf{e}_l)'/2
\end{align}

\begin{figure}[h!]
  \centering
  \begin{subfigure}[b]{0.144\linewidth}
    \includegraphics[width=\linewidth]{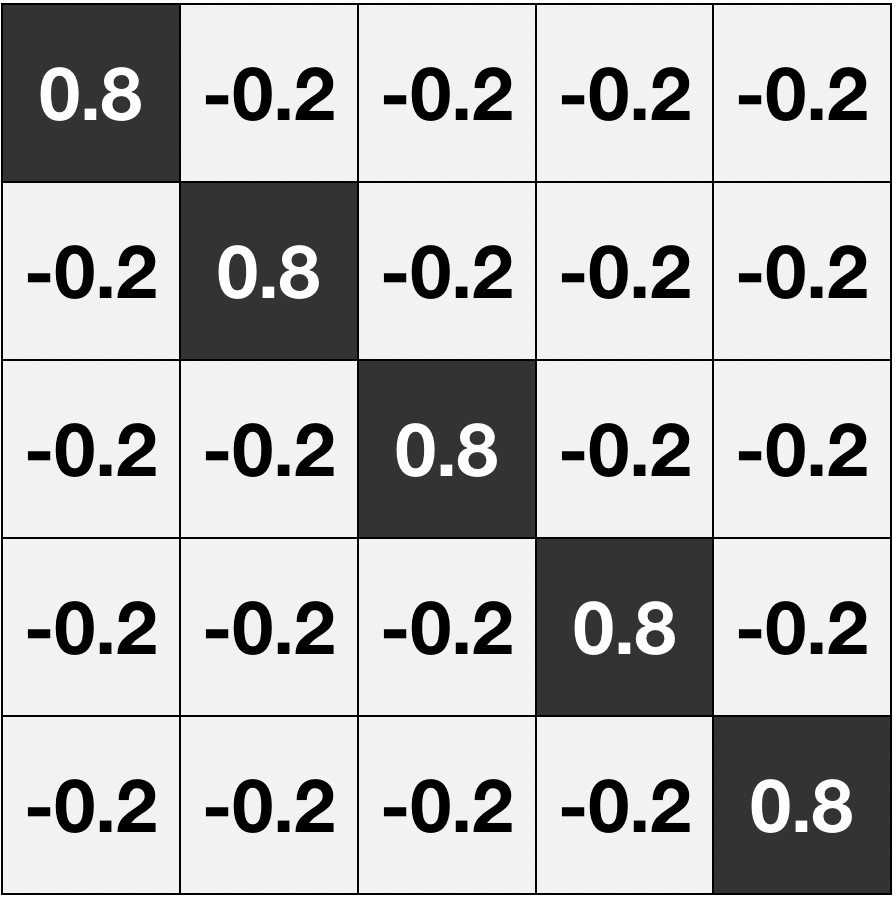}
    \caption{}
    \label{fig:cov_mat_all}
  \end{subfigure}
  ~~~~~  ~~~~~
  \begin{subfigure}[b]{0.6\linewidth}
    \includegraphics[width=0.24\linewidth]{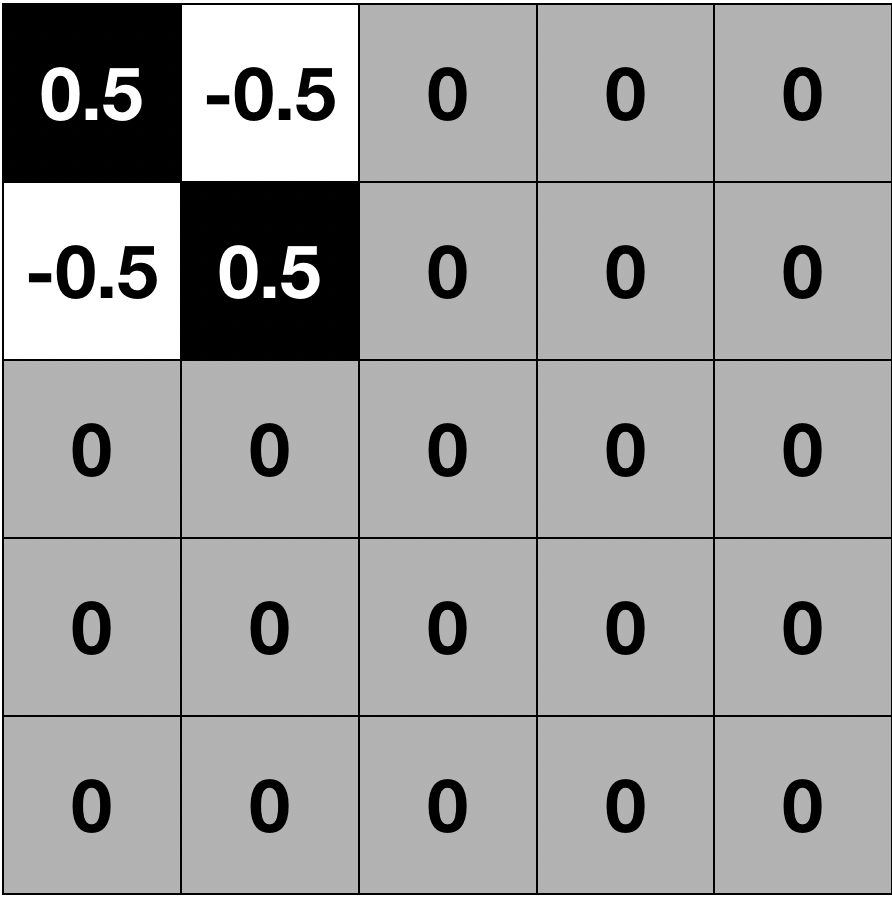}
    \includegraphics[width=0.24\linewidth]{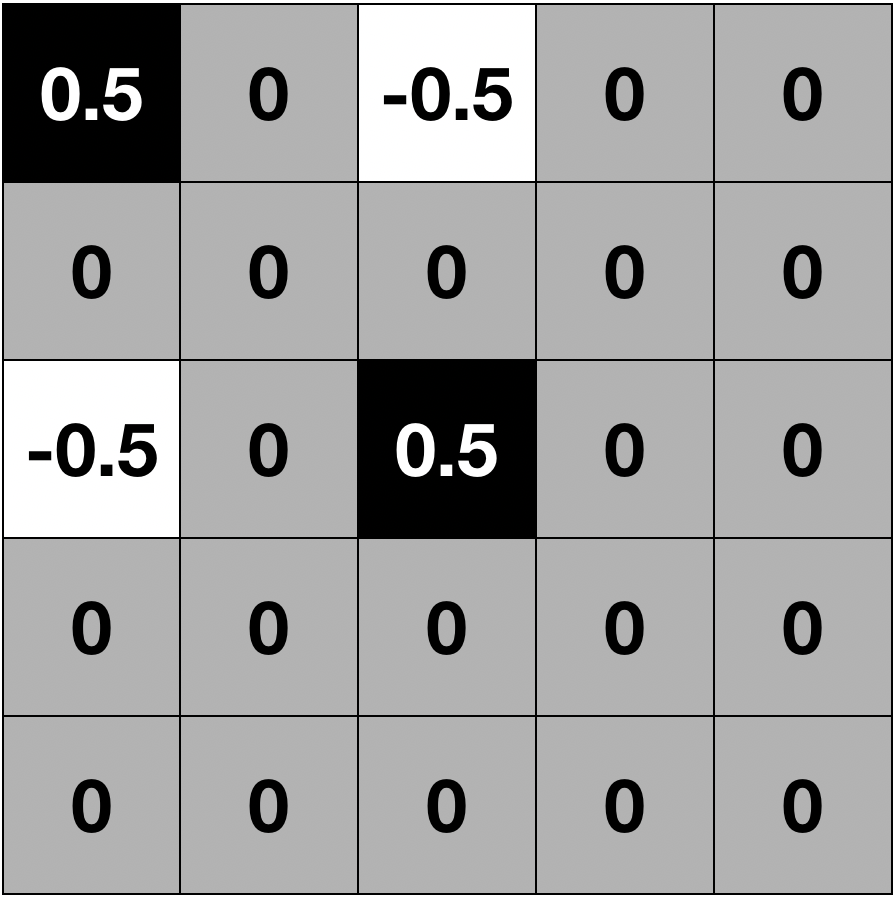}
    \includegraphics[width=0.24\linewidth]{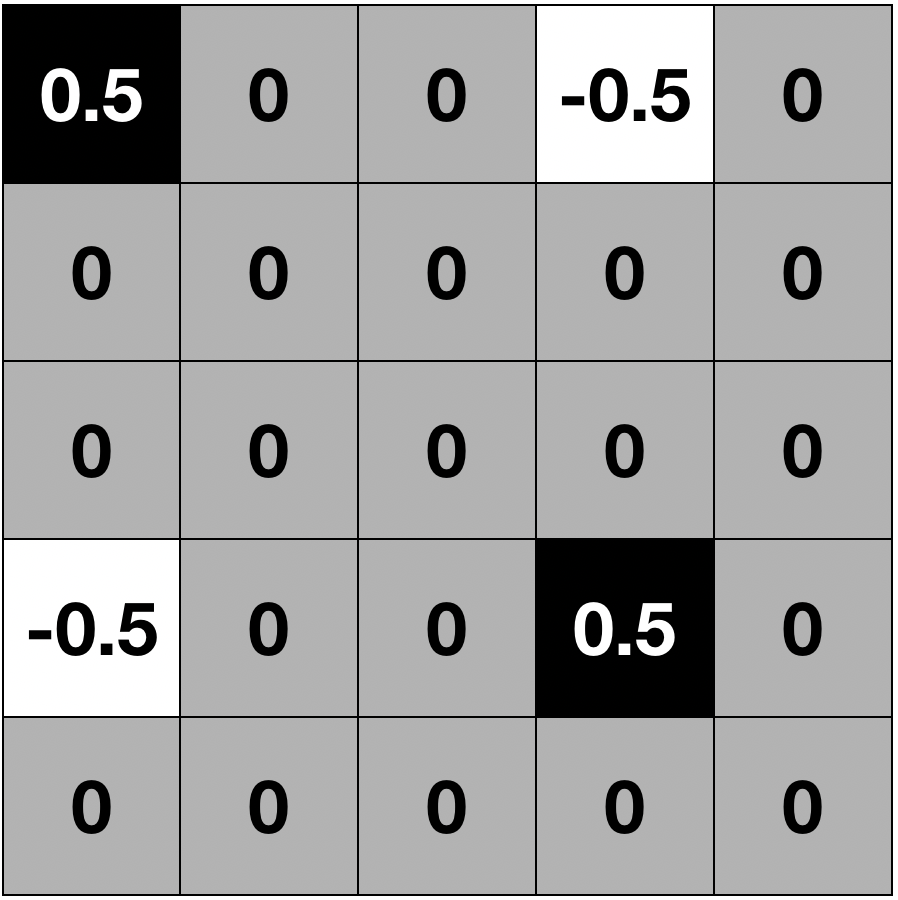}
    \includegraphics[width=0.24\linewidth]{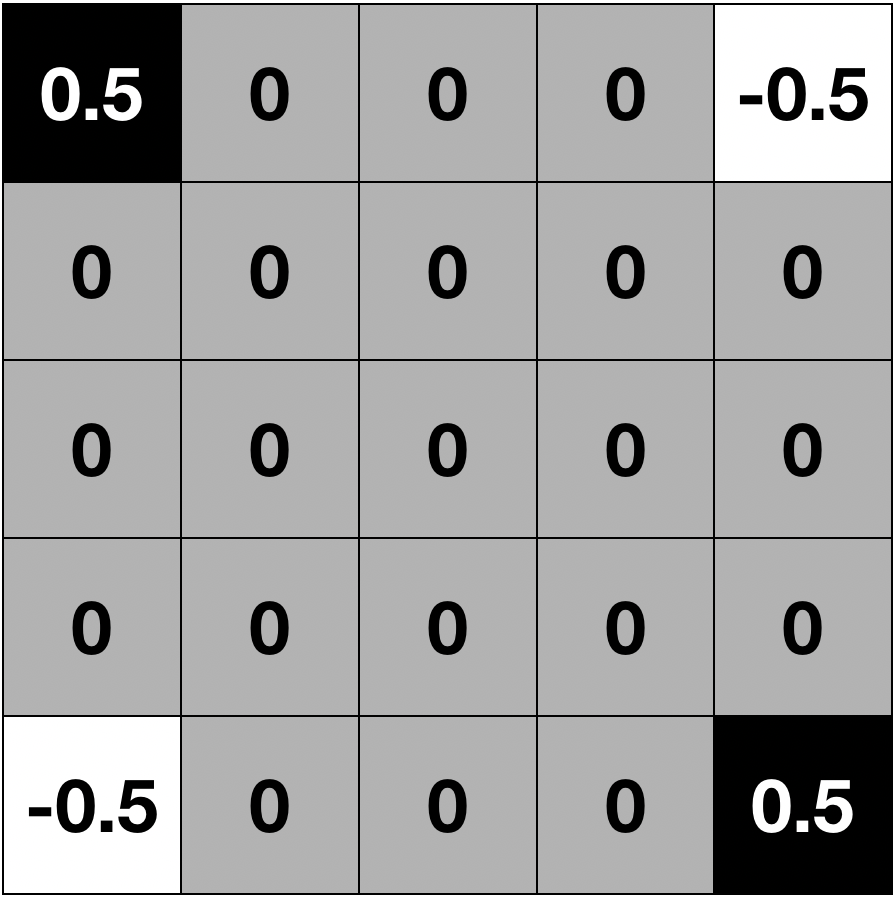}
    \caption{}
    \label{fig:cov_mat_each}
  \end{subfigure}
\caption{Variance matrices, $n=5$. \ref{fig:cov_mat_all}: The desired variance matrix for $\boldsymbol \delta^*$, $\mathbf I - \frac{\mathbbm 1\mathbbm 1'}{5}$. \ref{fig:cov_mat_each}: (Scaled) variance matrices of $\boldsymbol \delta^{1,2},\boldsymbol \delta^{1,3},\boldsymbol \delta^{1,4},\boldsymbol \delta^{1,5}$, and we choose the coefficients $\{\theta^{i,j}\}$ such that the linear combination of $\{Var[\boldsymbol \delta^{i,j}]\}$ gives $\mathbf I - \frac{\mathbbm 1\mathbbm 1'}{5}$.}
\label{fig:cov_mat}
\end{figure}

The variance of each Dirichlet random vector $\boldsymbol\delta^{l,i}$ as stated in \eqref{eq:cov_mat} is of a very simple form, where there are only 4 nonzero entries (see Figure \ref{fig:cov_mat_each} for the case when $n=5$). Since $\boldsymbol \delta^*$ is a linear combination of $\{\boldsymbol\delta^{l,i},\boldsymbol\delta^n\}$ which are independent, its variance is also a linear combination of the variances of $\{\boldsymbol\delta^{l,i},\boldsymbol\delta^n\}$. And in our design, we choose the coefficients $\{\theta^i\}$ and the $\{\alpha^{l,i}_0,\alpha^n_0\}$ in a way such that the linear combination of these simple matrices in Figure \ref{fig:cov_mat_each} gives $\mathbf I - \frac{\mathbbm 1\mathbbm 1'}{n}$ as in Figure \ref{fig:cov_mat_all}.

To proceed, we have the following monotonicity property on $\theta^i$'s.

\begin{lemma}
For the $\theta^i$ sequence depicted above, 
\begin{equation*}
\theta^1\leq \theta^2\leq \ldots \leq \theta^{n-1}
\end{equation*}
\label{lemma:theta_increasing}
\end{lemma}
Recall that we need to choose $C(\theta^{l})^2>1$ for all $l=1,\ldots,n-1$, and this lemma shows that $\theta^1 = \min_{1\leq j\leq n-1}\theta^j$. Thus, any $C>1/(\theta^1)^2=(n-1)^2/(4p_1^2)$ is a legitimate choice to attain the designated values of $\theta^i$ and $\boldsymbol\alpha^{l,i}$, $ \alpha_0^{l,i}$, $\boldsymbol\alpha^n$ above. Next, we show that using $\boldsymbol\delta^*$ indeed satisfies all three conditions \eqref{first order}, \eqref{second order} and \eqref{third order}.

\begin{theorem}[Eligible Dirichlet mixture]
The Dirichlet mixture $\boldsymbol\delta^*$ with $\gamma =\frac{4C}{n}$ where $C> \frac{(n-1)^2}{4p_1^2}$ satisfies \eqref{first order}, \eqref{second order}, and \eqref{third order} with $\mu =0$.
\label{thm:delta1}
\end{theorem}

\subsection{Another Eligible Dirichlet Mixture}\label{sec:2estimator}
Notice that Theorems \ref{thm:SFE}, \ref{thm:FFE} and \ref{thm:CFE} give that the variances of the gradient estimators are of order $O(\frac{\gamma n}{Rc^2})$ for all of $\hat\psi_{SFE}$, $\hat\psi_{FFE}$ and $\hat\psi_{CFE}$. The constant $\gamma$ can scale with the dimension $n$. Because of this, it is at times better to reduce $\gamma$, and hence the variance, at the cost of a $O(c)$ bias. It turns out that we can achieve so with another Dirichlet mixture, with the property that this $O(c)$ term vanishes when the problem dimension $n$ increases, and thus giving a beneficial balance of bias and variance in high dimension. Consider $\boldsymbol\delta^{**}=\sum_{l=1}^n\theta^l\boldsymbol\delta ^l,~\boldsymbol\delta ^l\sim Dir(\alpha^l)$ where 
\begin{equation*}
\theta^l=\begin{cases}
np_1~~~\text{\ for\ }l=1\\
p_l-p_1~~~\text{\ for\ }l=2,3,\ldots,n
\end{cases}
\end{equation*}
and
\begin{equation*}
\alpha^l=\begin{cases}
n^{\eta}\mathbbm{1}~~~\text{\ for\ }l=1\\
\mathbf{e}_l~~~\text{\ for\ }l=2,3,...,n
\end{cases}
\end{equation*}
for some $\eta\in \mathbbm R$ and $\gamma = \frac{n^{\eta+1}+1}{np_1^2}$. We have the following result for $\boldsymbol\delta^{**}$.
\begin{theorem}[Eligible Dirichlet mixture]
The Dirichlet mixture $\boldsymbol\delta^{**}$ satisfies \eqref{first order} and \eqref{second order}. In addition, $\exists \epsilon_{0} \in \mathbbm R$ such that $ \lVert \frac{c}{2} E[(\boldsymbol\delta - \mathbf p)'\nabla ^2 Z(\mathbf p)  (\boldsymbol\delta - \mathbf p)S(\mathbf p,\boldsymbol\delta)] -\epsilon_{0} \mathbbm 1\rVert_{\infty}=O(c\lVert \nabla ^2 Z \rVert_{\infty} n^{-\eta-1})$. 
\label{thm:delta2}
\end{theorem}

Here, the $\lVert \cdot \rVert_{\infty}$ is the entrywise maximum of a vector. The rationale behind this design is that the variance of $\boldsymbol \delta^l$ for $l=2,3,\ldots,n$ is $\mathbf 0$, while the variance of $\boldsymbol\delta^1\propto diag(\mathbbm 1) - \mathbbm 1\mathbbm 1'/n = \mathbf I - \mathbbm 1\mathbbm 1'/n$. The coefficients $\{\theta^{l}\}$ are chosen such that the condition \eqref{first order} is satisfied, and $\gamma$ is chosen such that \eqref{second order} is satisfied. By tuning $\eta$, we can control the decay of the $O(c)$ in the bias. 

Theorem \ref{thm:delta2} is useful in two senses. First, it shows that $\boldsymbol\delta^{**}$ satisfies \eqref{first order} and \eqref{second order}, so
that SFE, FFE, CFE with $\boldsymbol\delta^{**}$ has bias $E[\psi_{SFE}] = E[\psi_{FFE}] = \nabla Z(\mathbf p) -\frac{\mathbbm 1'\nabla Z(\mathbf p)}{n}\mathbbm 1+ \frac{c }{2}E[(\boldsymbol\delta - \mathbf p)'\nabla^2Z(\mathbf p)(\boldsymbol\delta - \mathbf p)S(\mathbf p,\boldsymbol\delta)]+O(nc^2L_2)$, i.e. the bias is $O(c)$ (Theorem \ref{thm:SFE} and \ref{thm:FFE}). 

Second, it provides a bound on the bias due to the Hessian term (the third term in \eqref{interim}, i.e. $\frac{c }{2}E[(\boldsymbol\delta - \mathbf p)'\nabla^2Z(\mathbf p)(\boldsymbol\delta - \mathbf p)S(\mathbf p,\boldsymbol\delta)]$). To be precise, we have the following two direct consequences: 
\begin{corollary}
Assuming  $\lVert \nabla^2 Z\rVert_{\infty}$ is bounded, when $\eta>-1$, as $n\to \infty$, $\exists \epsilon_{0,n}\in \mathbbm R$ such that $\lVert \frac{c}{2}E[(\boldsymbol\delta - \mathbf p)'\nabla^2Z(\mathbf p)(\boldsymbol\delta - \mathbf p)S(\mathbf p,\boldsymbol\delta)]-\epsilon_{0,n}\mathbbm 1 \rVert_{\infty} = O(c\lVert \nabla Z\rVert_{\infty} n^{-\eta-1})\to 0$. That is, for fixed $c$, the entrywise bias due to the Hessian approaches 0 as the dimension of the problem approaches infinity.
\label{cor:delta2_to_zero}
\end{corollary}
\begin{corollary}
Assuming  $\lVert \nabla^2 Z\rVert_{\infty}$ is bounded, taking $\eta=-1$, $\exists \epsilon_{0,n}\in \mathbbm R$, such that $\lVert \frac{c}{2}E[(\boldsymbol\delta - \mathbf p)'\nabla ^2 Z(\mathbf p)  (\boldsymbol\delta - \mathbf p)S(\mathbf p,\boldsymbol\delta)] -\epsilon_{0,n} \mathbbm 1\rVert = O(cM_2\sqrt{n})$. 
\label{cor:delta2_comparable}
\end{corollary}

Next, we will provide comparisons for $\eta$, using Finite Difference (FD) methods as baseline. Using the bound between the $L_2$ norm and the $L_{\infty}$ that $\lVert \mathbf x \rVert\leq \sqrt{n}\lVert \mathbf x \rVert_{\infty}$ for all $\mathbf x\in \mathbbm R^n$, Theorem \ref{thm:delta2} gives us $\lVert \frac{c}{2}E[(\boldsymbol\delta - \mathbf p)'\nabla^2Z(\mathbf p)(\boldsymbol\delta - \mathbf p)S(\mathbf p,\boldsymbol\delta)]-\epsilon_{0,n}\mathbbm 1 \rVert= O(c\lVert \nabla Z\rVert_{\infty} n^{-\eta-1/2})$. Thus, regarding the bias' dependency on the problem dimension, for $\eta = -1$, the bias when using $\boldsymbol \delta^{**}$ is $O(\sqrt{n})$ (the same as FD from Theorems \ref{thm:FD} and \ref{thm:FDrd} in the sequel). Also, $\gamma = O(n)$ for $\boldsymbol \delta^{**}$, and so the variances are $O(n^2)$ for both FD and for $\boldsymbol \delta^{**}$. On the other hand, if we choose $\eta>-1$, the bias when using $\boldsymbol \delta^{**}$ is $O(n^{-\eta-1/2})$ which is better than the $O(\sqrt{n})$ bias of FD. However, $\gamma = O(n^{\eta+2})$, making the variance using $\boldsymbol \delta^{**}$ $O(n^{\eta+3})$ worse than the variance for FD. 

In short, when we choose $\eta=-1$, the bias and variance using $\boldsymbol \delta^{**}$ are the same as using FD. If we choose $\eta>-1$, we have better bias but worse variance. In the case of small or moderate simulation budget $R$ (e.g., when we run gradient-descent optimization), so that the variance is the dominant term in the mean squared error, it appears choosing $\eta = -1$ is appropriate.

\section{Comparisons among Gradient Estimators and with Standard Finite Differences}\label{sec:comparisons}
We statistically compare among our constructed gradient estimators and also with naive FD schemes, based on the results in Sections \ref{sec:analysis} and \ref{sec:Dir}. For FD, we consider standard coordinate-wise FD (FD$_{standard}$) that divides the budget equally into each dimension, and its randomized variant FD$_{random}$ where each repetition uniformly randomly selects a dimension and applies FD on it. We will focus on forward FD, since a single function evaluation in this case would not be able to cancel out the first term in Taylor series and obtain a consistent gradient estimator, while central FD may contain evaluation input that shoots outside the probability simplex as we have discussed earlier. 

To compare with FD, we present the biases and variances of FD$_{standard}$ and FD$_{random}$. For FD$_{standard}$, suppose the total number of simulation replications $R$ is divided into the $n$ dimensions given by $R = nR_p$. We define 
\begin{equation*}
(\hat \psi_{FD,standard})_i := \frac{1}{R_p}\sum_{j=1}^{R_p}\frac{\hat{Z}^{2j-1}(\mathbf p +c(\mathbf e_i - \mathbf p))-\hat{Z}^{2j}(\mathbf p)}{c}
\end{equation*}
Denote $\psi_{FD,standard}:= E[\hat \psi_{FD,standard}]$. We have the following:
\begin{theorem}[Bias and variance for standard finite difference]
Under Assumptions \ref{assumption:smooth}  and \ref{assumption:smooth_noise} and setting $R=nR_p$, $\psi_{FD,standard}$ has bias given by
\begin{equation*}
    \psi_{FD,standard}-\nabla Z(\mathbf p)-\epsilon_0 \mathbbm 1=\epsilon
\end{equation*}
where $\epsilon_0 = -\nabla Z(\mathbf p)'\mathbf p$ and $\lVert \epsilon \rVert \leq 4cM_2\sqrt{n}$, and variance given by
\begin{equation*}
    E[\lVert \hat \psi_{FD,standard} - \psi_{FD,standard}\rVert^2] = \frac{2n^2\sigma^2}{Rc^2}
\end{equation*}
\label{thm:FD}
\end{theorem}

Different from FD$_{standard}$ where each direction is perturbed exactly $R_p$ times and in total there are $nR_p=R$ function evaluations, for FD$_{random}$, we assume that for $i=1,\ldots,R$, $l_i$ is a random ``index'' variable that takes value from $[n]$ with equal probability, and $l_i$'s are independent. Define
\begin{equation*}
\hat \psi_{FD,random}^i := n\frac{\hat{Z}(\mathbf p +c(\mathbf e_{l_i} - \mathbf p))-\hat{Z}(\mathbf p)}{c}\mathbf e_{l_i}
\end{equation*}
and
\begin{equation*}
\hat \psi_{FD,random} := \frac{1}{R}\sum_{i=1}^{R} \hat \psi_{FD,random}^i
\end{equation*}
Denote $\psi_{FD,random}:=E[\hat \psi_{FD,random}]$. We have the following:
\begin{theorem}[Bias and variance for random-dimension finite difference]
Under Assumptions \ref{assumption:smooth} and \ref{assumption:smooth_noise}, $\hat \psi_{FD,random}$ has bias given by
\begin{equation*}
    \psi_{FD,random}-\nabla Z (\mathbf p)-\epsilon_0 \mathbbm 1=\epsilon
\end{equation*}
where $\epsilon_0 = - \nabla Z(\mathbf p)'\mathbf p$ and $\lVert \epsilon\rVert \leq 4cM_2\sqrt{n} $, and variance given by
\begin{equation*}
    E[\lVert \hat \psi_{FD,random}-\psi_{FD,random} \rVert ^2] \leq \frac{2n^2}{Rc^2}[\sigma^2+4c^2M_1^2+16c^4L_1^2]
\end{equation*}
\label{thm:FDrd}
\end{theorem}

Tables \ref{tb:bias} and \ref{tb:var} show the biases and variances per simulation run or pair of simulation runs for SFE, FFE and CFE, for both simultaneous-perturbation-type estimators satisfying combinations of \eqref{first order}, \eqref{second order} and \eqref{third order}, and also FD$_{standard}$ and FD$_{random}$. For variance, the satisfaction of \eqref{third order} does not play a role in the bound, and thus we skip the case  \eqref{first order} +  \eqref{second order} + \eqref{third order} in Table \ref{tb:var}. Moreover, in the case that $\frac{1}{np_1} = O(1)$, that is, $np_1$ is bounded from below, or $\mathbf p$ is bounded away from the boundary of the simplex, from Theorem \ref{thm:delta1} we can set $C = O(\frac{n^2}{p_1^2})$ for $\boldsymbol\delta^{*}$, and thus $\gamma = \frac{4C}{n}=O(n^3)$. For $\boldsymbol\delta^{**}$, Theorem \ref{thm:delta2} shows that we can choose $\eta=-1$, and thus $\gamma = O(n)$. These give the second and third rows of Table \ref{tb:var}. We omit the $O(\cdot)$ in both tables. In addition, $R$ for FD$_{standard}$ is the total number of function evaluations, which is divided into $R/n$ evaluations for each dimension (assuming divisibility).

\begin{table}[h!]
  \begin{center}
    \caption{Bias $\lVert E[\hat {\psi}]-\nabla Z -\epsilon_0 \mathbbm 1\rVert$ for some $\epsilon_0\in\mathbbm R$ for $\boldsymbol \delta$ that satisfies \eqref{first order} and \eqref{second order} in general and its special case $\boldsymbol \delta^{**}$ with $\eta = -1$, and that satisfies \eqref{first order}, \eqref{second order}, and \eqref{third order} in general and its special case $\boldsymbol \delta^{*}$.}
    \begin{tabular}{c||c|c|c|c|c}
      &SFE&FFE&CFE&FD$_{standard}$&FD$_{random}$\\
      \hline
 \eqref{first order} + \eqref{second order}& $ncM_2+nc^2L_2$ &$ncM_2+nc^2L_2$ &$nc^2L_2$&\multirow{4}{*}{$\sqrt{n}cM_2$}&\multirow{4}{*}{$\sqrt{n}cM_2$}\\
$\boldsymbol\delta^{**}, \eta = -1$& $\sqrt{n}cM_2+nc^2L_2$ &$\sqrt{n}cM_2+nc^2L_2$ &$nc^2L_2$&&\\\cline{1-4}
\eqref{first order} + \eqref{second order}+\eqref{third order}&$nc^2L_2$ &$nc^2L_2$&$nc^2L_2$&&\\
$\boldsymbol\delta^{*}$&$nc^2L_2$ &$nc^2L_2$&$nc^2L_2$&&
    \end{tabular}
    \label{tb:bias}
  \end{center}
\end{table}

\begin{table}[h!]
  \begin{center}
    \caption{Variance $E[\lVert \hat \psi - \psi\rVert^2]$ for $\boldsymbol \delta$ that satisfies \eqref{first order} and \eqref{second order} in general, for $\boldsymbol \delta^{**}$ with $\eta = -1$, and for $\boldsymbol \delta^*$, assuming $\frac{1}{np_1} = O(1)$.}
    \begin{tabular}{c||c|c|c|c|c}
      &SFE&FFE&CFE&FD$_{standard}$&FD$_{random}$\\
      \hline
\eqref{first order} + \eqref{second order}&$\frac{\gamma n(\sigma^2+M_0^2)}{Rc^2}$ & $\frac{\gamma n(\sigma^2+c^2M_1^2+c^4L_1^2)}{Rc^2}$& $\frac{\gamma n(\sigma^2+c^2M_1^2+c^4L_1^2)}{Rc^2}$ &\multirow{3}{*}{$\frac{n^2\sigma^2}{Rc^2}$}&\multirow{3}{*}{$\frac{ n^2(\sigma^2+c^2M_1^2+c^4L_1^2)}{Rc^2}$}\\
$\boldsymbol \delta^{**}, \eta = -1, \gamma = O(n)$& $\frac{ n^2(\sigma^2+M_0^2)}{Rc^2}$ & $\frac{ n^2(\sigma^2+c^2M_1^2+c^4L_1^2)}{Rc^2}$& $\frac{ n^2(\sigma^2+c^2M_1^2+c^4L_1^2)}{Rc^2}$ &&\\
$\boldsymbol \delta^{*}, \gamma = O(n^3)$& $\frac{ n^4(\sigma^2+M_0^2)}{Rc^2}$ & $\frac{ n^4(\sigma^2+c^2M_1^2+c^4L_1^2)}{Rc^2}$& $\frac{ n^4(\sigma^2+c^2M_1^2+c^4L_1^2)}{Rc^2}$ &&

    \end{tabular}
    \label{tb:var}
  \end{center}
\end{table}

From Table \ref{tb:var} we see that the variances of SFE, SFE, and CFE are all linear in $\gamma$ for general $\boldsymbol \delta$ that satisfies \eqref{first order} and \eqref{second order}. An immediate question is whether we can choose $\gamma$ arbitrarily small to reduce the variance. The answer is negative, as seen from the following lemma:
\begin{lemma}\label{lm:gamma}
To satisfy \eqref{second order}, $\gamma\geq \frac{n-1}{4}$. 
\end{lemma}

Thus, from Lemma \ref{lm:gamma}, we see that $\gamma$ is at least linear in the dimension $n$, and with such $\gamma$, we have $O(n^2)$ variance for SFE/FFE/CFE, the same as $\boldsymbol\delta^{**}$ with $\eta = -1$, and also as FD$_{standard}$ and FD$_{random}$.



\subsection{Comparisons of Bias and Variance Among $\hat \psi_{SFE},\hat \psi_{FFE},\hat \psi_{CFE},\hat \psi_{FD}$}
Now we compare the biases and variances among all estimators. The biases in Table \ref{tb:bias} are all $O(c)$ or $O(c^2)$, while variances in Table \ref{tb:var} are all $O(c^{-2})$. The main differences lie in the dependence on dimension $n$. Using a small $c$ such that $\sqrt nc <<1$, we see that SFE/FFE/CFE that satisfy \eqref{first order}, \eqref{second order} and \eqref{third order}, or using $\boldsymbol \delta^*$, have the smallest bias of order $nc^2$. However, using $\boldsymbol\delta^*$ has a larger variance, of order $O(n^4c^{-2})$, than all the other estimators with order $n^2c^{-2}$ (assuming $\gamma$ is chosen as $O(n)$ for general estimators in view of Lemma \ref{lm:gamma}). On the other hand, the $O(nc)$ bias of SFE/FFE that satisfy \eqref{first order} and \eqref{second order} only is larger than all other estimators. Their variances are comparable to FD$_{standard}$ and FD$_{random}$ when $\gamma$ is set to be $O(n)$, if we assume  $\sqrt{c^2M_1^2+c^4L_1^2} << \sigma $ and $ M_0<<\sigma$. Note that CFE has consistently lower bias of order $O(nc^2)$, and its variance is similar to SFE, FFE, FD$_{standard}$ and FD$_{random}$. However, as discussed before, CFE may not be implementable in practice. SFE and FFE are comparable in both bias and variance in general, but as noted before, FFE has a lower variance than SFE because of the canceling of the first term in the Taylor series, an advantage that is not reflected in Table \ref{tb:var}.


Focusing on the explicit schemes $\boldsymbol\delta^*$ and $\boldsymbol\delta^{**}$ on the implementable estimators SFE and FFE, we see that both are comparable to FD$_{standard}$ and FD$_{random}$ in that none of them strictly dominate others.  $\boldsymbol \delta^*$ has a smaller bias but larger variance than FD$_{standard}$ and FD$_{random}$, while $\boldsymbol\delta^{**}$ has similar bias and variance as FD$_{standard}$ and FD$_{random}$. Given that the simulation budget for gradient estimators when used in descent-type optimization algorithms (such as the Frank-Wolfe and mirror descent procedures to be presented in Section \ref{sec:SA}) is typically small or moderate, we suggest to use $\boldsymbol \delta^{**}$ in practice. Also, although $\boldsymbol \delta^{**}$ has a comparable sampling complexity as FD$_{standard}$ and FD$_{random}$, empirically they perform differently when used in the descent procedures. We will discuss how and in what experimental conditions $\boldsymbol \delta^{**}$ appears to outperform simple FD schemes in Sections \ref{sec:opt num} and \ref{sec:summary}. Moreover, while not reflected in the developed theorems, simultaneous-perturbation-type schemes such as $\boldsymbol\delta^*$ and $\boldsymbol\delta^{**}$ have the obvious strength that they give at least some information about all components of the gradient with one or very few simulation runs.



\subsection{Comparisons with Gradient Estimators for Unconstrained Objectives}

Our results resemble, though also bear differences with, gradient estimation for unconstrained functions (i.e., the domain of the function $Z(\cdot)$ is not limited to the probability simplex) obtained in \cite{berahas2019theoretical}. There the authors consider gradient estimators of the form:
\begin{equation}\label{eq:unconstraind_1}
\hat \psi(\mathbf p) = \frac{1}{R}\sum_{j=1}^R\frac{\hat Z^{2j-1}(\mathbf p+c\boldsymbol \delta^j)-\hat Z^{2j}(\mathbf p)}{c}S(\mathbf p,\boldsymbol \delta^j)
\end{equation}
or the central version
\begin{equation}\label{eq:unconstraind_2}
\hat \psi(\mathbf p) = \frac{1}{R}\sum_{j=1}^R\frac{\hat Z^{2j-1}(\mathbf p+c\boldsymbol \delta^j)-\hat Z^{2j}(\mathbf p-c\boldsymbol \delta^j)}{2c}S(\mathbf p,\boldsymbol \delta^j)
\end{equation}
They consider two types of gradient estimators: interpolation type and smoothing type. For interpolation type estimators, they analyze Forward Finite Difference (FD$_{standard}^{uc}$) which resembles our FD$_{standard}$ but for unconstrained objective, and Linear Interpolation (LI). For smoothing type, they consider Gaussian Smoothing (GS), Central Gaussian Smoothing (CGS), Sphere Smoothing (SS), and Central Sphere Smoothing (CSS).

With Assumption \ref{assumption:smooth_noise}, results of \cite{berahas2019theoretical} can be adapted to show the bias and variance bounds in Table \ref{tb:unconstrained}. For FD$_{standard}^{uc}$ and LI, we perturb $R=R_pn$ times, i.e. $R_p$ times in each direction. For LI, $Q$ is a nonsingular $n\times n$ matrix such that $\lVert Q_{j,:}\rVert \leq 1$ and we use $A_{i,:},A_{:,i}$ to denote the $i$-th row and column of matrix $A$ respectively. (Results for FD$_{standard}^{uc}$ are from Theorem 2.1, LI from Theorem 2.3, GS from Equation (2.10) and Lemma 2.4, CGS from Equation (2.11) and Lemma 2.4, SS from Equation (2.35) and Lemma 2.9, and CSS from Equation (2.36) and Lemma 2.9 from \cite{berahas2019theoretical}). 

\begin{table}[h!]
  \begin{center}
    \caption{Bias $\lVert E[\hat \psi(\mathbf p) ]- \nabla Z(\mathbf p)\rVert $ and variance $E[\lVert \hat \psi(\mathbf p) - E[\hat \psi(\mathbf p)]\rVert^2 ]$ for unconstrained objectives, $O(\cdot)$ omitted}
    \begin{tabular}{c||c|c|c|c}
      Method&$\boldsymbol \delta^j$&$S(\mathbf p,\boldsymbol \delta^j)$&bias&variance\\
      \hline
    FD$_{standard}^{uc}$ \eqref{eq:unconstraind_1}& $\boldsymbol \delta^j = \mathbf e_l$, $j=(l-1)R_p+1,\ldots,lR_p$& $ n\boldsymbol \delta^j$&$\sqrt{n}cL_1$&$\frac{n^2\sigma^2}{c^2R}$\\
    \hline
      LI \eqref{eq:unconstraind_1}&$\boldsymbol \delta^j = Q_{l,:}$, $j=(l-1)R_p+1,\ldots,lR_p$&$nQ^{-1}_{:,l}$&$\lVert Q^{-1}\rVert \sqrt{n}cL_1$&$\frac{\lVert Q^{-1}\rVert^2n^2\sigma^2}{c^2R}$\\
      \hline
      GS \eqref{eq:unconstraind_1}&\multirow{2}{*}{$\boldsymbol \delta^j\sim \mathcal N(\mathbf 0,\mathbf I/n)$}&\multirow{2}{*}{$n\boldsymbol\delta^j$}&$cL_1$&$\frac{n(n\sigma^2+c^2M_1^2+c^4L_1^2n)}{c^2R}$\\
      CGS \eqref{eq:unconstraind_2}&&&$c^2L_2$&$\frac{n(n\sigma^2+c^2M_1^2+c^6L_2^2n)}{c^2R}$\\
      \hline
      SS \eqref{eq:unconstraind_1}&\multirow{2}{*}{$\boldsymbol \delta^j\sim Uni(\mathbbm S^{n-1})$}&\multirow{2}{*}{$n\boldsymbol\delta^j$}&$cL_1$&$\frac{n(n\sigma^2+c^2M_1^2+c^4L_1^2n)}{c^2R}$\\
      CSS \eqref{eq:unconstraind_2}&&&$c^2L_2$&$\frac{n(n\sigma^2+c^2M_1^2+c^6L_2^2n)}{c^2R}$\\
   
    \end{tabular}
    \label{tb:unconstrained}
  \end{center}
\end{table}

Notice that in Table \ref{tb:unconstrained}, for GS and CGS, we rescale the perturbation variance from $\mathbf I$ in \cite{berahas2019theoretical} to $\mathbf I/n$. This rescaling makes the comparison in Table \ref{tb:unconstrained} fairer, in the sense that now the ``lengths'' of the perturbation vectors $c\boldsymbol \delta^j$ are on the same scale. To be precise, for FD$_{standard}^{uc}$, LI (with orthonormal $Q$), SS, and CSS, the perturbation vector has length $\lVert c\boldsymbol \delta^j\rVert = c$, while for GS and CGS, $E[\lVert c\boldsymbol \delta^j\rVert^2] = c^2$. This is important especially when we have limited room for perturbation (e.g. $E[\lVert c\boldsymbol \delta^j\rVert^2] \leq \xi$ is bounded for some $\xi$). Now, with the same choice of $c$, we see that the bias for FD$_{standard}^{uc}$, LI (with orthonormal $Q$) are $O(\sqrt{n}cL_1)$, while that for GS, CGS, SS, and CSS are $O(c L_1)$. Thus, in terms of bias, GS, CGS, SS, and CSS are better than FD$_{standard}^{uc}$ and LI.

We provide more details on the rationale behind the rescaling mentioned above. The underlying reasoning is that gradient estimation might become harder (e.g., larger variance) as the allowed perturbation ``region'' becomes smaller. First, it is clear that the variance of the gradient estimator depends on the variance of the score function. The core idea behind LI, GS, SS, as well as our \eqref{second order} is that $\boldsymbol \delta^j S(\mathbf p,\boldsymbol \delta^j)'$ (or $(\boldsymbol \delta^j-\mathbf p) S(\mathbf p,\boldsymbol \delta^j)'$ for the simplex-constrained case), when taken average of, has an expectation of $\mathbf I$ (or $\mathbf I-\mathbbm 1\mathbbm 1'/n$). If we require $E[\lVert \boldsymbol \delta^j\rVert^2] $ or $E[\lVert \boldsymbol \delta^j-\mathbf p\rVert^2]$ to be small, then in order to achieve the desired sum or expectation, we need to compensate for this by making $\lVert S(\mathbf p,\boldsymbol\delta^j)\rVert $ larger. Notice that the observation noise is multiplied by this score function directly, which means that a larger score function will magnify the observation noise more. As an example, we see from LI that if we shrink the matrix $Q$, i.e., $ \boldsymbol\delta^j = s Q_{j,:}$ for some $s\in (0,1]$, we have variance $O(\frac{1}{s^2}\cdot \frac{\lVert Q^{-1}\rVert^2n^2\sigma^2}{c^2R})$. Thus, the smaller $s$ is, the smaller the perturbation is, but the larger the variance will be. For our simplex-constrained case, Lemma \ref{lm:gamma} tells us $\gamma\geq \frac{n-1}{4}$, i.e., $E[\lVert S(\boldsymbol \delta,\mathbf p)\rVert^2] = \gamma (n-1)= \Theta(n^2)$. In contrast, without the simplex constraint, methods such as GS with perturbation $\boldsymbol \delta^j$ and score function $S(\boldsymbol \delta^j,\mathbf p) = \boldsymbol \delta^j$ where $\boldsymbol \delta^j\sim \mathcal N(\mathbf 0,\mathbf I)$ gives that $E[\lVert  S(\boldsymbol \delta^j,\mathbf p) \rVert^2] = n$. Thus, the observation noise will be enlarged more by the score function of the gradient estimator for the simplex constrained objective than by that for the GS (with the same $c$).

In the experiments in \cite{berahas2019theoretical}, LI and FD$_{standard}^{uc}$ appear more effective. It is explained that when there is no noise in the observation, i.e., $\sigma = 0$, the variances of interpolation-type estimators (FD$_{standard}^{uc}$, LI) are 0, and their biases go to $0$ as $c\to 0$. In contrast, smoothing-type estimators (GS, CGS, SS, CSS) have biases going to $0$ as $c\to 0$, but their variances do not go to $0$ due to the randomness in the perturbation direction $\boldsymbol\delta^j$. It can also be seen from Table \ref{tb:unconstrained} that their variances have an $O(\frac{nM_1^2}{R})$ component that does not depend on $c$. 

Back to the case of simplex-constrained function, our results in Tables \ref{tb:bias} and \ref{tb:var} have some similarities with Table \ref{tb:unconstrained}. For the bias, FD-type estimators (and LI) are the same between constrained and unconstrained cases. However, our simultaneous-perturbation-type estimators with \eqref{first order} and \eqref{second order} have biases $O(ncM_2+nc^2L_2)$, which are worse than the counterparts in the unconstrained case such as GS and SS which have biases $O(cL_1)$ in terms of dimension dependency. Nevertheless, we believe that such good bias behaviors for GS and SS are only valid for special types of distribution of $\boldsymbol \delta^j$, say, Gaussian as in GS, and uniform on the sphere as in SS. For a generic simultaneous-perturbation-type gradient estimator for an unconstrained objective, the bias might be worse. For the variance, for small enough $c$, the variances of all gradient estimators in Tables \ref{tb:var} and \ref{tb:unconstrained} are dominated by the $O(c^{-2})$ term, and so the variances are $O(\frac{n^2\sigma^2}{c^2R})$ for all unconstrained gradient estimators listed in Table \ref{tb:unconstrained}, SFE/FFE/CFE with $\boldsymbol\delta^{**},\eta = -1$, and FD$_{standard}$ and FD$_{random}$. 

In summary, from Tables \ref{tb:bias} and \ref{tb:unconstrained}, the simplex constraint does not exert too much impact on the biases of FD-type methods, and the biases of simultaneous-perturbation-type methods may depend on the distribution of the perturbation vector $\boldsymbol \delta^j$. For the variances in Tables \ref{tb:var} and \ref{tb:unconstrained}, after ``normalizing'' the length of the perturbation vector $\boldsymbol\delta^j$, we see that the (leading terms of) variances of constrained and unconstrained FD, (C)GS, (C)SS, and SFE, FFE, CFE with $\boldsymbol\delta^{**}$ are the same, i.e., $O(\frac{n^2\sigma^2}{Rc^2})$.



\section{Applications in Constrained Stochastic Approximation}\label{sec:SA}
Our main motivation to study the gradient estimators in probability simplices comes from gradient descent schemes for stochastic optimization problems that involve probability distributions as decision variables.
To make our discussion concrete, let us use the first application in Section \ref{sec:intro}, namely DRO, to describe a relevant class of optimization problems. Here, let us expand our definition of $Z$ to a multi-distribution input, $Z(P^1,P^2,\ldots,P^m)$, where $P^1,\ldots,P^m$ are independent probability distributions. This quantity $Z(P^1,P^2,\ldots,P^m)$ can denote a finite-horizon expectation-type performance measure generated from i.i.d. replications from $m$ input models. For example, for a $G/G/1$ system, $P^1,P^2$ can represent the interarrival and service time distributions,  and $Z(P^1,P^2)$ can represent the average waiting time of a range of customers. Suppose $P^1,\ldots,P^m$ are uncertain but some information is available, the distributionally robust framework considers computing the worst-case bounds for the performance measure given by
\begin{equation}
\min_{P^i\in \mathcal U^i,i\in[m]}Z(P^1,\ldots,P^m)\text{\ \ \ \ and\ \ \ \ }\max_{P^i\in \mathcal U^i,i\in[m]}Z(P^1,\ldots,P^m)\label{opt}
\end{equation}
where the set $\mathcal U^i$ is known as the uncertainty set or ambiguity set, which encodes the modeler's partial information about the input distribution. If each $\mathcal U^i$ contains the true distribution of $P^i$, then \eqref{opt} results in valid lower and upper bounds for the true performance measure. On the other hand, if $\mathcal U^i$ is a data-driven set that contains the true distribution of $P^i$ with high confidence (i.e., $\mathcal U^i$ is a confidence region), then \eqref{opt} will translate the confidence into statistically correct bounds for the true performance measure (e.g., \citealt{ghosh2019robust,ben2013robust,delage2010distributionally}). 

A common example of uncertainty sets comprises moment and support constraints, i.e., $\mathcal U^i = \{P^i:E_{P^i}[f_l^i(X^i)]\leq \mu_l^i,l=1,\ldots,s^i,\text{supp }P^i = A^i\}$ where $X^i$ is a generic random variable under distribution $P^i$, $A^i \subset \mathcal X^i=\text{dom}(X^i)$ and $f_l^i:\mathcal X^i\to \mathbb R$. For example, when $\mathcal X^i = \mathbb R$, $f_l^i(x)=x,x^2$ denotes the upper bound constraints for the first two moments respectively, and $f_l^i(x)=-x,-x^2$ can be used to specify lower bounds for the moments. As another example, one could consider a neighborhood of a baseline model measured by statistical distances such as the $\phi$-divergence, in which case $\mathcal U^i = \{P^i:d_{\phi}(P^i,P^i_b)\leq \eta^i\}$ where $d_{\phi}(P^i,P^i_b)$ denotes the $\phi$-divergence from a baseline distribution $P^i_b$ given by $d_{\phi}(P^i,P^i_b) = \int \phi(\frac{dP^i}{dP_b^i})dP^i_b$. Here $\phi(\cdot)$ is a convex function satisfying $\phi(1)=0$, e.g., the Kullback-Leibler divergence is given by $\phi(x) = x\log x -x+1$. 
Depending on the information, these constraints or combinations of them can be used in a given context.

We aim to solve problems like \eqref{opt} by integrating our derived gradient estimators with SA, when $Z$ is accessible only via black-box function evaluation. Before we show the details, we discuss two issues. One is that, in many simulation problems, $P^i$ can be a continuous distribution. However, one can suitably discretize the distribution by using support points generated from a heavy-tailed distribution and solving the problem on the sampled support points, in such a way that statistical guarantee is retained (see Theorem 1 in \citealt{ghosh2019robust}). This reduces to solving
\begin{equation}
\min_{\mathbf{p}^i\in \mathcal{U}^i,i\in[m]} Z(\mathbf{p}^1,\ldots,\mathbf{p}^m)\text{\ \ \ \ and\ \ \ \ }\max_{\mathbf{p}^i\in \mathcal{U}^i,i\in[m]} Z(\mathbf{p}^1,\ldots,\mathbf{p}^m)\label{opt1}
\end{equation}
where $\mathbf p^i$ is a finite discrete probability vector and $\mathcal U^i$ is a suitably defined uncertainty set. Second, when $m>1$, i.e., multiple models, the update procedure at each iteration of gradient descent schemes like FW and MD involves a separable optimization problem among each of the $m$ input distributions, thus retaining tractability.

In the following, we detail our zeroth-order FWSA and MDSA to solve \eqref{opt1}. For convenience, we denote $\mathbf p = vec(\mathbf p^i:i=1,\ldots,m)$ and $\mathbf q = vec(\mathbf q^i:i=1,\ldots,m)$, where $vec(\cdot)$ denotes vectorized concatenation. We use $\mathbf p_k, \mathbf q_k$ to denote the variable and the solution to the subproblem at the $k$-th iteration, and denote $\boldsymbol\delta = vec( \boldsymbol\delta^i:i=1,\ldots,m)$, $\nabla Z(\mathbf p)=vec(\nabla_{\mathbf p^i} Z(\mathbf p):i=1,\ldots,m)$ and similarly for the estimator $\hat\psi(\mathbf p)$. In addition, we let $\mathcal U = \prod_{i=1}^m \mathcal U^i$ so that $\mathbf p_k\in \mathcal U$ for all $k$. Since the maximization problem can be written as a minimization problem via simple negation, to avoid redundancy we only consider the minimization case.

\subsection{Multi-Distribution Extension of $\boldsymbol \delta^{*}$ and $\boldsymbol \delta^{**}$}
To proceed, we need to extend the formulations for $\boldsymbol \delta^{*}$ and $\boldsymbol \delta^{**}$ in Sections \ref{sec:1estimator} and \ref{sec:2estimator} from single-distribution to multi-distribution inputs. This is facilitated by noting the flexibility that for $\boldsymbol \delta^*$ there is only a lower bound requirement for $\gamma$ (Theorem \ref{thm:delta1}) and for $\boldsymbol \delta^{**}$ we have a tuning parameter $\eta$. 

For $\boldsymbol \delta^*$, we can choose $\gamma = max(\gamma_1,\ldots,\gamma_m)$ where $\gamma_l>\frac{(n^l-1)^2}{n^l(\min_{j\in [n^l]} p_j^l)^2}$, with the corresponding $C_l = \frac{n^l\gamma}{4}$. For each distribution, we get the set of $\{\theta^{n_l}\}$ and $\{\boldsymbol\alpha^{i,j}\}$ using $C_l$, and generate the perturbation $\boldsymbol \delta_l$ according to the these parameters. Then the final
perturbation is $\boldsymbol \delta = vec(\boldsymbol \delta_i:i=1,\ldots,m)$, with $\gamma$ as above. By the independence of $\{\boldsymbol\delta_l\}$, it is easy to verify that \eqref{first order}, \eqref{second order}, and \eqref{third order} are still satisfied by $\boldsymbol \delta$.

Similarly, for $\boldsymbol \delta^{**}$, we can choose $\gamma = max(\gamma_1,\ldots,\gamma_m)$ where $\gamma_l = \frac{2}{n^l(\min_{j\in [n^l]} p_j^l)^2}$ (i.e., corresponding to the case $\eta = -1$). Then we can choose $\eta_l$ such that $(n^l)^{\eta_l} = \frac{\gamma\cdot n^l(\min_{j\in [n^l]} p_j^l)^2-1}{n^l}$. This gives us the parameter $\eta_l$ for each $\boldsymbol \delta_l$, and the remaining steps are the same as $\boldsymbol \delta^{*}$ (generate each $\boldsymbol \delta_l$ independently and then concatenate them to get $\boldsymbol \delta$). 

\subsection{Frank-Wolfe Stochastic Approximation with Dirichlet Mixture Gradient Estimators}
The classical FW method uses a linear function as an approximation to the objective function at each iteration. In particular, suppose we can calculate the gradient $\nabla Z$ exactly. At each iteration $k$, given a current solution $\mathbf p_k$, FW seeks a solution to 
\begin{equation*}
\min_{\mathbf q\in \mathcal U} \nabla Z(\mathbf p_k)'(\mathbf q-\mathbf p_k)
\end{equation*}
to obtain $\mathbf q_k$, and then update the solution to $\mathbf p_{k+1}=\mathbf p_k+\epsilon_k(\mathbf q_k-\mathbf p_k)=(1-\epsilon_k)\mathbf p_k+\epsilon_k\mathbf q_k$ for some stepsize $\epsilon_k$. Note that, within the probability simplex, this update rule can be interpreted as a mixture between the current solution $\mathbf p_k$ and the subproblem solution $\mathbf q_k$.

However, in our formulation, the exact gradient $\nabla Z$ is unobserved but only noisy function evaluation is accessible. Thus, we replace $\nabla Z$ with an estimated gradient developed in the previous sections. In particular, we choose our gradient estimator $\hat \psi$ as one of $\hat \psi_{SFE}$, $\hat \psi_{FFE}$, $\hat \psi_{CFE}$, $\hat \psi_{FD,standard}$, $\hat \psi_{FD,random}$ (and with different choices of Dirichlet mixtures), and we let $\psi=E[\hat \psi]$ be its expectation. Then the optimization problem can be reformulated as 
\begin{equation}
\min_{\mathbf q\in \mathcal U} \hat \psi(\mathbf p_k)'(\mathbf q-\mathbf p_k)
\label{eq:objective_est_grad}
\end{equation}
Formulation \eqref{eq:objective_est_grad} can be decomposed into $m$ problems for each input model, so what we really need to solve is
\begin{equation*}
\min_{\mathbf q^i\in \mathcal U^i} \hat \psi_{\mathbf p^i}(\mathbf p_k)'(\mathbf q^i-\mathbf p^i_k)
\end{equation*}

\begin{algorithm}
\caption{FWSA/MDSA using FFE and Dirichlet mixture $\boldsymbol\delta^{**}$ with $\eta= -1$}
\begin{algorithmic}
\STATE \textbf{Initialization:} $\mathbf p_1\in \mathcal U$.
\STATE \textbf{Input:} Parameters $a,b,\beta,\theta,R_0$ (and $\alpha$ if use MDSA). 
\STATE \textbf{Procedure:} For each iteration $k=1,2,\ldots$, given $\mathbf p_k$: 
\STATE \textbf{1. Set parameters for the Dirichlet mixtures ($\boldsymbol\delta^{**}$)}: Compute $\gamma_k = \gamma_0 (\min_{i\in[m],j\in[n^i]} p_{k,j}^i)^{-2}$. For $l=1,\ldots, m$, compute the $\Theta^l$ and $A^l$ where $\Theta^l = (\theta^{1},\ldots,\theta^{n^l})$ and $A^l = (\boldsymbol\alpha^{1},\ldots,\boldsymbol \alpha^{n^l})$.

\STATE \textbf{2. Gradient estimation:} Set the perturbation size $c_k = b/k^{\theta}$, and the number of perturbation $R_k = [R_0k^{\beta}]+1$.
\bindent
\FOR{$j=1,\ldots,R_{k}$}
\STATE Generate $\boldsymbol\delta_l=\sum_{i=1}^{n^l}  \Theta^l_i\boldsymbol\delta^l_i$ where $\boldsymbol\delta_i^l \sim Dir(A^l_i)$. Let $\boldsymbol\delta = vec(\boldsymbol\delta_l:l=1,\ldots,m)$.
\STATE Compute the finite difference
$$\bar \psi^{j} = \frac{\gamma_k}{c_k}(\hat Z((1-c_k)\mathbf p_k+c_k\boldsymbol\delta)-\hat Z(\mathbf p_k))(\boldsymbol\delta-\mathbf p_k).$$
\ENDFOR
\eindent
\STATE Calculate $\hat \psi_{k} = \frac{1}{R_k}\sum_{j=1}^{R_k}\bar \psi^{j}$.
\STATE \textbf{3.(FWSA) Solve subproblem and update solution:} Set $\epsilon_k = a/k$. Split $\hat \psi_{k} = vec(\hat \psi^l_k,l=1,\ldots, m)$. For $l=1,\ldots,m$, solve
$$\mathbf q_{k+1}^l \in \arg\min_{\mathbf q^l\in\mathcal U^l} (\hat \psi^l_k)'(\mathbf q^l-\mathbf p_k^l).$$
Let $\mathbf q_{k+1} = vec(\mathbf q_{k+1}^l:l=1,\ldots,m)$, update $\mathbf p_{k+1} = (1-\epsilon_k)\mathbf p_{k}+ \epsilon_k \mathbf q_{k+1}$.

\STATE \textbf{3.(MDSA)  Solve subproblem and update solution:} Set $\rho_k = a/k^{\alpha}$. Split $\hat \psi_{k} = vec(\hat \psi^l_k,l=1,\ldots, m)$. For $l=1,\ldots,m$, solve 
$$\mathbf p_{k+1}^l \in \arg\min_{\mathbf q^l\in\mathcal U^l} \rho_k(\hat \psi^l_k)'(\mathbf q^l-\mathbf p_k^l)+V(\mathbf p_k^l,\mathbf q^l).$$
Update $\mathbf p_{k+1} = vec(\mathbf p_{k+1}^l:l=1,\ldots,m)$.
\end{algorithmic} 
\label{algo_FWSA}
\end{algorithm}

Algorithm \ref{algo_FWSA} describes FWSA (and MDSA to be presented later) in detail, using FFE with $\boldsymbol\delta^{**}$ as the gradient estimator. Notice that the definition of each $\gamma_k$ (where $\gamma_k$ is the $\gamma$ used in the score function in the gradient estimator at iteration $k$) involves the inverse square of the smallest entry of $\mathbf p$, thus can blow up if the latter is too small. However, the update rule indicates that each new $\mathbf p_{k+1}$ is a convex mixture of the current $\mathbf p_k$ and another vector $\mathbf q_k$ lying inside the product of probability simplices. This turns out to make the growth of $\gamma_k$ controllable. The following result presents this controllable growth.
\begin{lemma}
Let $Z(\cdot):\prod_{i=1}^m\mathcal P^{n^i}\to \mathbbm R$, and let the gradient update rule be $\mathbf p_{k+1}=(1-\epsilon_k)\mathbf p_{k}+\epsilon_k\mathbf q_{k}$ where $\mathbf q_k\in\prod_{i=1}^m\mathcal P^{n^i}$. In addition, we set the update stepsize $\epsilon_k = \frac{a}{k}$ for some $a<1/2$, and denote $p_0 = \min_{i\in[m],j\in n^i} p^i_{1,j}$. Then 
\begin{equation*}
\min_{i\in[m],j\in n^i} p^i_{k,j}\geq \frac{A}{k^a}p_0 
\end{equation*}
for some $A>0$.
\label{lm:bdd_FWSA}
\end{lemma}

Based on the above result, we show the convergence of FWSA. We first state the following assumptions:
\begin{assumption}[Smoothness of the Frank-Wolfe gap]
Let $g(\cdot):\prod_{i=1}^m\mathcal P^{n^i}\to \mathbbm R$ be the Frank-Wolfe gap $g(\mathbf p) = -\min_{\mathbf q \in \mathcal U}\nabla Z(\mathbf p)(\mathbf q-\mathbf p)$. Assume that $g(\cdot)$ is continuous and an optimal solution exists.\label{assumption:FW gap}
\end{assumption}

\begin{assumption}[Solution uniqueness] 
There exists a unique minimizer $\mathbf p^*$ for $\min_{\mathbf p\in \mathcal U} Z(\mathbf p)$, and $\mathbf p^*$ is the only feasible solution such that $g(\mathbf p^*)=0$.\label{assumption:uniqueness}
\end{assumption}

Assumptions \ref{assumption:FW gap} and \ref{assumption:uniqueness} ensure the smoothness of the FW gap function, and the uniqueness of the optimal solution which is the only feasible solution with a zero FW gap.


\begin{theorem}[FWSA convergence]
Let $\hat \psi$ be one of our considered gradient estimators (SFE, FFE, CFE, FD$_{standard}$, FD$_{random}$) with expectation $E[\hat \psi] =\psi$, and a corresponding set of parameters, that for $c$ small enough satisfies
\begin{equation*}
\lVert \psi(\mathbf p) - \nabla Z(\mathbf p) -\epsilon_0(\mathbf p) \mathbbm 1\rVert \leq B_1c
\end{equation*}
for some $\epsilon_0(\mathbf p) \in \mathbbm R$ that depends on $\mathbf p$, and
\begin{equation*}
E[\lVert \hat \psi - \psi \rVert^2]\leq B_2\frac{\gamma}{Rc^2}
\end{equation*}

Suppose we choose the update stepsize as $\epsilon_k = a/k$ for some $a\leq 1/2$, and the perturbation size as $c_k = b/k^{\theta}$ for some $\theta$. For each iteration, the number of simulation replications is $R_{k} = [R_0k^{\beta}]+1$ where $\beta > 2(a+\theta)$. In addition, we choose $\gamma_k = \gamma_0 (\min_{i\in[m]}\min_{j\in[n^i]} p_{k,j}^i)^{-2}$.

Then, with the above configuration, under Assumptions \ref{assumption:smooth}, \ref{assumption:smooth_noise}, and \ref{assumption:FW gap}, $D(Z(\mathbf p_k),\mathcal L^*)\to 0$ a.s. where $\mathcal L^* = \{Z(\mathbf p):~\mathbf p~\text{satisfies}~g(\mathbf p)=0\}$ and $D(x,A)=\inf_{y\in A}\lVert x-y \rVert$ for any point $x$ and set $A$ in the Euclidean space.

In addition, if Assumption \ref{assumption:uniqueness} also holds, then $\mathbf p_k$ converges to $\mathbf p^*$ a.s.
\label{main_FWSA}
\end{theorem}

We make three comments:
\begin{enumerate}[leftmargin=*]
\item Regarding the assumption for $\gamma_k$, note that when there is only one distribution (with dimension $n$), $\boldsymbol\delta^*$ gives $\gamma^* = n (\min_{j\in[n]} p_{j})^{-2}$, and $\boldsymbol\delta^{**}$ with $\eta = -1$ gives $\gamma^* = \frac{2}{n} (\min_{j\in[n]} p_{j})^{-2}$. So $\gamma_k = O((\min_{j\in[n]} p_{j})^{-2})$ appears to be a reasonable assumption. 

\item We assume that the bias is $O(c)$ because Tables \ref{tb:bias} and \ref{tb:var} and the accompanied discussions show that FFE with $\boldsymbol\delta^{**}$ is the best among SFE, FFE and CFE in terms of variance, and for this setting we have bias $O(c)$. For the case where bias is $O(c^2)$, results like Theorem \ref{main_FWSA} still hold.

\item For the choice of parameters $a, \theta, \beta$, we notice that the squared bias by assumption is $O(c^2)$ and variance is $O(\frac{\gamma }{Rc^2})$. To minimize the sum of these two values, $c = O((\frac{\gamma_k }{R_k})^{1/4})$ seems to be a good choice, and this translates into $-4\theta = 2a-\beta$. For example, if we choose $\theta = 1/8$, $a = 1/4$, then $\beta = 1$.
\end{enumerate}



\subsection{Mirror Descent Stochastic Approximation with Dirichlet Gradient Estimator}
Next we consider integrating our gradient estimators into MDSA.
To begin, we recall MD in the setting of non-stochastic, single-distributional problems. That is, we consider $Z(\cdot):\mathcal P^n\to \mathbbm R$. At each iteration $k$, we update the variable $\mathbf p_k$ using $\mathbf p_{k+1} = prox_{\mathbf p_k}(\rho_k\nabla Z(\mathbf p_k))$, where $\rho_k>0$ and the prox-mapping is defined as 
\begin{equation*}
prox_{\mathbf x}(\mathbf y)=\arg\min_{\mathbf z\in \mathcal U}\mathbf y'(\mathbf z-\mathbf x)+V(\mathbf x,\mathbf z)
\end{equation*}
Here $V(\cdot,\cdot):\mathcal P^n\times \mathcal P^n\to \mathbbm R$ is the remainder of the first order approximation of a strongly convex function $w(\cdot):\mathcal P^n\to \mathbbm R$. To be precise, we have, for all $\mathbf x, \mathbf y\in \mathcal P^n$,
\begin{equation*}
w(\mathbf y)-w(\mathbf x)\geq \nabla w(\mathbf x)'(\mathbf y - \mathbf x)+\frac{\alpha}{2}\lVert \mathbf y - \mathbf x\rVert^2
\end{equation*}
In particular, we take $w(\cdot)$ as the entropic distance-generating function (with $\alpha=1$): $w(\mathbf p) = \sum_{i=1}^n p_i\log(p_i)$. Given this $w(\cdot)$ we define
\begin{equation*}
V(\mathbf p,\mathbf q) = w(\mathbf q)-w(\mathbf p)-\nabla w(\mathbf p)'(\mathbf q-\mathbf p) = \sum_{i=1}^n q_i\log(\frac{q_i}{p_i})
\end{equation*}
The update rule is 
\begin{equation*}
\mathbf p_{k+1} = prox_{\mathbf p_k}(\rho_k \nabla Z(\mathbf p_k))=\arg\min_{\mathbf q\in \mathcal U}\rho_k\nabla Z(\mathbf p_k)'(\mathbf q-\mathbf p_k)+V(\mathbf p_k,\mathbf q)
\end{equation*}
This procedure is known as the entropic descent algorithm \citep{beck2003mirror}. Like in the previous subsection, here we consider situations where we only have noisy, biased gradient estimator $\hat \psi(\mathbf p_k)$. Thus, similar to before, we replace the true gradient with its estimator, so the update rule becomes
\begin{equation*}
\mathbf p_{k+1} = prox_{\mathbf p_k}(\rho_k \hat\psi (\mathbf p_k))=\arg\min_{\mathbf q\in \mathcal U}\rho_k\hat \psi(\mathbf p_k)'(\mathbf q-\mathbf p_k)+V(\mathbf p_k,\mathbf q)
\end{equation*}

We show the a.s. convergence of MDSA. For this, we need two alternate assumptions:
\begin{assumption}[Solution properties]
The problem $\min_{\mathbf q\in \mathcal U} Z(\mathbf q)$ has a unique optimal solution $\mathbf p^*\in \mathcal U$ such that for any feasible $\mathbf q \in \mathcal U,\mathbf q \neq \mathbf p^*$, it holds that $\nabla Z(\mathbf q)'(\mathbf q -\mathbf p^*)>0$.\label{assumption:properties}
\end{assumption}

\begin{assumption}[Decay of smallest component]
We have $\min_{l\in[n],j\in[n^l]} p_{k,j}^l \geq \frac{\bar A}{k^{d}}$ for some $\bar A,d>0$.\label{assumption:decay}
\end{assumption}

Assumption \ref{assumption:properties} ensures the uniqueness of the optimal solution that satisfies a standard first-order condition. Assumption \ref{assumption:decay} ensures that the decay of the smallest component of $\mathbf p_k$ is controllable.


\begin{theorem}[MDSA convergence]



Let $\hat \psi$ be one of our gradient estimators (SFE, FFE, CFE) with expectation $E[\hat \psi] =\psi$, and a corresponding set of parameters, that for $c$ small enough satisfies
\begin{equation*}
\lVert \psi(\mathbf p) - \nabla Z(\mathbf p) -\epsilon_0(\mathbf p) \mathbbm 1\rVert \leq B_1c
\end{equation*}
for some $\epsilon_0(\mathbf p) \in \mathbbm R$ that depends on $\mathbf p$
\begin{equation*}
E[\lVert \hat \psi - \psi \rVert^2]\leq B_2\frac{\gamma}{Rc^2}
\end{equation*}

Suppose we choose the update stepsize as $\rho_k = a/k^{\alpha}$ for some $a, \frac{1}{2}< \alpha\leq1 $, and the perturbation size as $c_k = b/k^{\theta}$ such that $\alpha+\theta>1$. For each iteration, the number of simulation replications is $R_{k} = [R_0k^{\beta}]+1$ where $2\alpha+ \beta -2d-2\theta>1$. In addition, we choose $\gamma_k = \gamma_0 (\min_{i\in[m]}\min_{j\in[n^i]} p_{k,j}^i)^{-2}$.

Then, with the above configuration, under Assumptions \ref{assumption:smooth},  \ref{assumption:smooth_noise}, \ref{assumption:properties} and \ref{assumption:decay}, $\mathbf p_k$ converges to $\mathbf p^*$ a.s.
\label{main_MDSA}
\end{theorem}


\section{Numerical Results}\label{sec:numerics}
We conduct numerical experiments to compare the performances of our gradient estimators and their integration into FWSA and MDSA schemes. We divide our numerical investigations into two parts. Section \ref{section:exp1} first studies the statistical performances of our gradient estimators using SFE, FFE, CFE and the two Dirichlet mixtures $\boldsymbol\delta^*$ and $\boldsymbol\delta^{**}$. Section \ref{sec:opt num} compares our gradient estimators with benchmark FD schemes on two optimization settings, one on Rosenbrock objectives (Section \ref{section:exp2}), and one on a single-server-queue model (Section \ref{section:exp3}). In these experiments we will consider $\phi$-divergence and moment-based constraints. 
Also, recall in Section \ref{sec:Dir} that our Dirichlet mixture construction allows some mixture components to have a smaller support than the original base probability distribution. In our experiments, for convenience we add a small positive adjustment to all entries of $\boldsymbol\alpha^{l,i}$ and $\boldsymbol\alpha^{n}$ in these Dirichlet mixtures such that the supports are all the same (as the original). This small adjustment has very negligible effect on all our calculations, and \eqref{eqn:delta_moment1}, \eqref{eqn:delta_moment2} and Lemma \ref{lemma_third_order} which hold for regular Dirichlet random vectors still hold approximately for these modified Dirichlet random vectors.



\subsection{Performances of Gradient Estimators}
\label{section:exp1}
\subsubsection{Moment Properties of Dirichlet Mixtures.}
The key properties for the perturbation vector $\boldsymbol\delta$ are the moment properties \eqref{first order}, \eqref{second order} and \eqref{third order}. The formulations of $\boldsymbol\delta^*$ and $\boldsymbol\delta^{**}$ are designed to achieve these properties. 

\begin{figure}[h!]
  \centering
  \begin{subfigure}[b]{0.25\linewidth}
    \includegraphics[width=\linewidth]{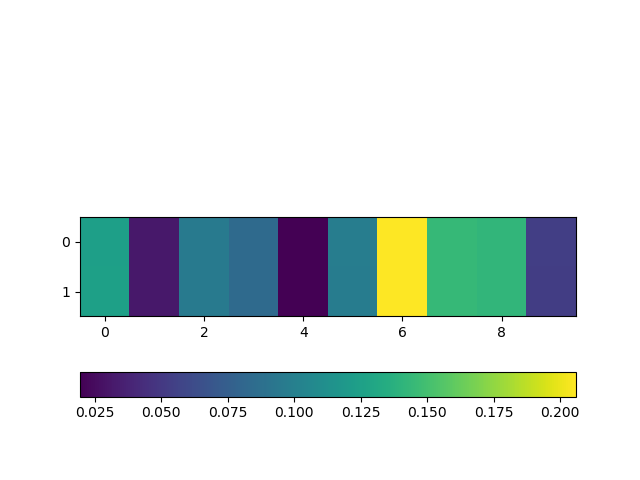}
    \caption{}
  \end{subfigure}
  \begin{subfigure}[b]{0.25\linewidth}
    \includegraphics[width=\linewidth]{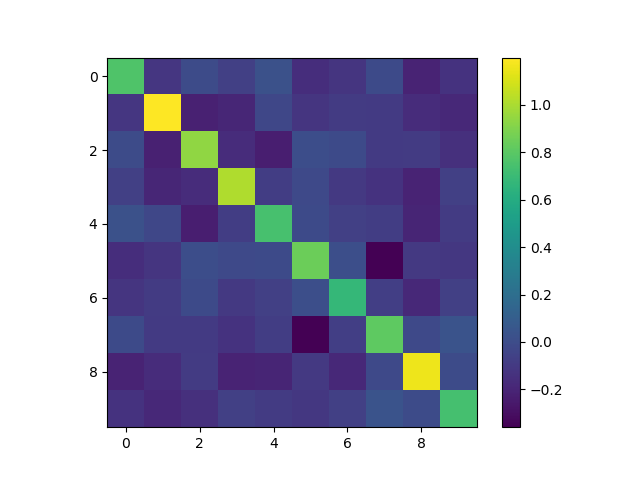}
    \caption{}
  \end{subfigure}
  \begin{subfigure}[b]{0.25\linewidth}
    \includegraphics[width=\linewidth]{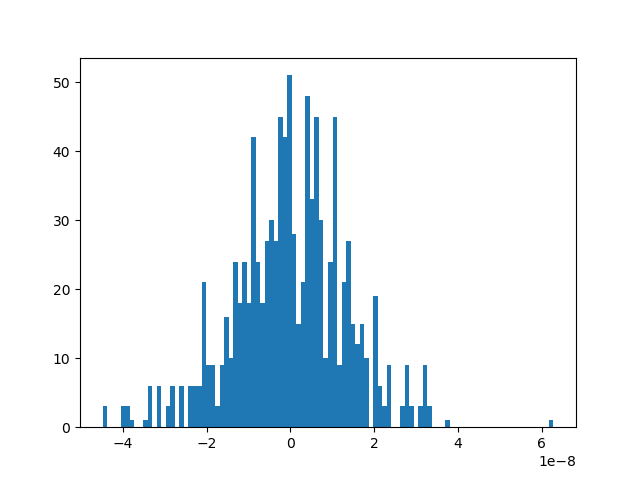}
    \caption{}
  \end{subfigure}
  \begin{subfigure}[b]{0.25\linewidth}
    \includegraphics[width=\linewidth]{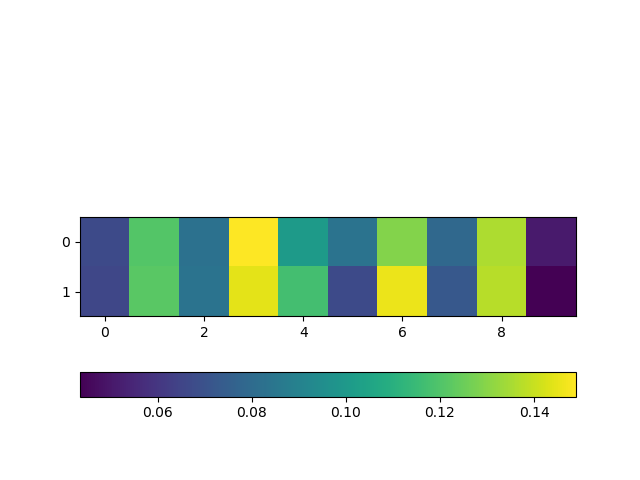}
    \caption{}
  \end{subfigure}
  \begin{subfigure}[b]{0.25\linewidth}
    \includegraphics[width=\linewidth]{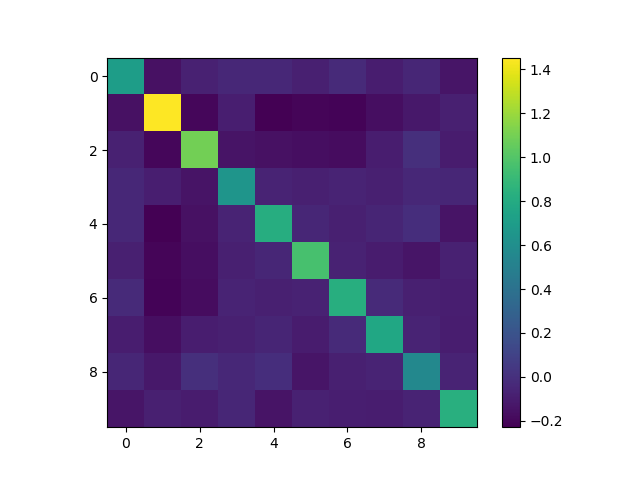}
    \caption{}
  \end{subfigure}
  \begin{subfigure}[b]{0.25\linewidth}
    \includegraphics[width=\linewidth]{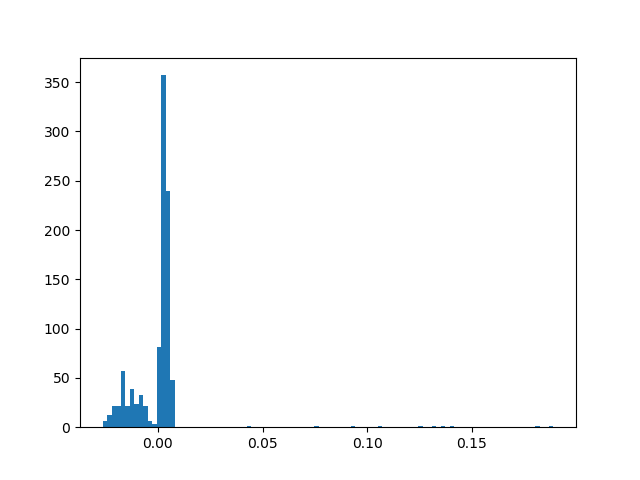}
    \caption{}
  \end{subfigure}
\caption{Distributions for (a,b,c): $\boldsymbol\delta^*$ and (d,e,f): $\boldsymbol\delta^{**}$}
\label{fig:fig0}
\end{figure}

Our first experiment shows the first three moments of $\boldsymbol\delta^*$ and $\boldsymbol\delta^{**}$. We first set the base probability vector $\mathbf p$ by randomly generating each component from $Uniform(0,1)$ and then normalizing them. In each of Figure \ref{fig:fig0} (a) and (d), the first row shows this base vector $\mathbf p$, and the second row is an approximation of $E[\boldsymbol\delta^*]$ or $E[\boldsymbol\delta^{**}]$ by averaging 100 repetitions. We see that for both $\boldsymbol\delta^*$ and $\boldsymbol\delta^{**}$, the first moment is close to $\mathbf p$, validating the condition \eqref{first order}. Figure \ref{fig:fig0} (b) and (e) plot the 100-average approximations of $\gamma E[(\boldsymbol\delta^*-\mathbf{p})(\boldsymbol\delta^*-\mathbf{p})']$ and $\gamma E[(\boldsymbol\delta^{**}-\mathbf{p})(\boldsymbol\delta^{**}-\mathbf{p})']$, which appear close to $\mathbf I - \mathbbm 1\mathbbm 1'/n$, showing that \eqref{second order} is also satisfied. Figure \ref{fig:fig0} (c) and (f) present the histograms of the 100-average approximations of $E[(\boldsymbol\delta^*-\mathbf{p})_i(\boldsymbol\delta^*-\mathbf{p})_j(\boldsymbol\delta^*-\mathbf{p})_k]$ and $E[(\boldsymbol\delta^{**}-\mathbf{p})_i(\boldsymbol\delta^{**}-\mathbf{p})_j(\boldsymbol\delta^{**}-\mathbf{p})_k]$ for all $(i,j,k)\in[10]^3$. Recall that for our gradient estimators to achieve $O(c)$ bias, conditions \eqref{first order} and \eqref{second order} suffice. But to achieve $O(c^2)$ bias (for $\psi_{SFE}$ and $\psi_{FFE}$), \eqref{third order} needs to be satisfied. One difference between $\boldsymbol\delta^{*}$ and $\boldsymbol\delta^{**}$ is that $\boldsymbol\delta^{*}$ also satisfies the third moment constraints \eqref{third order} and here we visualize this. We see that the third moments of $\boldsymbol\delta^*$ vanish ($\sim 10^{-8}$) which is as expected, while the third moments of $\boldsymbol\delta^{**}$ are $\sim 0.01$ and have outliers since it is designed to have non-constant third moments (in exchange for a better variance). The above confirms our calculations on the moment properties of $\boldsymbol\delta^*$ and $\boldsymbol\delta^{**}$.


\subsubsection{Estimation Variances of Dirichlet Mixtures.}
The moment properties determine the performance of the gradient estimators in expectation. Another important statistical property is the variance, and this is influenced by the type of gradient estimator plus many  other parameters, including the noise level $\sigma$, number of function evaluation $R$, perturbation size $c$, problem dimension $n$, and the multiplier $\gamma$ which is used in the score function $S(\mathbf p,\boldsymbol\delta)=\gamma(\boldsymbol\delta-\mathbf p)$.

In this set of experiments, we test 5 gradient estimators $\psi_{SFE,\boldsymbol\delta^*}$, $\psi_{FFE,\boldsymbol\delta^*}$, $\psi_{FFE,\boldsymbol\delta^{**}}$, $\psi_{CFE,\boldsymbol\delta^{**}}$, and $\psi_{FD,random}$ for the simple quadratic objective function (Here we use the notation $\psi_{\text{method},\text{perturbation}}$ where ``method" refers to SFE, FFE, CFE or FD, and ``perturbation" refers to the Dirichlet mixture $\boldsymbol\delta^*$ or $\boldsymbol\delta^{**}$, or in the case of FD, the random-dimension specification that is shorthanded as ``random"):
\begin{equation*}
    Z(\mathbf p) = \lVert \mathbf p - \mathbbm 1/n \rVert^2
\end{equation*}
We take our observation as $\hat Z(\mathbf p) = Z(\mathbf p) +X$, where $X\sim \mathcal N(0,\sigma^2)$.

For this quadratic objective, we have by our definition of the gradient that 
\begin{align*}
    \nabla Z(\mathbf p)_i &= \lim_{\epsilon\to 0^+}\frac{Z((1-\epsilon)\mathbf p+\epsilon\mathbf e_i)-Z(\mathbf p)}{\epsilon} \nonumber\\
    &= \langle 2(\mathbf p-\mathbbm 1/n),\mathbf e_i-\mathbf p\rangle =2p_i-2\lVert \mathbf p\rVert^2
\end{align*}

Thus we have $\nabla Z(\mathbf p) = 2\mathbf p - 2\lVert \mathbf p\rVert^2\mathbbm 1 $. The constant $M_1$ and $L_1$ are bounded above by
\begin{equation*}
    M_1\leq \sup_{\mathbf p\in \mathcal P}\inf_{\epsilon\in \mathbbm R}\lVert \nabla Z(\mathbf p)-\epsilon \mathbbm 1\rVert \leq \sup_{\mathbf p\in \mathcal P}2\lVert \mathbf p\rVert \leq 2
\end{equation*}
and
\begin{equation}
       \inf_{\epsilon \in \mathbbm R} \lVert \nabla Z(\mathbf p) - \nabla Z(\mathbf q)-\epsilon \mathbbm 1\rVert  \leq 2\lVert \mathbf p - \mathbf q \rVert
\end{equation}
and so $L_1 \leq 2$. 


Based on the comments below Tables \ref{tb:bias} and \ref{tb:var}, we have that $\gamma = O(n^3)$ for $\boldsymbol\delta^*$ and $\gamma = O(n)$ for $\boldsymbol\delta^{**}$. With $M_1,L_1\leq 2$, we expect the results in Table \ref{tb:exp_expectation}.

\begin{table}[h!]
  \begin{center}
    \caption{Order of variance for different estimators}
    \begin{tabular}{c||c|c|c|c|c}
      &$\psi_{SFE,\boldsymbol\delta^*}$&$\psi_{FFE,\boldsymbol\delta^*}$& $\psi_{FFE,\boldsymbol\delta^{**}}$& $\psi_{CFE,\boldsymbol\delta^{**}}$& $\psi_{FD,random}$\\
      \hline
Variance&$O(\frac{n^4(\sigma^2+1)}{Rc^2})$&$O(\frac{n^4(\sigma^2+c^2+c^4)}{Rc^2})$&$O(\frac{n^2(\sigma^2+c^2+c^4)}{Rc^2})$&$O(\frac{n^2(\sigma^2+c^2+c^4)}{Rc^2})$&$O(\frac{n^2(\sigma^2+c^2+c^4)}{Rc^2})$
    \end{tabular}
    \label{tb:exp_expectation}
  \end{center}
\end{table}

The setup of the experiment is the following: we vary $\sigma$, $R$, $c$, $n$ one at a time while keeping all other parameters fixed. Table \ref{tb:parameters} shows our choices of the parameters. And for each set of parameters $(\sigma, R,c,n)$, we set 20 base probability vector by randomly generating $\mathbf p_i\sim Dir(10\cdot\mathbbm 1)$ (note that for Dirichlet distribution $Dir(\alpha \mathbbm 1)$, $\alpha = 1$ corresponds to a uniform distribution over the simplex, and the larger $\alpha$ is, the more concentrated the distribution is toward the center of the simplex). For each $\mathbf p_i$, we try the 5 gradient estimators mentioned above for 50 times and get $\psi_s^{i,j},~i=1,\ldots,20,~j=1,\ldots,50$ and $s\in \{(SFE,\boldsymbol\delta^*),(FFE,\boldsymbol\delta^*),(FFE,\boldsymbol\delta^{**}),(CFE,\boldsymbol\delta^{**}) ,(FD,random)\}$. We calculate the mean $\bar{\psi}_{s,i} = \frac{1}{50}\sum_{j=1}^{50} \psi^{i,j}_s$, and use $v_{s,i} = \frac{1}{49}\sum_{j=1}^{50}\lVert \psi^{i,j}_s - \bar{\psi}_{s,i} \rVert^2$ as an estimate for the variance at this particular $\mathbf p_i$. We also use $v_s = \frac{1}{20}\sum_{i=1}^{20} v_{s,i}$ as an estimate for the variance of a particular gradient estimator $s$.

\begin{figure}[h!]
  \centering
  \begin{subfigure}[b]{0.45\linewidth}
    \includegraphics[width=\linewidth]{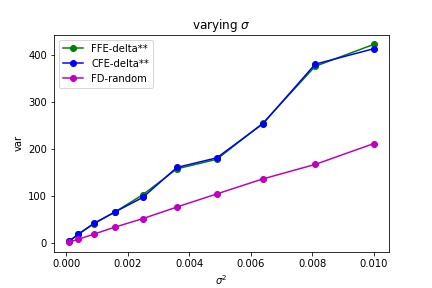}
    \caption{}
  \end{subfigure}
  \begin{subfigure}[b]{0.45\linewidth}
    \includegraphics[width=\linewidth]{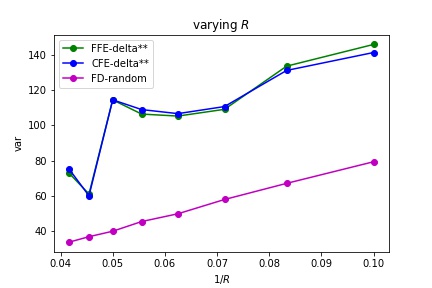}
    \caption{}
  \end{subfigure}
  \newline
  \begin{subfigure}[b]{0.45\linewidth}
    \includegraphics[width=\linewidth]{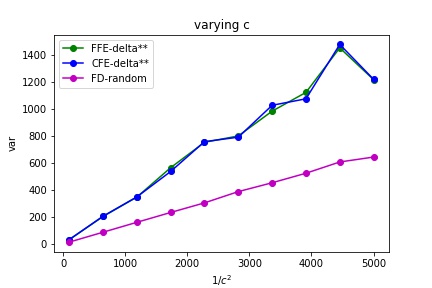}
    \caption{}
  \end{subfigure}
  \begin{subfigure}[b]{0.45\linewidth}
    \includegraphics[width=\linewidth]{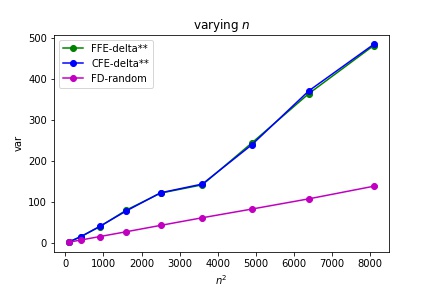}
    \caption{}
  \end{subfigure}
\caption{Variance against $\sigma,R,c,N$ for $(FFE,\boldsymbol\delta^{**}),(CFE,\boldsymbol\delta^{**}) ,(FD,random)$}
\label{fig:var}
\end{figure}

\begin{figure}[h!]
  \centering
  \begin{subfigure}[b]{0.45\linewidth}
    \includegraphics[width=\linewidth]{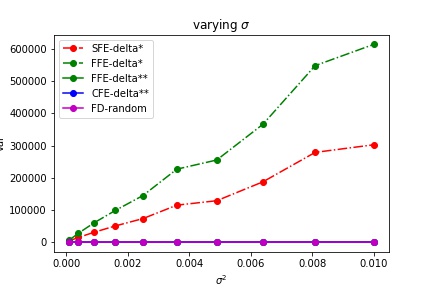}
    \caption{}
  \end{subfigure}
  \begin{subfigure}[b]{0.45\linewidth}
    \includegraphics[width=\linewidth]{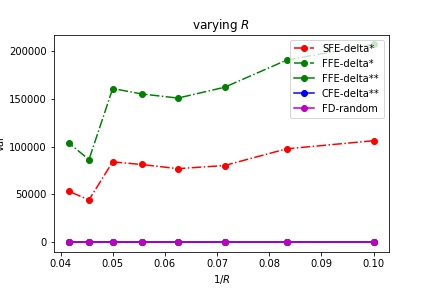}
    \caption{}
  \end{subfigure}
    \newline
  \begin{subfigure}[b]{0.45\linewidth}
    \includegraphics[width=\linewidth]{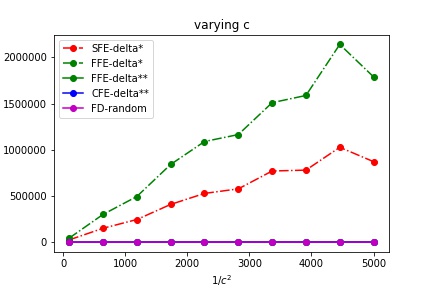}
    \caption{}
  \end{subfigure}
  \begin{subfigure}[b]{0.45\linewidth}
    \includegraphics[width=\linewidth]{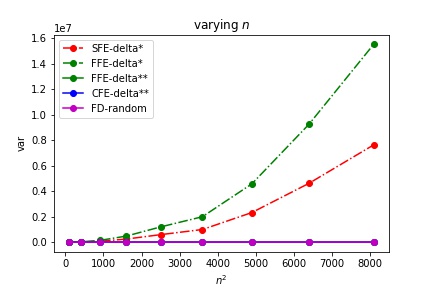}
    \caption{}
  \end{subfigure}
\caption{Variance against $\sigma,R,c,n$ for all considered estimators}
\label{fig:var_full}
\end{figure}

Figure \ref{fig:var} shows that the resulting variances are proportional to $ \sigma^2,1/R,c^{-2}$ and $n^2$ for $\psi_{FFE,\boldsymbol\delta^{**}}$, $\psi_{CFE,\boldsymbol\delta^{**}}$, and $\psi_{FD,random}$. This validates our analysis that for $c<<\sigma $, the variance is $O(\frac{\gamma n\sigma^2}{Rc^2})$. The variances of $\psi_{CFE,\boldsymbol\delta^{*}}$ and $\psi_{FFE,\boldsymbol\delta^{*}}$ turn out to be much larger than variances of the other three gradient estimators. For this reason we do not show them in Figure \ref{fig:var}, but in Figure \ref{fig:var_full}.

\begin{table}[h!]
  \begin{center}
    \caption{Choices of $\sigma, R,c,n$ for experiments on the gradient estimation variance}
    \begin{tabular}{c||c|c|c|c}
      &$\sigma$&$R$& $c$& $n$\\
      \hline
      vary $\sigma$&$\{0.01,0.02,\ldots,0.1\}$&15&0.05&20\\
      vary $R$&0.05&$\{10,12,\ldots,24\}$&0.05&20\\
      vary $c$&0.05&15&$\{0.0141,0.0150,\ldots,0.1000\}$&20\\
      vary $n$&0.05&30&0.1&$\{10, 20,\ldots, 90\}$
    \end{tabular}
    \label{tb:parameters}
  \end{center}
\end{table}

In Figure \ref{fig:var_full}, variances of estimators with $\boldsymbol\delta^{**}$ are on the scale of $10^3$ while the variances of $\boldsymbol\delta^*$ are on the scale of $10^5\sim10^6$. This is in agreement with our theoretical result that $\gamma = O(n^3)$ for $\boldsymbol\delta^*$ and $\gamma = O(n)$ for $\boldsymbol\delta^{**}$. From this set of experiments, among all the six combinations of gradient estimators and Dirichlet mixtures we propose, it appears that the best one is $\psi_{FFE,\boldsymbol\delta^{**}}$. We will use $\psi_{FFE,\boldsymbol\delta^{**}}$ in our experiments regarding optimization performances next. 

\subsection{Performances of FWSA and MDSA}\label{sec:opt num}
We investigate the numerical performances in integrating $\psi_{FFE,\boldsymbol\delta^{**}}$ into FWSA and MDSA. We also compare with $\psi_{FD,random}$ as a benchmark. 
\subsubsection{Rosenbrock Objective with Divergence Constraint}
\label{section:exp2}
In this set of experiments, we aim to solve for the problem
\begin{align*}
    &~~~~~\min_{\mathbf p \in \mathcal P} Z(\mathbf p):= f(\mathbf p +(1-1/n)\mathbbm 1)\\
    &s.t.~ \mathbf p \in \mathcal U := \{\mathbf q \in \mathcal P| d_{\phi}(\mathbf q, \mathbf p_b)\leq 100\}
\end{align*}
where $f$ is the $n$-dimensional Rosenbrock objective
\begin{equation*}
    f(\mathbf x) := \sum_{i=1}^{n-1}100(x_{i+1} - x_i^2)^2+(1-x_i)^2
\end{equation*}
and $\mathbf p_b$ is the baseline generated by $\mathbf p_b = \frac{\mathbf q_b}{\mathbbm 1'\mathbf q_b}$, $(\mathbf q_b)_i \sim 1+Uniform(0,1),~i=1,2,\ldots,n$. Here we add 1 to the uniform distribution so that $\min_i (\mathbf q_b)_i$ would not be too small, i.e., $\mathbf p_b$ is relatively far from the boundary of the probability simplex. In this case, the uncertainty set $\mathcal U$ and the $k-$iteration variable $\mathbf p_k$ are bounded away from the boundary, and so $\gamma$ for $\boldsymbol\delta^{*}$ and $\boldsymbol\delta^{**}$ are guaranteed to be not too large for all iterations. We further assume that the observation $\hat Z(\mathbf p) = Z(\mathbf p) + X$, where $X\sim \mathcal N(0,\sigma^2)$. 

We test $\psi_{FFE,\boldsymbol\delta^{**}}$ and 
$\psi_{FD,random}$ with different dimension $n$, number of function evaluations in each round $R$, and noise level $\sigma$. For each set of parameters $(n,R,\sigma)$, we generate 12 $\mathbf p_b$'s and apply $\psi_{FFE,\boldsymbol\delta^{**}}$ and $\psi_{FD,random}$ with MDSA to the optimization problem. In the MDSA, we set $a = 0.005$, $\alpha = 1$, $b=4/n$, $\theta = 0.25$, $\beta = 0$, so $R_k = R_0 = R$, and we set the maximum number of iterations to 50. We initialize with $\mathbf p_0 = \mathbf p_b\in \mathcal U$.
\begin{figure}[h!]
  \centering
  \begin{subfigure}[b]{0.3\linewidth}
    \includegraphics[width=\linewidth]{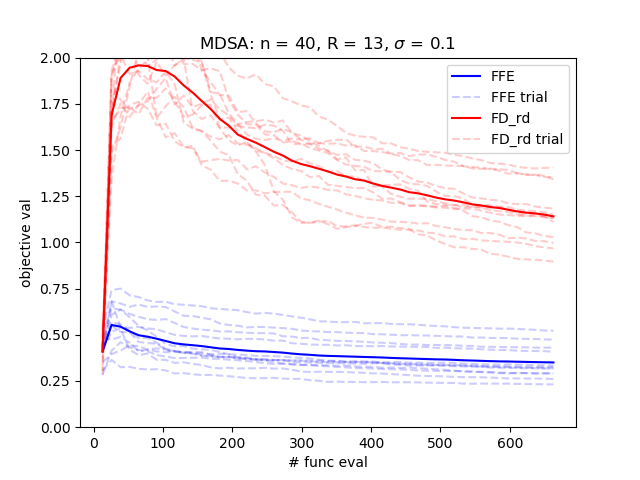}
    \caption{}
  \end{subfigure}
    \begin{subfigure}[b]{0.3\linewidth}
    \includegraphics[width=\linewidth]{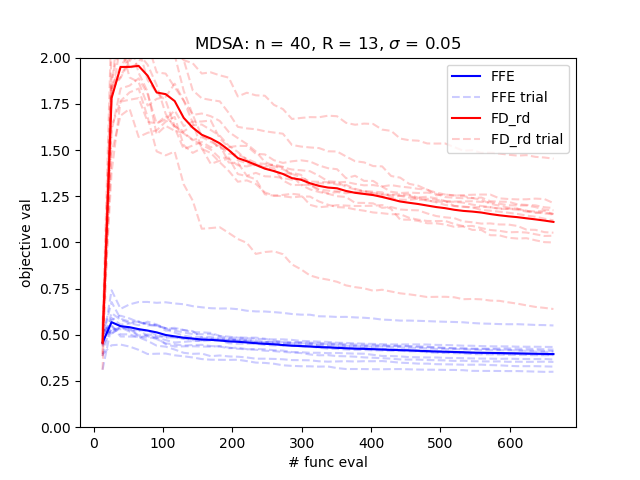}
    \caption{}
  \end{subfigure}
    \begin{subfigure}[b]{0.3\linewidth}
    \includegraphics[width=\linewidth]{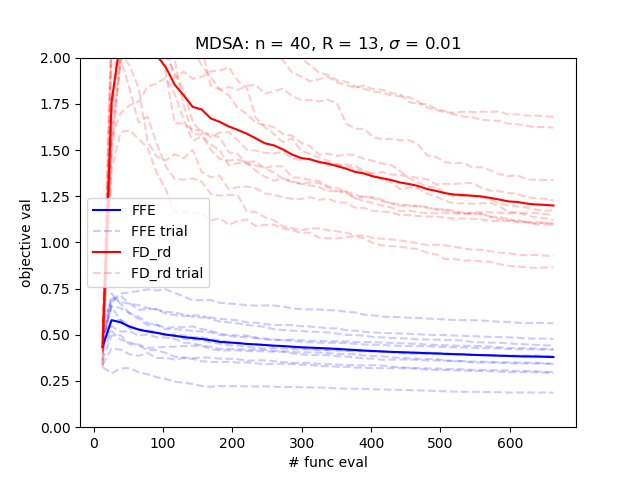}
    \caption{}
  \end{subfigure}
  
    \begin{subfigure}[b]{0.3\linewidth}
    \includegraphics[width=\linewidth]{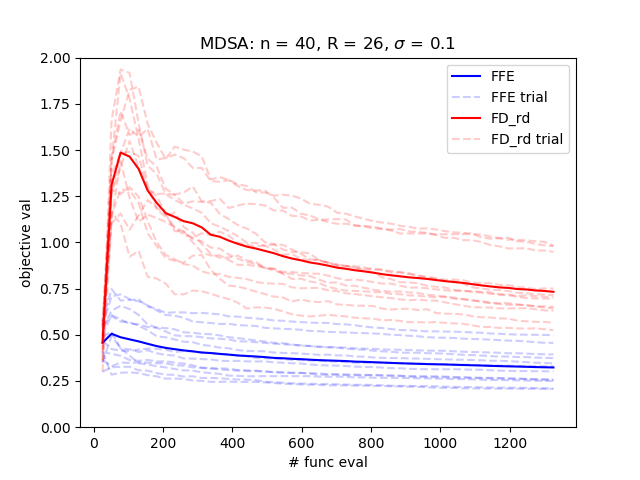}
    \caption{}
  \end{subfigure}
    \begin{subfigure}[b]{0.3\linewidth}
    \includegraphics[width=\linewidth]{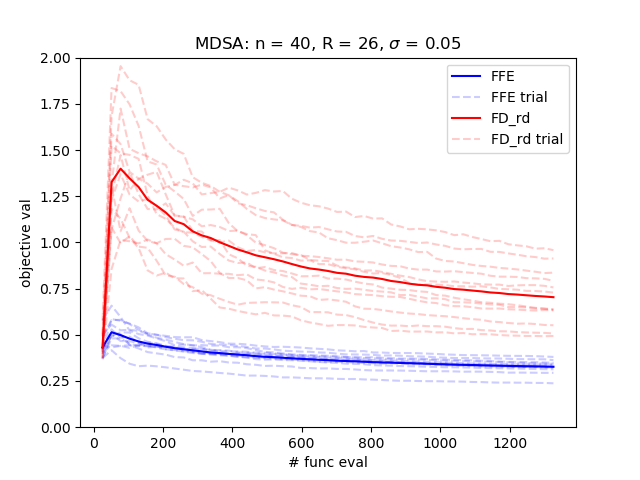}
    \caption{}
  \end{subfigure}
    \begin{subfigure}[b]{0.3\linewidth}
    \includegraphics[width=\linewidth]{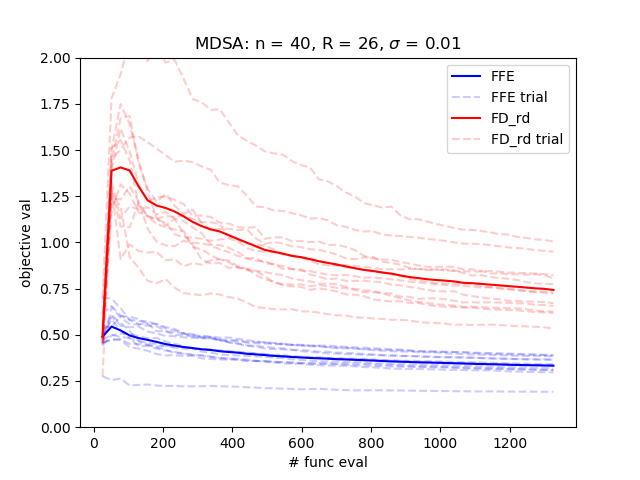}
    \caption{}
  \end{subfigure}
\caption{Performances of MDSA using $\psi_{FFE,\boldsymbol\delta^{**}}$ and $\psi_{FD,random}$ for Rosenbrock objectives, $n=40$, dashed line: 12 trials, solid line: average of 12 trials}
\label{fig:rosen40}
\end{figure}

\begin{figure}[h!]
  \centering
  \begin{subfigure}[b]{0.3\linewidth}
    \includegraphics[width=\linewidth]{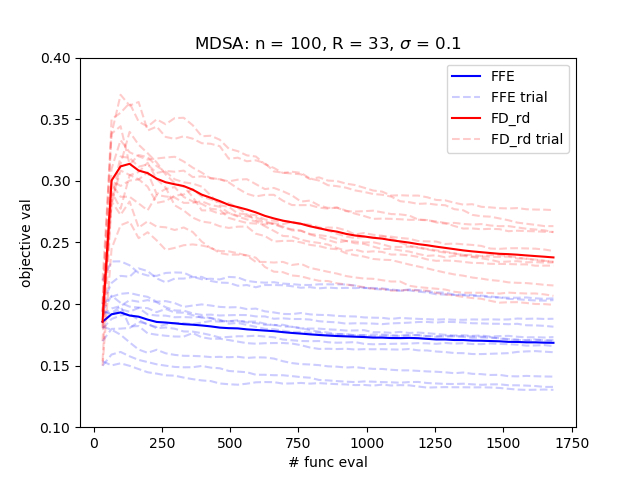}
    \caption{}
  \end{subfigure}
   \begin{subfigure}[b]{0.3\linewidth}
    \includegraphics[width=\linewidth]{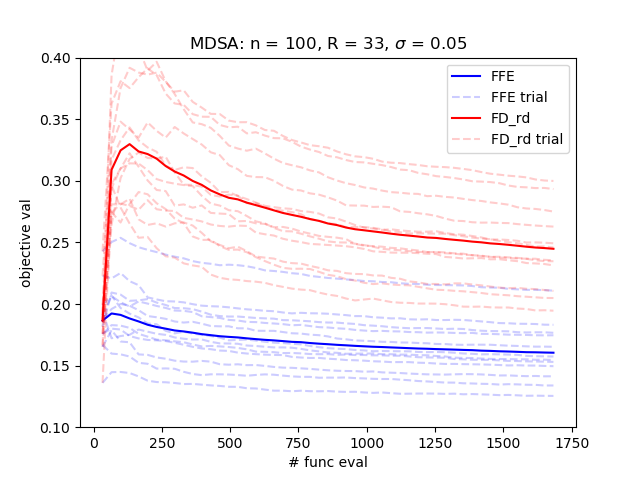}
    \caption{}
  \end{subfigure}
  \begin{subfigure}[b]{0.3\linewidth}
    \includegraphics[width=\linewidth]{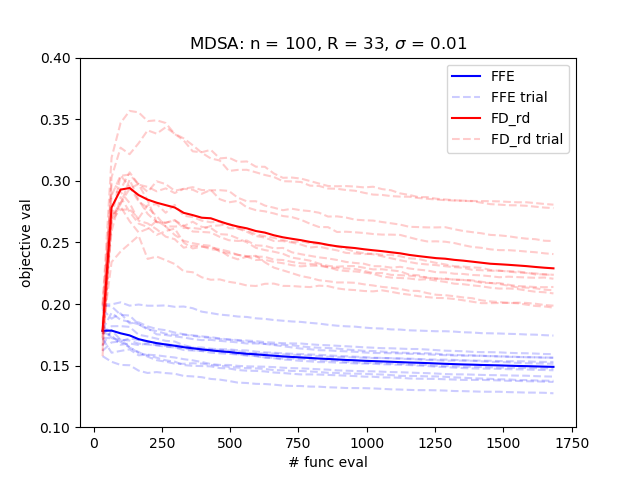}
    \caption{}
  \end{subfigure}
  \begin{subfigure}[b]{0.3\linewidth}
    \includegraphics[width=\linewidth]{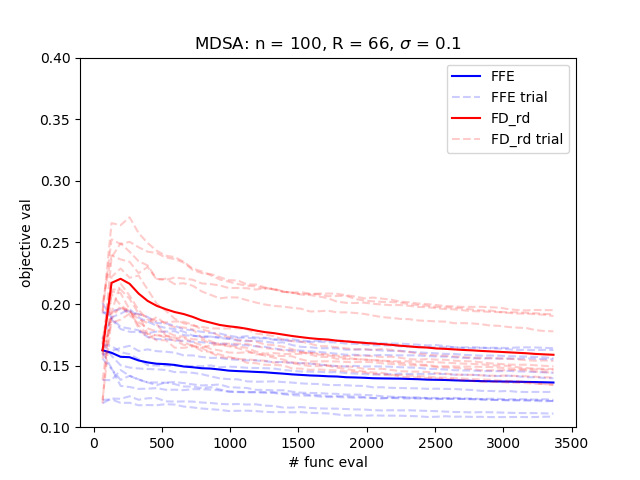}
    \caption{}
  \end{subfigure}
  \begin{subfigure}[b]{0.3\linewidth}
    \includegraphics[width=\linewidth]{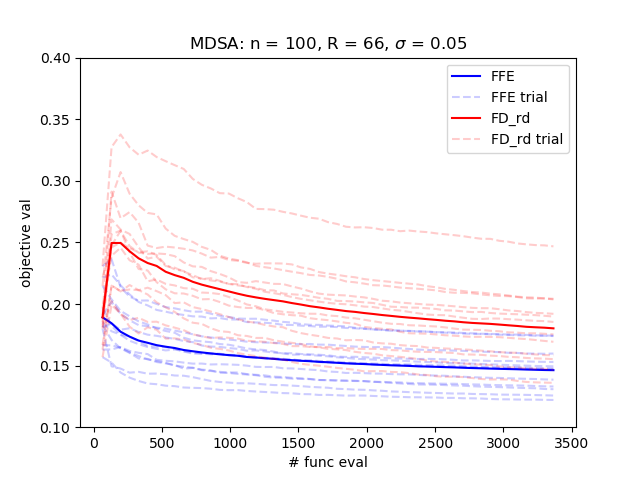}
    \caption{}
  \end{subfigure}
  \begin{subfigure}[b]{0.3\linewidth}
    \includegraphics[width=\linewidth]{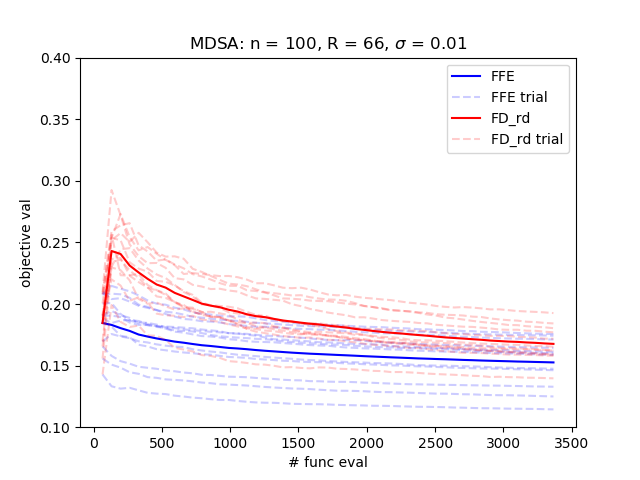}
    \caption{}
  \end{subfigure}
\caption{Performances of MDSA using $\psi_{FFE,\boldsymbol\delta^{**}}$ and $\psi_{FD,random}$ for Rosenbrock objectives, $n=100$, dashed line: 12 trials, solid line: average of 12 trials}
\label{fig:rosen}
\end{figure}

\begin{figure}[h!]
  \centering
  \begin{subfigure}[b]{0.3\linewidth}
    \includegraphics[width=\linewidth]{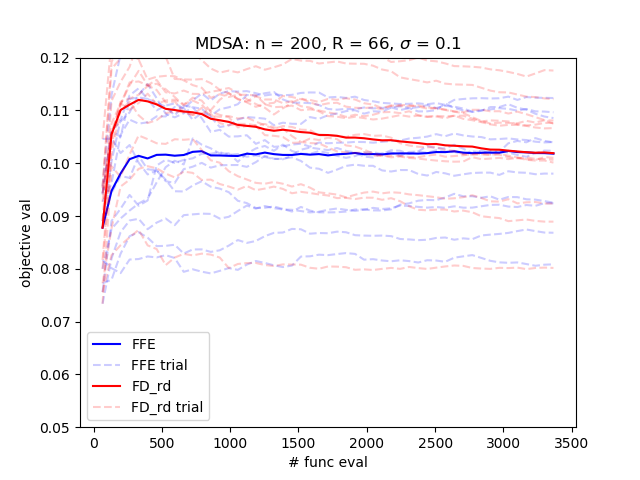}
    \caption{}
  \end{subfigure}
    \begin{subfigure}[b]{0.3\linewidth}
    \includegraphics[width=\linewidth]{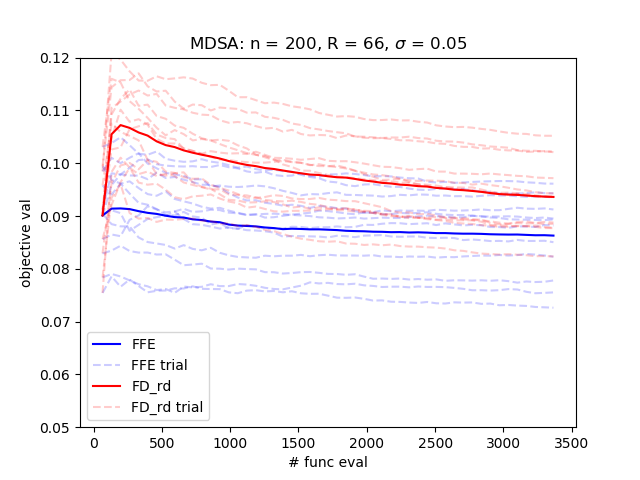}
    \caption{}
  \end{subfigure}
    \begin{subfigure}[b]{0.3\linewidth}
    \includegraphics[width=\linewidth]{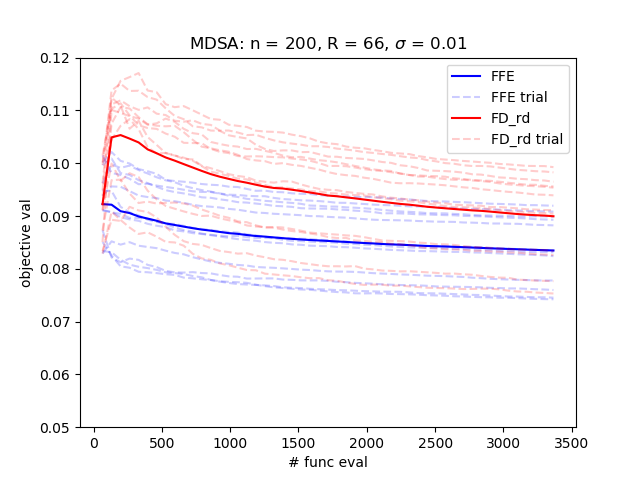}
    \caption{}
  \end{subfigure}
  
    \begin{subfigure}[b]{0.3\linewidth}
    \includegraphics[width=\linewidth]{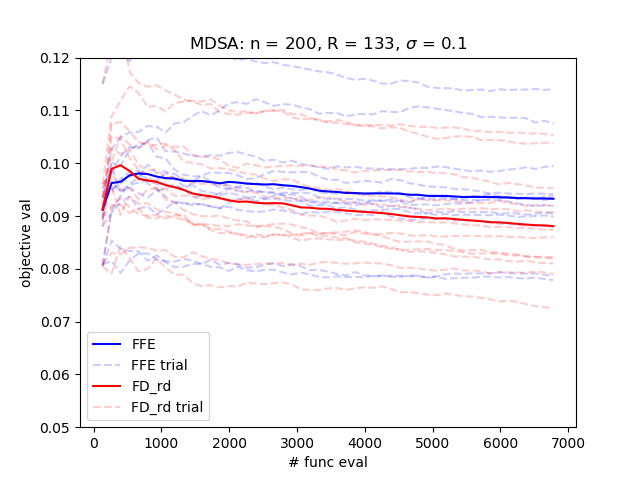}
    \caption{}
  \end{subfigure}
    \begin{subfigure}[b]{0.3\linewidth}
    \includegraphics[width=\linewidth]{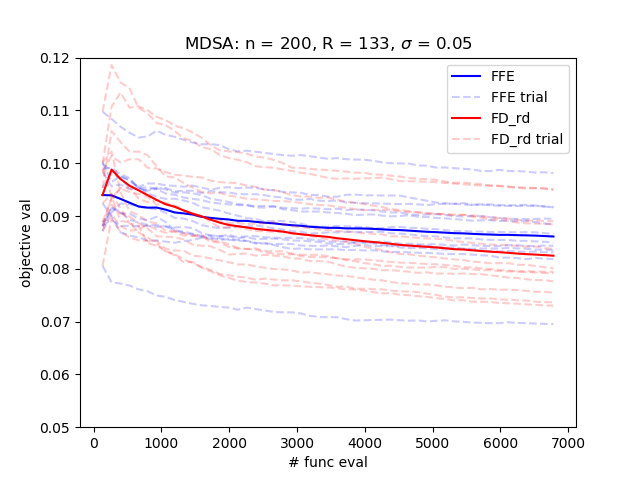}
    \caption{}
  \end{subfigure}
    \begin{subfigure}[b]{0.3\linewidth}
    \includegraphics[width=\linewidth]{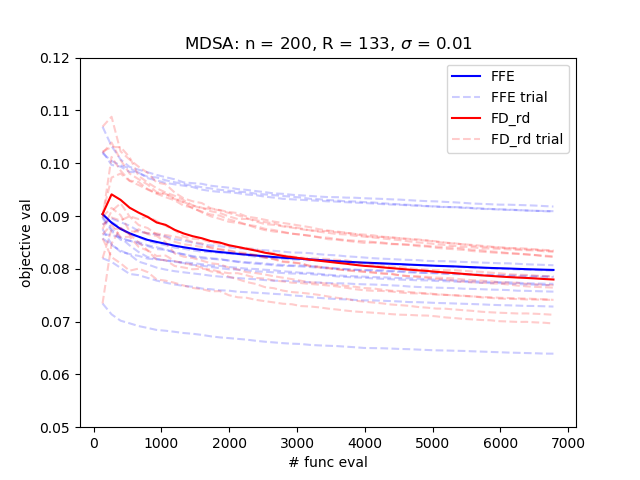}
    \caption{}
  \end{subfigure}
\caption{Performances of MDSA using $\psi_{FFE,\boldsymbol\delta^{**}}$ and $\psi_{FD,random}$ for Rosenbrock objectives, $n=200$, dashed line: 12 trials, solid line: average of 12 trials}
\label{fig:rosen200}
\end{figure}

Figures \ref{fig:rosen40}, \ref{fig:rosen} and \ref{fig:rosen200} show the results for $n = 40, 100, 200$ respectively. We make two observations. The first is on the two gradient estimators: $\psi_{FD,random}$ tends to increase the objectives by a large amount at the first few iterations, and then decrease the objectives. In contrast, $\psi_{FFE,\boldsymbol\delta^{**}}$ also increases the objectives initially but only by a small amount (except in the case for $n=200,\sigma = 0.1$, but in this case, the noise in the observation might be too large compared to the objective function value or the magnitude of the gradient to make meaningful optimization for both gradient estimators). After 50 iterations, the objective value using $\psi_{FFE,\boldsymbol\delta^{**}}$ decreases by a small amount compared with the initial value, and due to the initial increase, $\psi_{FD,random}$ cannot reduce the objective. This shows that $\psi_{FFE,\boldsymbol\delta^{**}}$ does a better job in the minimization of Rosenbrock objective with $\phi$-divergence constraint. 

Our second observation is on the parameters $n,R,\sigma$: From each column in Figures \ref{fig:rosen40}, \ref{fig:rosen} and \ref{fig:rosen200} we see that, for fixed $n$ and $\sigma$, increasing the number of function evaluations $R$ significantly decreases the initial increase of FD, but has little effect on FFE. One possible explanation for the performance of $\psi_{FD, random}$ is that by design, randomly perturbing only $R$ directions and setting all other directions to 0 may result in very large errors for the unperturbed directions. This might significantly deteriorate the accuracy of the proximal updates  compared to updates using the true gradient. 

\subsubsection{$M/G/1$ Queue Objective}
\label{section:exp3}

Our next experimental set considers the service quality of a single-server queue, on which we will apply both FWSA and MDSA and compare between $\psi_{FFE,\boldsymbol\delta^{**}}$ and $\psi_{FD,random}$. The objective is the average waiting time of the first 500 customers in an M/G/1 queue, which can be simulated by Lindley's recursion. We suppose the interarrival time is known and follows an exponential distribution with rate 1. On the other hand, the service time follows a discrete distribution with known support points $x_i = 0.1+1.1\frac{i-1}{n-1},i=1,\ldots,n$. The probability mass function $\mathbf p$ on these support points is the unknown variable (as discussed in Section \ref{sec:SA}, if the distribution is unknown and continuous, we can discretize it first by using a suitable randomized scheme proposed in \cite{ghosh2019robust}).

For the first set of experiment, we run FWSA on a moment constrained problem. The uncertainty set is defined relative to the baseline distribution $\mathbf p_b = \mathbbm 1/n$, given by
\begin{align*}
    \mathcal U = \{\mathbf p \in \mathcal P|&  0.8E_{\mathbf p_b, \mathbf x}[X] \leq E_{\mathbf p, \mathbf x}[X] \leq 1.2E_{\mathbf p_b, \mathbf x}[X],\\
    & 0.8E_{\mathbf p_b, \mathbf x}[X^2] \leq E_{\mathbf p, \mathbf x}[X^2] \leq 1.2E_{\mathbf p_b, \mathbf x}[X^2]\}
\end{align*}
where $E_{\mathbf q,\mathbf x}[f(X)]:= \sum_{i=1}^n q_if(x_i)$. The parameters we choose for FWSA are $a = 0.25$, $b = 0.3$, $\theta = 0.125$, $\beta = 1$, and the initial $\mathbf p_0=\mathbf p_b\in\mathcal U$. We test the algorithms for $n=20, 50, 100$ with corresponding $R = \frac{1}{5}n$. Five trials are carried out for each $n$ and the maximum iteration for each trial is 50. 

We evaluate the performances of the algorithm using two criteria: First is the objective value, approximated by the average of 10 simulation runs. Second is an approximation to the FW optimality gap $-\hat \psi(\mathbf p_k)'(\mathbf q_k - \mathbf p_k)$, where $\mathbf q_k$ is the solution to the $k$-th iteration subproblem.

\begin{figure}[h!]
  \centering
  \begin{subfigure}[b]{0.3\linewidth}
    \includegraphics[width=\linewidth]{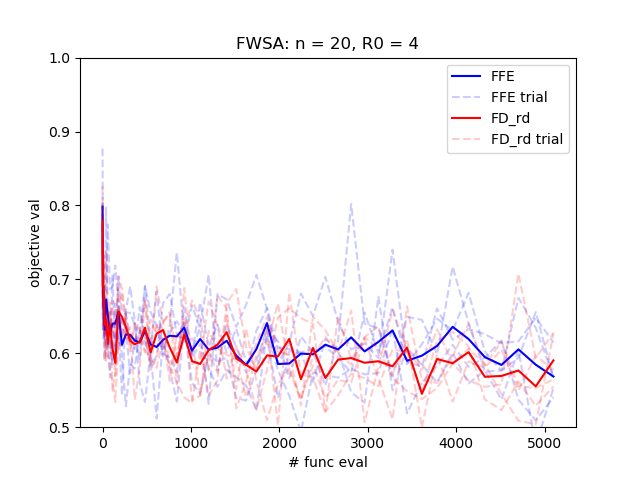}
    \caption{}
  \end{subfigure}
    \begin{subfigure}[b]{0.3\linewidth}
    \includegraphics[width=\linewidth]{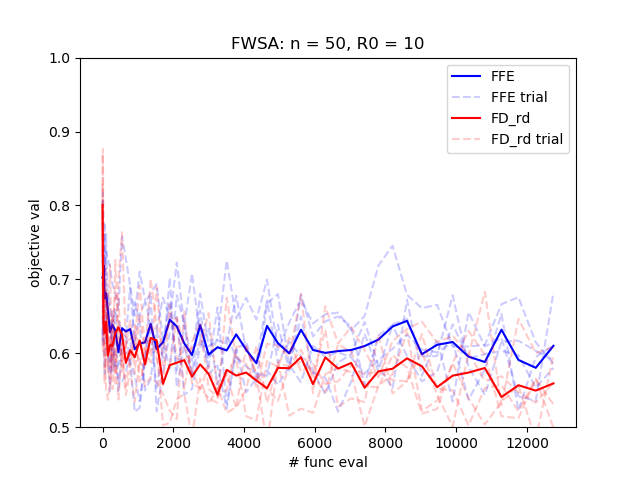}
    \caption{}
  \end{subfigure}
  \begin{subfigure}[b]{0.3\linewidth}
    \includegraphics[width=\linewidth]{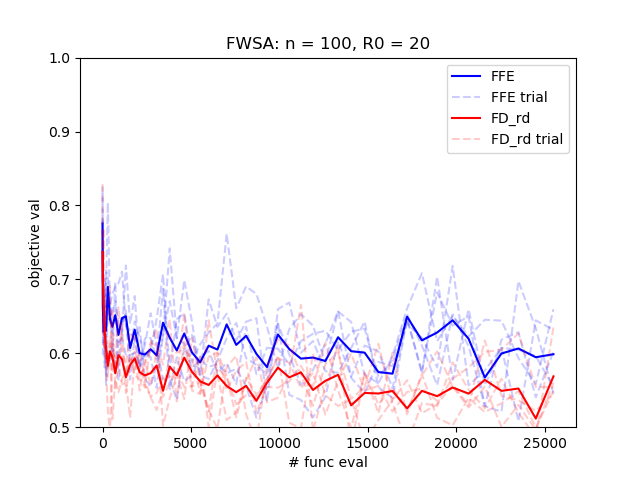}
    \caption{}
  \end{subfigure}
  \begin{subfigure}[b]{0.3\linewidth}
    \includegraphics[width=\linewidth]{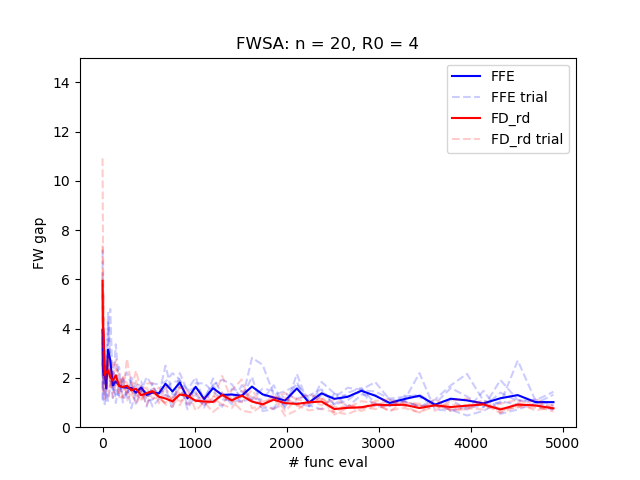}
    \caption{}
  \end{subfigure}
  \begin{subfigure}[b]{0.3\linewidth}
    \includegraphics[width=\linewidth]{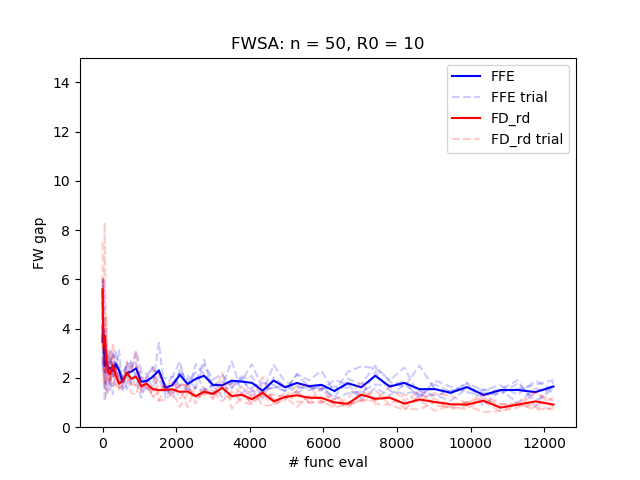}
    \caption{}
  \end{subfigure}
  \begin{subfigure}[b]{0.3\linewidth}
    \includegraphics[width=\linewidth]{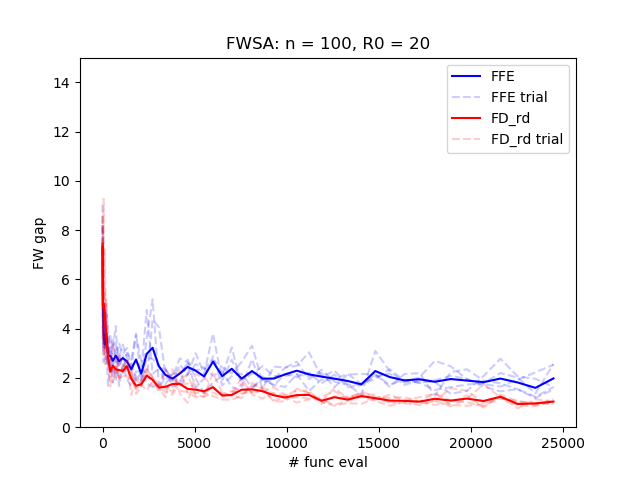}
    \caption{}
  \end{subfigure}
\caption{Performances of FWSA using $\psi_{FFE,\boldsymbol\delta^{**}}$ and $\psi_{FD,random}$ for M/G/1 queue, dashed line: 5 trials, solid line: average of 5 trials}
\label{fig:fwsa_exp}
\end{figure}

Figure \ref{fig:fwsa_exp} shows the trends of the function value and the approximate FW gap against the number of function evaluations for different $n$ and $R_0$. From the first row of Figure \ref{fig:fwsa_exp}, we see that both $\psi_{FFE,\boldsymbol\delta^{**}}$ and $\psi_{FD,random}$ can decrease the objective value, but $\psi_{FD,random}$ tends to converge to smaller objective values. And from the second row, we notice that for both gradient estimators, the FW optimality gap decreases but not down to 0, which could be due to slow convergence and the algorithm has not reached close enough to 0 yet at the termination iteration. 

For the second set of experiment, we run MDSA on a $\phi$-divergence constrained problem, and in particular we choose the Kullback-Leibler divergence. The uncertainty set is defined relative to the baseline distribution $\mathbf p_b = \mathbbm 1/n$ given by
\begin{equation*}
    \mathcal U = \left\{\mathbf p \in \mathcal P\Bigg| d_{\phi}(\mathbf p, \mathbf p_b) = \sum_{i=1}^n p_{i}\log(\frac{p_i}{p_{b,i}}) \leq 0.05\right\}
\end{equation*}
The parameters we choose for MDSA are $a = 0.3$, $\alpha = 1$, $b = 0.3$, $\theta = 0.25$, $\beta = 0$, and the initial $\mathbf p_0=\mathbf p_b\in\mathcal U$. We test the algorithms for $n=20, 50, 100$ with corresponding $R = \frac{1}{5}n$. Five trials are done for each $n$ and the maximum iteration for each trial is 50. 

Similar to FWSA, we evaluate the performances of our algorithm using two criteria. First is the objective value, approximated by the average of 50 simulation runs. Second is $V(\mathbf p_k, \mathbf p_{k+1}) = \sum_{i=1}^n p_{k+1,i}\log(\frac{p_{k+1,i}}{p_{k,i}})$, which is used in the proximal mapping in MDSA and serves as a distance measure between $\mathbf p_k$ and $\mathbf p_{k+1}$.

\begin{figure}[h!]
  \centering
  \begin{subfigure}[b]{0.3\linewidth}
    \includegraphics[width=\linewidth]{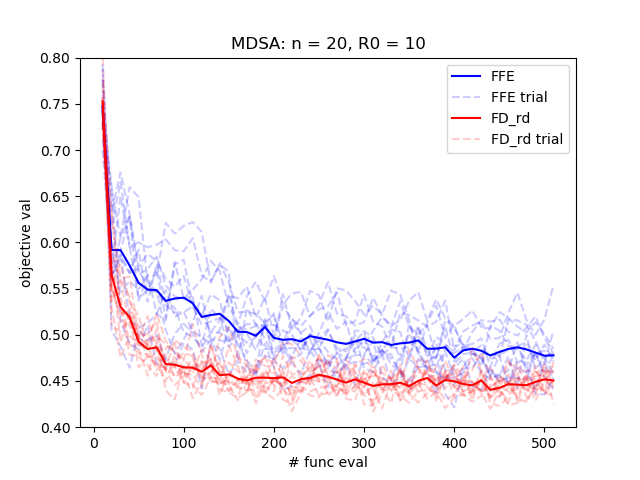}
    \caption{}
  \end{subfigure}
    \begin{subfigure}[b]{0.3\linewidth}
    \includegraphics[width=\linewidth]{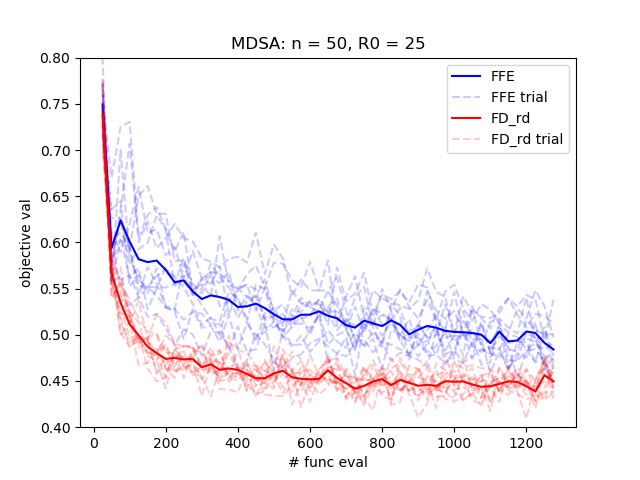}
    \caption{}
  \end{subfigure}
  \begin{subfigure}[b]{0.3\linewidth}
    \includegraphics[width=\linewidth]{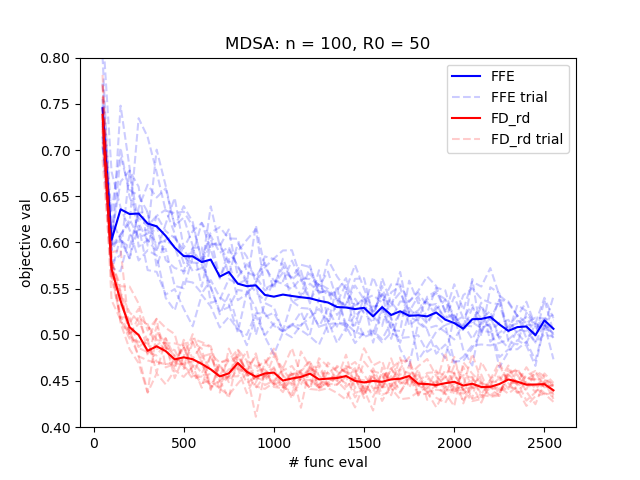}
    \caption{}
  \end{subfigure}
  \begin{subfigure}[b]{0.3\linewidth}
    \includegraphics[width=\linewidth]{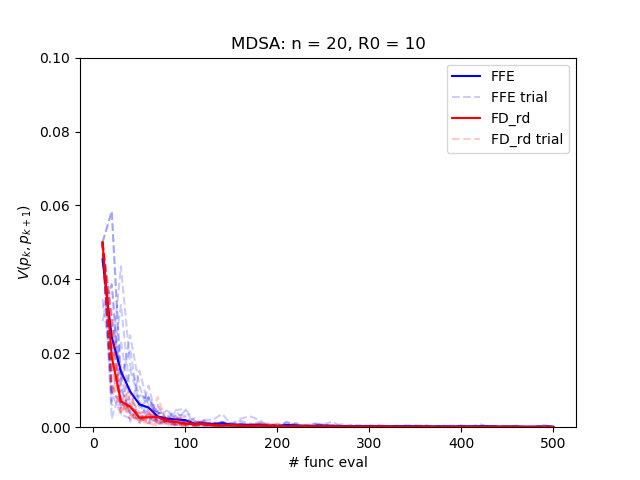}
    \caption{}
  \end{subfigure}
  \begin{subfigure}[b]{0.3\linewidth}
    \includegraphics[width=\linewidth]{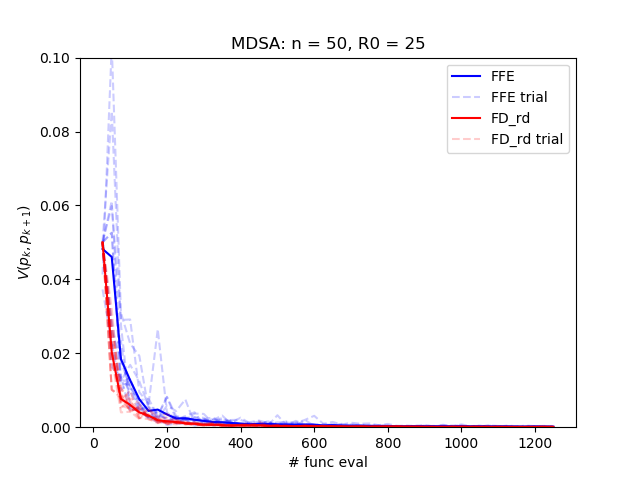}
    \caption{}
  \end{subfigure}
  \begin{subfigure}[b]{0.3\linewidth}
    \includegraphics[width=\linewidth]{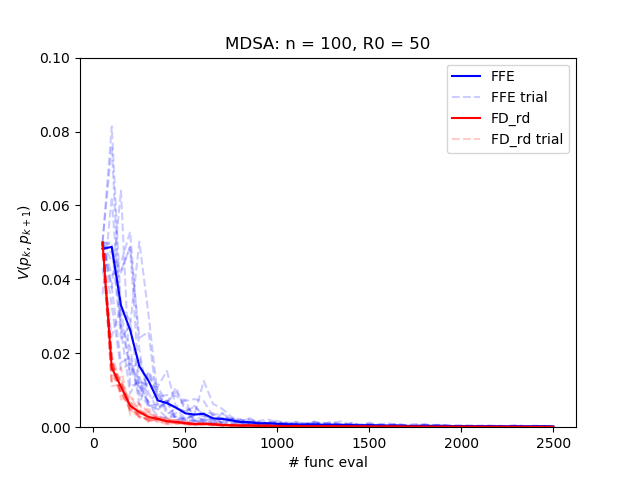}
    \caption{}
  \end{subfigure}
\caption{Performances of MDSA using $\psi_{FFE,\boldsymbol\delta^{**}}$ and $\psi_{FD,random}$ for M/G/1 queue, dashed line: 12 trials, solid line: average of 12 trials}
\label{fig:mdsa_exp}
\end{figure}

Figure \ref{fig:mdsa_exp} shows the trends of the function value and $V(\mathbf p_k,\mathbf p_{k+1})$ against the number of function evaluations for different $n$ and $R_0$. From the first row of Figure \ref{fig:mdsa_exp}, we see that both gradient estimators work very well in minimizing the objective function, but $\psi_{FD,random}$ is faster and can converge to a better solution (lower objective). In the second row, we see that $V(\mathbf p_k, \mathbf p_{k+1})$ converges to 0 as expected, and $\psi_{FD,random}$ converges faster. 

\subsection{Summary and Recommendations}\label{sec:summary}
First, from the experiments in Section \ref{section:exp1}, we see that in terms of bias, $\boldsymbol\delta^*$ is more advantageous than $\boldsymbol\delta^{**}$ since the former's third moments attain a constant value while the latter does not, and so the biases of SFE, FFE, and CFE with $\boldsymbol \delta^*$ are $O(c^2)$, better than the $O(c)$ biases with $\boldsymbol \delta^{**}$. However, in terms of variance, $\boldsymbol\delta^{**}$ is smaller than $\boldsymbol\delta^*$ as $\boldsymbol\delta^{**}$ is designed to have a smaller $\gamma$, and this benefit appears in our experiments to outweigh the impact from the bias. In practice, $\psi_{FFE,\boldsymbol\delta^{**}}$ appears to be the most preferable choice among the six combinations formed by the three gradient estimators and the two Dirichlet mixtures. 

From the experiments in Section \ref{sec:opt num}, we see that $\psi_{FFE,\boldsymbol\delta^{**}}$ outperforms $\psi_{FD, random}$ in the Rosenbrock objective, while $\psi_{FD, random}$ does a better job in the queueing model. In addition, we see that in the queuing problem, MDSA has a better convergence behavior than FWSA, as $V(\mathbf p_k,\mathbf p_{k+1})$ converges to 0 for MDSA while the FW gap does not converge to 0 for FWSA. Thus, we believe that a combination of MDSA and $\psi_{FFE,\boldsymbol\delta^{**}}$ or $\psi_{FD,random}$ works best in practice. And in terms of choosing between $\psi_{FFE,\boldsymbol\delta^{**}}$ and $\psi_{FD,random}$, we expect that when the objective value is sensitive to changes in decision variables, $\psi_{FFE,\boldsymbol\delta^{**}}$ performs better, and when the objective is insensitive, $\psi_{FD,random}$ is preferable. 

We provide some intuition behind the last recommendation. By design,
\begin{equation*}
    (\psi_{FD,random})_i = \frac{n}{R} \sum_{j=1}^{N_i} \frac{\hat{Z}^{2j-1}(\mathbf p +c(\mathbf e_{i} - \mathbf p))-\hat{Z}^{2j}(\mathbf p)}{c}\mathbf e_{i}
\end{equation*}
where $N_i$ is the number of times the $i$-th direction gets perturbed, and for those directions that are never perturbed, we set $(\psi_{FD,random})_i = 0$. However, it is possible that setting the latter to 0 is not too problematic. In particular, if the gradient $\nabla Z_i(\mathbf p) \approx 0$, then setting them to 0 is in fact a good choice. We suspect that this is the case in the queuing model: The support points are very close to each other, so that the distribution described by the probability weights $\mathbf p$ might be insensitive to some directions of change. On the other hand, for Rosenbrock, it is not the case that many components of the gradient are 0, so that missing even one direction can make the result of the proximal update step meaningless (if not harmful). This explains the bad performances of $\psi_{FD,random}$. In contrast, the main idea of $\psi_{FFE}$ is to use $S(\mathbf p, \mathbf \boldsymbol\delta)(\mathbf \boldsymbol\delta - \mathbf p)' \nabla Z$ to approximate the gradient (where $S(\mathbf p, \mathbf \boldsymbol\delta) = \gamma (\mathbf \boldsymbol\delta - \mathbf p)$). For fixed $\lVert \mathbf \boldsymbol\delta - \mathbf p \rVert$ and $\lVert \nabla Z\rVert$, $\lvert (\mathbf \boldsymbol\delta - \mathbf p)' \nabla Z\rvert$ measures the ``angle" between $(\mathbf \boldsymbol\delta - \mathbf p)$ and $\nabla Z$. The smaller this angle is, the closer $(\mathbf \boldsymbol\delta - \mathbf p)$ is to $\nabla Z$, and the more ``weights" we put on $S(\mathbf p, \mathbf \boldsymbol\delta)$. Thus, in the case where $S(\mathbf p, \mathbf \boldsymbol\delta)$ does not align well with $\nabla Z$, $\psi_{FFE}$ is smaller in magnitude and the proximal update step will automatically put more penalty on the $V(\cdot,\cdot)$ term, making sure that the update will not be hurt by the poorly estimated gradient by too much.

\ACKNOWLEDGMENT{We gratefully acknowledge support from the National Science Foundation under grants CAREER CMMI-1834710 and IIS-1849280. 
}

\bibliographystyle{informs2014} 
\bibliography{bibliography} 

%
%
%





\ECSwitch


\ECHead{Proofs of Statements}


\section{Proofs for Section \ref{sec:problem}}



\proof{Proof of Lemma \ref{lemma:taylor}.}
Define function $f:[0,1] \to \mathbb R$
\begin{equation*}
    f(t) = Z(\mathbf p_1 + t(\mathbf p_2-\mathbf p_1))
\end{equation*}
Then, applying the standard Taylor expansion to $f(t)$ around $t=0$, we get the desired result.\hfill\Halmos
\endproof

\section{Proofs for Section \ref{sec:analysis}}
\proof{Proof of Theorem \ref{thm:SFE}.}
For the expectation, 
\begin{align}
\psi_{SFE} &= E[\hat \psi_{SFE}] = E[\frac{1}{R}\sum_{i=1}^{R}E[\frac{\hat Z^{i}((1-c)\mathbf p+c\boldsymbol\delta^i)}{c}S(\mathbf p,\boldsymbol\delta ^i)|\boldsymbol\delta^i]]\nonumber \\
&=E[\frac{1}{R}\sum_{i=1}^{R}\frac{Z((1-c)\mathbf p+c\boldsymbol\delta^i)}{c}S(\mathbf p,\boldsymbol\delta ^i)]\nonumber \\
&=E[\gamma \frac{Z((1-c)\mathbf p+c\boldsymbol\delta)}{c}(\boldsymbol\delta -\mathbf p)]~~~\text{noting that $S(\mathbf p,\boldsymbol\delta)=\gamma(\boldsymbol\delta -\mathbf p)$}
\label{eq:sfe_bias1}
\end{align}

By Lemma \ref{lemma:taylor}, we can write the Taylor expansion of $Z(\mathbf p+c (\boldsymbol\delta-\mathbf p))$ as 
\begin{equation}
    Z(\mathbf p+c (\boldsymbol\delta-\mathbf p)) = Z(\mathbf p)+c\nabla Z(\mathbf p)'(\boldsymbol\delta-\mathbf p)+\frac{c^2}{2}(\boldsymbol\delta-\mathbf p)'\nabla^2Z(\mathbf p+\eta c (\boldsymbol\delta-\mathbf p))(\boldsymbol\delta-\mathbf p)
\label{eq:sfe_taylor}
\end{equation}
for some $\eta\in[0,1]$ that depends on $\mathbf p, \boldsymbol \delta$. Combining \eqref{eq:sfe_bias1} and \eqref{eq:sfe_taylor}, we get
\begin{align}
\psi_{SFE}&=\gamma Z(\mathbf p)E[ \boldsymbol\delta-\mathbf p]/c+ \gamma E[\nabla Z(\mathbf p)'(\boldsymbol\delta-\mathbf p)(\boldsymbol\delta-\mathbf p)]+\frac{\gamma c}{2}E[(\boldsymbol\delta - \mathbf p)'\nabla^2Z(\mathbf p+c\eta(\boldsymbol\delta - \mathbf p))(\boldsymbol\delta - \mathbf p)(\boldsymbol\delta - \mathbf p)]\nonumber \\
&=0\cdot\gamma Z(\mathbf p) + (\mathbf I-\mathbbm 1\mathbbm 1'/n) \nabla Z(\mathbf p)+\frac{\gamma c}{2}E[(\boldsymbol\delta - \mathbf p)'\nabla^2Z(\mathbf p)(\boldsymbol\delta - \mathbf p)(\boldsymbol\delta - \mathbf p)]+D\nonumber ~~~\text{by \eqref{first order} and \eqref{second order}}\\
&=\nabla Z(\mathbf p)-\frac{\mathbbm 1'\nabla Z(\mathbf p)}{n}\mathbbm 1+\frac{\gamma c}{2}E[(\boldsymbol\delta - \mathbf p)'\nabla^2Z(\mathbf p)(\boldsymbol\delta - \mathbf p)(\boldsymbol\delta - \mathbf p)]+D
\label{eq:sfe_bias2}
\end{align}
where $D = \frac{\gamma c}{2}E[(\boldsymbol\delta - \mathbf p)'(\nabla^2Z(\mathbf p+c\eta(\boldsymbol\delta - \mathbf p))-\nabla^2Z(\mathbf p))(\boldsymbol\delta - \mathbf p)(\boldsymbol\delta - \mathbf p)]$.

The first and second terms in \eqref{eq:sfe_bias2} give the correct gradient up to a translation of $\mathbbm 1$. For the third term, if only \eqref{first order} and \eqref{second order} are satisfied,
\begin{align}
\lVert \frac{\gamma c}{2}E[(\boldsymbol\delta - \mathbf p)'\nabla^2Z(\mathbf p)(\boldsymbol\delta - \mathbf p)(\boldsymbol\delta - \mathbf p)] \rVert&\leq \frac{\gamma c}{2} E[\lVert(\boldsymbol\delta - \mathbf p)'\nabla^2Z(\mathbf p)(\boldsymbol\delta - \mathbf p)(\boldsymbol\delta - \mathbf p)\rVert] \nonumber \\
& \leq \frac{\gamma c}{2} \lVert \nabla^2Z(\mathbf p)\rVert E[\lVert\boldsymbol\delta - \mathbf p\rVert^3] \nonumber \\
& \leq (n-1)c \lVert \nabla^2Z(\mathbf p)\rVert\nonumber\\
& \leq (n-1)cM_2\nonumber
\end{align}
where $M_2 = \sup_{\mathbf p\in \mathcal P} \lVert \nabla^2 Z(\mathbf p) \rVert$. (The sup is finite because $\nabla^2 Z$ is $L_2$-Lipschitz and $\mathcal P$ is bounded.) The second to last inequality holds because $\lVert\boldsymbol\delta - \mathbf p\rVert \leq 2$, and by \eqref{second order}, $\gamma E[\lVert\boldsymbol\delta - \mathbf p\rVert^2]=n-1$. 

If in addition to \eqref{first order} and \eqref{second order}, \eqref{third order} is also true, then the third term in \eqref{eq:sfe_bias2} becomes a translation of $\mathbbm 1$: 
\begin{equation*}
\frac{\gamma c}{2}E[(\boldsymbol\delta - \mathbf p)'\nabla^2Z(\mathbf p)(\boldsymbol\delta - \mathbf p)(\boldsymbol\delta - \mathbf p)] = \frac{c\mu}{2}\mathbbm 1'\nabla^2 Z(\mathbf p) \mathbbm 1 \mathbbm 1
\end{equation*}

For the remainder term $D$, we have
\begin{align}
\lVert D \rVert &\leq \frac{\gamma c}{2}E[\lVert (\boldsymbol\delta - \mathbf p)'(\nabla^2Z(\mathbf p+c\eta(\boldsymbol\delta - \mathbf p))-\nabla^2Z(\mathbf p))(\boldsymbol\delta - \mathbf p)(\boldsymbol\delta - \mathbf p)\rVert]\nonumber \\
&\leq \frac{\gamma c}{2}E[\lVert (\boldsymbol\delta - \mathbf p)\rVert^3 \lVert \nabla^2Z(\mathbf p+c\eta(\boldsymbol\delta - \mathbf p))-\nabla^2Z(\mathbf p)\rVert ]\nonumber \\
& \leq \frac{\gamma c^2L_2}{2}E[\lVert \boldsymbol\delta - \mathbf p\rVert^4]~~~~~\text{because $\nabla^2 Z$ is $L_2$-Lipschitz}\nonumber \\
& \leq 2\gamma c^2L_2E[\lVert \boldsymbol\delta - \mathbf p\rVert^2]~~~~~\text{because $\lVert \boldsymbol\delta - \mathbf p\rVert \leq 2$}\nonumber \\
& = 2(n-1)c^2L_2~~~~~\text{by \eqref{second order}}\nonumber
\end{align}

Thus, 
\begin{equation*}
\psi_{SFE} = \nabla Z(\mathbf p) +\epsilon_0 \mathbbm 1 +\mathbf \epsilon
\end{equation*}
where $(\epsilon_0, \lVert \epsilon \rVert) = \begin{cases}(-\nabla Z(\mathbf p)'\mathbbm 1/n,O(ncM_2+nc^2L_2))~~~\text{if \eqref{first order} and \eqref{second order} are true}\\
(\frac{c\mu}{2}\mathbbm 1'\nabla Z(\mathbf p)\mathbbm 1-\nabla Z(\mathbf p)'\mathbbm 1/n,O(nc^2L_2))~~~\text{if \eqref{third order} is also true}\end{cases}$. 

For the variance, we use the decomposition 
\begin{equation}
Var[\hat \psi_{SFE}]=E[Var[\hat \psi_{SFE}|\boldsymbol\delta^1,...,\boldsymbol\delta^{R}]]+Var[E[\hat \psi_{SFE}|\boldsymbol\delta^1,...,\boldsymbol\delta^{R}]]
\label{sfe_var0}
\end{equation}

For the first term,
\begin{align}
tr[E[Var[\hat \psi_{SFE}|\boldsymbol\delta^1,...,\boldsymbol\delta^{R}]]]&=tr[E[\frac{1}{R^2}\sum_{i=1}^{R}Var[\frac{\hat Z^{i}((1-c)\mathbf p +c\boldsymbol\delta^i)}{c}S(\mathbf p,\boldsymbol\delta^i)|\boldsymbol\delta^i]]]\nonumber \\
&=E[tr[\frac{\sigma^2((1-c)\mathbf p +c\boldsymbol\delta)}{Rc^2}S(\mathbf p,\boldsymbol\delta)S(\mathbf p,\boldsymbol\delta)']]\nonumber \\
&\leq \frac{\sigma^2}{Rc^2}E[\lVert S(\mathbf p,\boldsymbol\delta) \rVert^2]~~~~~\text{by Assumption \ref{assumption:smooth_noise}}\nonumber\\
&\leq \frac{(n-1)\gamma \sigma^2}{Rc^2}~~~~~\text{because $S(\mathbf p,\boldsymbol\delta) = \gamma (\boldsymbol \delta - \mathbf p)$ and by \eqref{second order}}
\label{sfe_var1}
\end{align}

For the second term,
\begin{align}
tr[Var[E[\hat \psi_{SFE}|\boldsymbol\delta^1,...,\boldsymbol\delta^{R}]]]&=\frac{1}{R}tr[Var[\frac{Z((1-c)\mathbf p+c\boldsymbol\delta)}{c}S(\mathbf p,\boldsymbol\delta)]]\nonumber \\
&\leq \frac{M_0^2}{Rc^2} E[ \lVert S(\mathbf p,\boldsymbol\delta) \rVert^2]\nonumber \\
& = \frac{\gamma M_0^2 (n-1)}{Rc^2}~~~~~\text{because $S(\mathbf p,\boldsymbol\delta) = \gamma (\boldsymbol \delta - \mathbf p)$ and by \eqref{second order}}
\label{sfe_var2}
\end{align}
where $M_0 = \sup_{\mathbf p\in \mathcal P} \lvert Z(\mathbf p) \rvert$. (The sup is finite because $Z$ is $L_0$-Lipschitz and $\mathcal P$ is bounded.)

Combining \eqref{sfe_var1} and \eqref{sfe_var2} with \eqref{sfe_var0}, we get
\begin{equation*}
tr[Var[\hat{\psi}_{SFE}]] \leq \frac{\gamma (n-1) }{Rc^2}[\sigma^2+M_0^2]
\end{equation*}\hfill\Halmos
\endproof

\proof{Proof of Theorem \ref{thm_necessary_conditions}.}
By Lemma \ref{lemma:taylor}, we get 
\begin{align}
\psi_{SFE} &= E[\frac{Z(\mathbf p+c (\boldsymbol\delta-\mathbf p))}{c}S(\mathbf p,\boldsymbol\delta)]\nonumber \\
&=\frac{1}{c} E[(Z(\mathbf p)+c \nabla Z(\mathbf p)'(\boldsymbol\delta-\mathbf p)+\frac{c^2}{2}(\boldsymbol\delta - \mathbf p)'\nabla^2Z(\mathbf p)(\boldsymbol\delta - \mathbf p))S(\mathbf p,\boldsymbol\delta)]+D\nonumber \\
&=\frac{ Z(\mathbf p)}{c}E[S(\mathbf p,\boldsymbol\delta)]+ E[S(\mathbf p,\boldsymbol\delta)(\boldsymbol\delta - \mathbf p)']\nabla Z(\mathbf p)+\frac{ c}{2}E[(\boldsymbol\delta - \mathbf p)'\nabla^2Z(\mathbf p)(\boldsymbol\delta - \mathbf p)S(\mathbf p,\boldsymbol\delta)]+D
\label{eq:thm2_taylor}
\end{align}
where $D = \frac{\gamma c}{2}E[(\boldsymbol\delta - \mathbf p)'(\nabla^2Z(\mathbf p+c\eta(\boldsymbol\delta - \mathbf p))-\nabla^2Z(\mathbf p))(\boldsymbol\delta - \mathbf p)(\boldsymbol\delta - \mathbf p)]$.

The $\lVert D\rVert$ in \eqref{eq:thm2_taylor} is $O(c^2)$ because
\begin{align}
    \lVert D \rVert &\leq \frac{\gamma c}{2}E[\lVert (\boldsymbol\delta - \mathbf p)'(\nabla^2Z(\mathbf p+c\eta(\boldsymbol\delta - \mathbf p))-\nabla^2Z(\mathbf p))(\boldsymbol\delta - \mathbf p)(\boldsymbol\delta - \mathbf p)\rVert]\nonumber\\
    &\leq \frac{\gamma c}{2}E[\lVert \nabla^2Z(\mathbf p+c\eta(\boldsymbol\delta - \mathbf p))-\nabla^2Z(\mathbf p)\rVert\cdot \lVert \boldsymbol\delta - \mathbf p\rVert^3]\nonumber\\
    &\leq \frac{\gamma c^2L_2}{2} E[\lvert \boldsymbol\delta -\mathbf p \rVert^4]~~~~~\text{because $\nabla^2 Z$ is $L_2$-Lipschitz}\nonumber\\
    &\leq 8\gamma c^2 L_2=O(c^2)~~~~~\text{because $\lVert \boldsymbol\delta -\mathbf p \rVert\leq 2$ and by \eqref{second order}}\nonumber
\end{align}

For the first term, we need $E[S(\mathbf p,\boldsymbol\delta)] =\gamma E[\boldsymbol\delta - \mathbf p]= \epsilon_1\mathbbm 1$ for some $\epsilon_1\in \mathbbm R$. Noting that $\mathbbm 1'(\gamma E[\boldsymbol\delta - \mathbf p]) = \gamma (1-1) = 0 =n\epsilon_1$, we have $\epsilon_1 = 0$ and so $E[S(\mathbf p,\boldsymbol\delta)] = 0$.

For the second term, we have $E[S(\mathbf p,\boldsymbol\delta)(\boldsymbol\delta - \mathbf p)'] = \gamma E[(\boldsymbol\delta - \mathbf p)(\boldsymbol\delta - \mathbf p)'] = \mathbf I-\mathbbm 1\mathbf v'$ for any $\mathbf v\in \mathbbm R^n$. However, $(\boldsymbol\delta - \mathbf p)(\boldsymbol\delta - \mathbf p)'$ is symmetric, so $\mathbf v = \epsilon_2\mathbbm 1$ for some $\epsilon_2 \in \mathbbm R$. And, again, we have $\mathbbm 1'(\boldsymbol\delta - \mathbf p)(\boldsymbol\delta - \mathbf p)'\mathbbm 1 = 0$, so that $\epsilon_2 = \frac{1}{n}$, and we get the second condition.

For the third term, its $k$-th component is
\begin{equation}
\gamma E[\sum_{i,j}\nabla^2Z(\mathbf p)_{i,j}(\boldsymbol\delta - \mathbf p)_i(\boldsymbol\delta - \mathbf p)_j(\boldsymbol\delta - \mathbf p)_k]=\gamma \sum_{i,j}\nabla^2Z(\mathbf p)_{i,j}E[(\boldsymbol\delta - \mathbf p)_i(\boldsymbol\delta - \mathbf p)_j(\boldsymbol\delta - \mathbf p)_k]
\label{eq_temp}
\end{equation}
and the requirement for the cancellation of this term up to translation of $\mathbbm 1$ is that \eqref{eq_temp} is the same for all $k\in[n]$. In addition, we want this to be true for all possible Hessian $\nabla^2 Z$ (i.e. all possible symmetric matrices). In particular, let's define $\mathbf E_{ab} = \mathbf e_a\mathbf e_b'$. 
We consider when $\nabla^2 Z = (\mathbf E_{ab}+\mathbf E_{ba})/2$,
\begin{equation*}
\gamma E[\sum_{i,j}\nabla^2Z(\mathbf p)_{i,j}(\boldsymbol\delta - \mathbf p)_i(\boldsymbol\delta - \mathbf p)_j(\boldsymbol\delta - \mathbf p)_k]=\gamma E[(\boldsymbol\delta - \mathbf p)_a(\boldsymbol\delta - \mathbf p)_b(\boldsymbol\delta - \mathbf p)_k]
\end{equation*}
is the same for all possible $a,b,k$, thus leading to \eqref{third order}.\hfill\Halmos
\endproof

\proof{Proof of Theorem \ref{thm:FFE}.}
For the expectation, 
\begin{align}
\psi_{FFE} &= E[\hat \psi_{FFE}] = E[\frac{1}{R}\sum_{i=1}^{R}E[\frac{\hat Z^{2i-1}((1-c)\mathbf p+c\boldsymbol\delta^i)-\hat Z^{2i}(\mathbf p)}{c}S(\mathbf p,\boldsymbol\delta ^i)|\boldsymbol\delta^i]]\nonumber \\
&=E[\frac{1}{R}\sum_{i=1}^{R}\frac{Z((1-c)\mathbf p+c\boldsymbol\delta^i)-Z(\mathbf p)}{c}S(\mathbf p,\boldsymbol\delta ^i)]\nonumber \\
&=E[\gamma \frac{Z((1-c)\mathbf p+c\boldsymbol\delta)-Z(\mathbf p)}{c}(\boldsymbol\delta -\mathbf p)]~~~\text{noting that $S(\mathbf p,\boldsymbol\delta)=\gamma(\boldsymbol\delta -\mathbf p)$}
\label{eq:ffe_bias1}
\end{align}

By Lemma \ref{lemma:taylor}, we can write $Z(\mathbf p+c (\boldsymbol\delta-\mathbf p))-Z(\mathbf p)$ as 
\begin{equation}
    Z(\mathbf p+c (\boldsymbol\delta-\mathbf p)) -Z(\mathbf p)= c\nabla Z(\mathbf p)'(\boldsymbol\delta-\mathbf p)+\frac{c^2}{2}(\boldsymbol\delta-\mathbf p)'\nabla^2Z(\mathbf p+\eta c (\boldsymbol\delta-\mathbf p))(\boldsymbol\delta-\mathbf p)
\label{eq:ffe_taylor}
\end{equation}
for some $\eta\in[0,1]$ that depends on $\mathbf p, \boldsymbol \delta$. Combining \eqref{eq:ffe_bias1} and \eqref{eq:ffe_taylor}, we get
\begin{align}
\psi_{FFE}&= \gamma E[\nabla Z(\mathbf p)'(\boldsymbol\delta-\mathbf p)(\boldsymbol\delta-\mathbf p)]+\frac{\gamma c}{2}E[(\boldsymbol\delta - \mathbf p)'\nabla^2Z(\mathbf p+c\eta(\boldsymbol\delta - \mathbf p))(\boldsymbol\delta - \mathbf p)(\boldsymbol\delta - \mathbf p)]\nonumber \\
&=(\mathbf I-\mathbbm 1\mathbbm 1'/n) \nabla Z(\mathbf p)+\frac{\gamma c}{2}E[(\boldsymbol\delta - \mathbf p)'\nabla^2Z(\mathbf p)(\boldsymbol\delta - \mathbf p)(\boldsymbol\delta - \mathbf p)]+D\nonumber ~~~\text{by \eqref{first order} and \eqref{second order}}\\
&=\nabla Z(\mathbf p)-\frac{\mathbbm 1'\nabla Z(\mathbf p)}{n}\mathbbm 1+\frac{\gamma c}{2}E[(\boldsymbol\delta - \mathbf p)'\nabla^2Z(\mathbf p)(\boldsymbol\delta - \mathbf p)(\boldsymbol\delta - \mathbf p)]+D
\label{eq:ffe_bias2}
\end{align}
where $D = \frac{\gamma c}{2}E[(\boldsymbol\delta - \mathbf p)'(\nabla^2Z(\mathbf p+c\eta(\boldsymbol\delta - \mathbf p))-\nabla^2Z(\mathbf p))(\boldsymbol\delta - \mathbf p)(\boldsymbol\delta - \mathbf p)]$.

The first and second terms in \eqref{eq:ffe_bias2} give the correct gradient up to a translation of $\mathbbm 1$. For the third term, if only \eqref{first order} and \eqref{second order} are satisfied,
\begin{align}
\lVert \frac{\gamma c}{2}E[(\boldsymbol\delta - \mathbf p)'\nabla^2Z(\mathbf p)(\boldsymbol\delta - \mathbf p)(\boldsymbol\delta - \mathbf p)] \rVert&\leq \frac{\gamma c}{2} E[\lVert(\boldsymbol\delta - \mathbf p)'\nabla^2Z(\mathbf p)(\boldsymbol\delta - \mathbf p)(\boldsymbol\delta - \mathbf p)\rVert] \nonumber \\
& \leq \frac{\gamma c}{2} \lVert \nabla^2Z(\mathbf p)\rVert E[\lVert\boldsymbol\delta - \mathbf p\rVert^3] \nonumber \\
& \leq (n-1)c \lVert \nabla^2Z(\mathbf p)\rVert\nonumber\\
& \leq (n-1)cM_2\nonumber
\end{align}
where $M_2 = \sup_{\mathbf p\in \mathcal P} \lVert \nabla^2 Z(\mathbf p) \rVert$. (The sup is finite because $\nabla^2 Z$ is $L_2$-Lipschitz and $\mathcal P$ is bounded.) The second to the last inequality is because $\lVert\boldsymbol\delta - \mathbf p\rVert \leq 2$, and by \eqref{second order}, $\gamma E[\lVert\boldsymbol\delta - \mathbf p\rVert^2]=n-1$. 

If in addition to \eqref{first order} and \eqref{second order}, \eqref{third order} is also true, then the third term in \eqref{eq:ffe_bias2} becomes a translation of $\mathbbm 1$: 
\begin{equation*}
\frac{\gamma c}{2}E[(\boldsymbol\delta - \mathbf p)'\nabla^2Z(\mathbf p)(\boldsymbol\delta - \mathbf p)(\boldsymbol\delta - \mathbf p)] = \frac{c\mu}{2}\mathbbm 1'\nabla^2 Z(\mathbf p) \mathbbm 1 \mathbbm 1
\end{equation*}

For the remainder term $D$, we have
\begin{align}
\lVert D \rVert &\leq \frac{\gamma c}{2}E[\lVert (\boldsymbol\delta - \mathbf p)'(\nabla^2Z(\mathbf p+c\eta(\boldsymbol\delta - \mathbf p))-\nabla^2Z(\mathbf p))(\boldsymbol\delta - \mathbf p)(\boldsymbol\delta - \mathbf p)\rVert]\nonumber \\
&\leq \frac{\gamma c}{2}E[\lVert (\boldsymbol\delta - \mathbf p)\rVert^3 \lVert \nabla^2Z(\mathbf p+c\eta(\boldsymbol\delta - \mathbf p))-\nabla^2Z(\mathbf p)\rVert ]\nonumber \\
& \leq \frac{\gamma c^2L_2}{2}E[\lVert \boldsymbol\delta - \mathbf p\rVert^4]~~~~~\text{because $\nabla^2 Z$ is $L_2$-Lipschitz}\nonumber \\
& \leq 2\gamma c^2L_2E[\lVert \boldsymbol\delta - \mathbf p\rVert^2]~~~~~\text{because $\lVert \boldsymbol\delta - \mathbf p\rVert \leq 2$}\nonumber \\
& = 2(n-1)c^2L_2~~~~~\text{by \eqref{second order}}\nonumber
\end{align}

Thus, 
\begin{equation*}
\psi_{FFE} = \nabla Z(\mathbf p) +\epsilon_0 \mathbbm 1 +\mathbf \epsilon
\end{equation*}
where $(\epsilon_0, \lVert \epsilon \rVert) = \begin{cases}(-\nabla Z(\mathbf p)'\mathbbm 1/n,O(ncM_2+nc^2L_2))~~~\text{if \eqref{first order} and \eqref{second order} are true}\\
(\frac{c\mu}{2}\mathbbm 1'\nabla Z(\mathbf p)\mathbbm 1-\nabla Z(\mathbf p)'\mathbbm 1/n,O(nc^2L_2))~~~\text{if \eqref{third order} is also true}\end{cases}$. 

For the variance, we use the decomposition 
\begin{equation}
Var[\hat \psi_{FFE}]=E[Var[\hat \psi_{FFE}|\boldsymbol\delta^1,...,\boldsymbol\delta^{R}]]+Var[E[\hat \psi_{FFE}|\boldsymbol\delta^1,...,\boldsymbol\delta^{R}]]
\label{ffe_var0}
\end{equation}

For the first term,
\begin{align}
tr[E[Var[\hat \psi_{FFE}|\boldsymbol\delta^1,...,\boldsymbol\delta^{R}]]]&=tr[E[\frac{1}{R^2}\sum_{i=1}^{R}Var[\frac{\hat Z^{2i-1}((1-c)\mathbf p +c\boldsymbol\delta^i)-Z^{2i}(\mathbf p)}{c}S(\mathbf p,\boldsymbol\delta^i)|\boldsymbol\delta^i]]]\nonumber \\
&=E[tr[\frac{\sigma^2((1-c)\mathbf p +c\boldsymbol\delta)+\sigma^2(\mathbf p)}{Rc^2}S(\mathbf p,\boldsymbol\delta)S(\mathbf p,\boldsymbol\delta)']]\nonumber \\
&\leq \frac{2\sigma^2}{Rc^2}E[\lVert S(\mathbf p,\boldsymbol\delta) \rVert^2]~~~~~\text{by Assumption \ref{assumption:smooth_noise}}\nonumber\\
&=\frac{2(n-1)\gamma \sigma^2}{Rc^2}~~~~~\text{because $S(\mathbf p,\boldsymbol\delta) = \gamma (\boldsymbol \delta - \mathbf p)$ and by \eqref{second order}}
\label{ffe_var1}
\end{align}

For the second term,
\begin{align}
&tr[Var[E[\hat \psi_{FFE}|\boldsymbol\delta^1,...,\boldsymbol\delta^{R}]]]\nonumber\\
&=\frac{1}{R}tr[Var[\frac{Z((1-c)\mathbf p+c\boldsymbol\delta)-Z(\mathbf p)}{c}S(\mathbf p,\boldsymbol\delta)]]\nonumber \\
&\leq \frac{1}{Rc^2}E[(Z((1-c)\mathbf p+c\boldsymbol\delta)-Z(\mathbf p))^2 \lVert S(\mathbf p,\boldsymbol\delta) \rVert^2]\nonumber \\
& = \frac{1}{Rc^2}E[(c\nabla Z(\mathbf p+c\eta(\boldsymbol\delta - \mathbf p))'(\boldsymbol \delta - \mathbf p))^2\lVert S(\mathbf p,\boldsymbol\delta) \rVert^2]~~~~~\text{for some $\eta \in [0,1]$}\nonumber\\
& = \frac{1}{Rc^2}E[(c\nabla Z(\mathbf p)'(\boldsymbol \delta - \mathbf p)+c(\nabla Z(\mathbf p+c\eta(\boldsymbol\delta - \mathbf p))-\nabla Z(\mathbf p))'(\boldsymbol \delta - \mathbf p))^2\lVert S(\mathbf p,\boldsymbol\delta) \rVert^2]\nonumber\\
&\leq \frac{2}{Rc^2}E[((c\nabla Z(\mathbf p)'(\boldsymbol \delta - \mathbf p))^2+(c(\nabla Z(\mathbf p+c\eta(\boldsymbol\delta - \mathbf p))-\nabla Z(\mathbf p))'(\boldsymbol \delta - \mathbf p))^2)\lVert S(\mathbf p,\boldsymbol\delta) \rVert^2]~~~\text{by\ }(a+b)^2\leq 2(a^2+b^2)\nonumber\\
&\leq \frac{2}{Rc^2}E[c^2\lVert \nabla Z(\mathbf p)\rVert^2\cdot\lVert \boldsymbol \delta - \mathbf p\rVert^2\cdot \lVert S(\mathbf p,\boldsymbol\delta) \rVert^2]\nonumber\\
&~~+\frac{2}{Rc^2}E[c^2\lVert \nabla Z(\mathbf p+c\eta(\boldsymbol\delta - \mathbf p))-\nabla Z(\mathbf p)\rVert^2 \cdot \lVert \boldsymbol \delta - \mathbf p\rVert^2\cdot \lVert S(\mathbf p,\boldsymbol\delta) \rVert^2]\nonumber\\
& \leq \frac{2}{Rc^2}(4c^2M_1^2+16c^4L_1^2)E[\lVert S(\mathbf p,\boldsymbol\delta) \rVert^2]\nonumber\\
& = \frac{ 8\gamma (n-1)(M_1^2+4c^2L_1^2)}{R}~~~~~\text{because $S(\mathbf p,\boldsymbol\delta) = \gamma (\boldsymbol \delta - \mathbf p)$ and by \eqref{second order}}
\label{ffe_var2}
\end{align}

Combining \eqref{ffe_var1} and \eqref{ffe_var2} with \eqref{ffe_var0}, we get
\begin{equation*}
tr[Var[\hat{\psi}_{FFE}]] \leq \frac{2\gamma (n-1) }{Rc^2}[\sigma^2+4c^2M_1^2+16c^4L_1^2]
\end{equation*}\hfill\Halmos
\endproof

\proof{Proof of Theorem \ref{thm:CFE}.}
For the expectation, 
\begin{align}
\psi_{CFE} &= E[\hat \psi_{CFE}] = E[\frac{1}{R}\sum_{i=1}^{R}E[\frac{\hat Z^{2i-1}((1-c)\mathbf p+c\boldsymbol\delta^i)-\hat Z^{2i}((1+c)\mathbf p-c\boldsymbol\delta^i)}{2c}S(\mathbf p,\boldsymbol\delta ^i)|\boldsymbol\delta^i]]\nonumber \\
&=E[\frac{1}{R}\sum_{i=1}^{R}\frac{Z((1-c)\mathbf p+c\boldsymbol\delta^i)-Z((1+c)\mathbf p-c\boldsymbol\delta^i)}{2c}S(\mathbf p,\boldsymbol\delta ^i)]\nonumber \\
&=E[\gamma \frac{Z((1-c)\mathbf p+c\boldsymbol\delta)-Z((1+c)\mathbf p-c\boldsymbol\delta)}{2c}(\boldsymbol\delta -\mathbf p)]~~~\text{noting that $S(\mathbf p,\boldsymbol\delta)=\gamma(\boldsymbol\delta -\mathbf p)$}
\label{eq:cfe_bias1}
\end{align}

By Lemma \ref{lemma:taylor}, we can write $Z(\mathbf p+c (\boldsymbol\delta-\mathbf p))-Z((1+c)\mathbf p-c\boldsymbol\delta)$ as 
\begin{align}
    Z(\mathbf p+c (\boldsymbol\delta-\mathbf p)) -Z((1+c)\mathbf p-c\boldsymbol\delta)&= 2c\nabla Z(\mathbf p)'(\boldsymbol\delta-\mathbf p)+\frac{c^2}{2}(\boldsymbol\delta-\mathbf p)'\nabla^2Z(\mathbf p+\eta_1 c (\boldsymbol\delta-\mathbf p))(\boldsymbol\delta-\mathbf p) \nonumber\\
    &- \frac{c^2}{2}(\boldsymbol\delta-\mathbf p)'\nabla^2Z(\mathbf p-\eta_2 c (\boldsymbol\delta-\mathbf p))(\boldsymbol\delta-\mathbf p)
\label{eq:cfe_taylor}
\end{align}
for some $\eta_1, \eta_2\in[0,1]$ that depends on $\mathbf p, \boldsymbol \delta$. Combining \eqref{eq:cfe_bias1} and \eqref{eq:cfe_taylor}, we get
\begin{align}
\psi_{CFE}&= \gamma E[\nabla Z(\mathbf p)'(\boldsymbol\delta-\mathbf p)(\boldsymbol\delta-\mathbf p)]+D\nonumber \\
&=(\mathbf I-\mathbbm 1\mathbbm 1'/n) \nabla Z(\mathbf p)+D\nonumber ~~~\text{by \eqref{first order} and \eqref{second order}}\\
&=\nabla Z(\mathbf p)-\frac{\mathbbm 1'\nabla Z(\mathbf p)}{n}\mathbbm 1+D
\label{eq:cfe_bias2}
\end{align}
where $D = \frac{\gamma c}{4}((\boldsymbol\delta-\mathbf p)'\nabla^2Z(\mathbf p+\eta_1 c (\boldsymbol\delta-\mathbf p))(\boldsymbol\delta-\mathbf p) - (\boldsymbol\delta-\mathbf p)'\nabla^2Z(\mathbf p-\eta_2 c (\boldsymbol\delta-\mathbf p))(\boldsymbol\delta-\mathbf p))(\boldsymbol\delta-\mathbf p)$.

The first and second terms in \eqref{eq:cfe_bias2} give the correct gradient up to a translation of $\mathbbm 1$. For the last term $D$, we have
\begin{align}
\lVert D \rVert &\leq \frac{\gamma c}{4}E[\lVert (\boldsymbol\delta - \mathbf p)'(\nabla^2Z(\mathbf p+c\eta_1(\boldsymbol\delta - \mathbf p))-\nabla^2Z(\mathbf p-c\eta_2(\boldsymbol\delta - \mathbf p)))(\boldsymbol\delta - \mathbf p)(\boldsymbol\delta - \mathbf p)\rVert]\nonumber \\
&\leq \frac{\gamma c}{4}E[\lVert (\boldsymbol\delta - \mathbf p)\rVert^3 \lVert \nabla^2Z(\mathbf p+c\eta_1(\boldsymbol\delta - \mathbf p))-\nabla^2Z(\mathbf p-c\eta_2(\boldsymbol\delta - \mathbf p))\rVert ]\nonumber \\
& \leq \frac{\gamma c^2L_2}{2}E[\lVert \boldsymbol\delta - \mathbf p\rVert^4]~~~~~\text{because $\nabla^2 Z$ is $L_2$-Lipschitz}\nonumber \\
& \leq 2\gamma c^2L_2E[\lVert \boldsymbol\delta - \mathbf p\rVert^2]~~~~~\text{because $\lVert \boldsymbol\delta - \mathbf p\rVert \leq 2$}\nonumber \\
& = 2(n-1)c^2L_2~~~~~\text{by \eqref{second order}}\nonumber
\end{align}

Thus, 
\begin{equation*}
\psi_{CFE} = \nabla Z(\mathbf p) +\epsilon_0 \mathbbm 1 +\mathbf \epsilon
\end{equation*}
where $\epsilon_0=-\mathbbm 1'\nabla Z(\mathbf p)/n, \lVert \epsilon \rVert = O(nc^2L_2)$. 

For the variance, we use the decomposition 
\begin{equation}
Var[\hat \psi_{CFE}]=E[Var[\hat \psi_{CFE}|\boldsymbol\delta^1,...,\boldsymbol\delta^{R}]]+Var[E[\hat \psi_{CFE}|\boldsymbol\delta^1,...,\boldsymbol\delta^{R}]]
\label{cfe_var0}
\end{equation}

For the first term,
\begin{align}
&tr[E[Var[\hat \psi_{CFE}|\boldsymbol\delta^1,...,\boldsymbol\delta^{R}]]]\nonumber\\
&=tr[E[\frac{1}{R^2}\sum_{i=1}^{R}Var[\frac{\hat Z^{2i-1}((1-c)\mathbf p +c\boldsymbol\delta^i)-Z^{2i}((1+c)\mathbf p -c\boldsymbol\delta^i)}{2c}S(\mathbf p,\boldsymbol\delta^i)|\boldsymbol\delta^i]]]\nonumber \\
&=E[tr[\frac{\sigma^2((1-c)\mathbf p +c\boldsymbol\delta)+\sigma^2((1+c)\mathbf p -c\boldsymbol\delta)}{4Rc^2}S(\mathbf p,\boldsymbol\delta)S(\mathbf p,\boldsymbol\delta)']]\nonumber \\
&\leq \frac{\sigma^2}{2Rc^2}E[\lVert S(\mathbf p,\boldsymbol\delta) \rVert^2]~~~~~\text{by Assumption \ref{assumption:smooth_noise}}\nonumber\\
&=\frac{\gamma \sigma^2(n-1)}{2Rc^2}~~~~~\text{because $S(\mathbf p,\boldsymbol\delta) = \gamma (\boldsymbol \delta - \mathbf p)$ and by \eqref{second order}}
\label{cfe_var1}
\end{align}

For the second term,
\begin{align}
&tr[Var[E[\hat \psi_{CFE}|\boldsymbol\delta^1,...,\boldsymbol\delta^{R}]]]\nonumber\\
&=\frac{1}{R}tr[Var[\frac{Z((1-c)\mathbf p+c\boldsymbol\delta)-Z((1+c)\mathbf p -c\boldsymbol\delta)}{2c}S(\mathbf p,\boldsymbol\delta)]]\nonumber \\
&\leq \frac{1}{4Rc^2}E[(Z((1-c)\mathbf p+c\boldsymbol\delta)-Z((1+c)\mathbf p -c\boldsymbol\delta))^2 \lVert S(\mathbf p,\boldsymbol\delta) \rVert^2]\nonumber \\
& = \frac{1}{4Rc^2}E[(2c\nabla Z(\mathbf p+c\eta(\boldsymbol\delta - \mathbf p))'(\boldsymbol \delta - \mathbf p))^2\lVert S(\mathbf p,\boldsymbol\delta) \rVert^2]~~~~~\text{for some $\eta \in [-1,1]$}\nonumber\\
& = \frac{1}{Rc^2}E[(c\nabla Z(\mathbf p)'(\boldsymbol \delta - \mathbf p)+c(\nabla Z(\mathbf p+c\eta(\boldsymbol\delta - \mathbf p))-\nabla Z(\mathbf p))'(\boldsymbol \delta - \mathbf p))^2\lVert S(\mathbf p,\boldsymbol\delta) \rVert^2]\nonumber\\
&\leq \frac{2}{Rc^2}E[((c\nabla Z(\mathbf p)'(\boldsymbol \delta - \mathbf p))^2+(c(\nabla Z(\mathbf p+c\eta(\boldsymbol\delta - \mathbf p))-\nabla Z(\mathbf p))'(\boldsymbol \delta - \mathbf p))^2)\lVert S(\mathbf p,\boldsymbol\delta) \rVert^2]~~~\text{by\ }(a+b)^2\leq 2(a^2+b^2)\nonumber\\
&\leq \frac{2}{Rc^2}E[c^2\lVert \nabla Z(\mathbf p)\rVert^2\cdot\lVert \boldsymbol \delta - \mathbf p\rVert^2\cdot \lVert S(\mathbf p,\boldsymbol\delta) \rVert^2]\nonumber\\
&~~+\frac{2}{Rc^2}E[c^2\lVert \nabla Z(\mathbf p+c\eta(\boldsymbol\delta - \mathbf p))-\nabla Z(\mathbf p)\rVert^2 \cdot \lVert \boldsymbol \delta - \mathbf p\rVert^2\cdot \lVert S(\mathbf p,\boldsymbol\delta) \rVert^2]\nonumber\\
& \leq \frac{2}{Rc^2}(4c^2M_1^2+16c^4L_1^2)E[\lVert S(\mathbf p,\boldsymbol\delta) \rVert^2]\nonumber\\
& = \frac{ 8\gamma (n-1)(M_1^2+4c^2L_1^2)}{R}~~~~~\text{because $S(\mathbf p,\boldsymbol\delta) = \gamma (\boldsymbol \delta - \mathbf p)$ and by \eqref{second order}}\nonumber\\
&\leq \frac{L_0^2\gamma^2}{R} E[ \lVert \boldsymbol\delta - \mathbf p \rVert^4]~~~~~\text{because $Z$ is $L_0$-Lipschitz}\nonumber \\
\label{cfe_var2}
\end{align}

Combining \eqref{cfe_var1} and \eqref{cfe_var2} with \eqref{cfe_var0}, we get
\begin{equation*}
tr[Var[\hat{\psi}_{CFE}]] \leq \frac{\gamma (n-1) }{Rc^2}[\sigma^2/2+ 8c^2M_1^2 + 32c^4L_1^2]
\end{equation*}\hfill\Halmos
\endproof

\section{Proofs for Section \ref{sec:Dir}}
\proof{Proof of Lemma \ref{lemma_third_order}.}
The proof is by direct calculation using the formula for higher moment \eqref{eq:higher_moments}. We divide into three cases.

Case 1: $i=j=k$.
\begin{align}
&E[(\boldsymbol\delta^l-E[\boldsymbol\delta^l])_i(\boldsymbol\delta^l-E[\boldsymbol\delta^l])_j(\boldsymbol\delta^l-E[\boldsymbol\delta^l])_k]\nonumber \\
&=E[(\boldsymbol\delta^l_i)^3]-3E[\boldsymbol\delta^l_i]E[(\boldsymbol\delta^l_i)^2]+2E[\boldsymbol\delta^l_i]^3\nonumber \\
&=\frac{\Gamma(\alpha^l_0)}{\Gamma(\alpha^l_0+3)}\cdot\frac{\Gamma(\alpha^l_i+3)}{\Gamma(\alpha^l_i)}-3\frac{\alpha^l_i}{\alpha^l_0}\frac{\Gamma(\alpha^l_0)}{\Gamma(\alpha^l_0+2)}\cdot\frac{\Gamma(\alpha^l_i+2)}{\Gamma(\alpha^l_i)}+2(\frac{\alpha^l_i}{\alpha^l_0})^3\nonumber \\
&=\frac{\alpha^l_i(\alpha^l_i+1)(\alpha^l_i+2)}{\alpha^l_0(\alpha^l_0+1)(\alpha^l_0+2)}-\frac{3(\alpha^l_i)^2(\alpha^l_i+1)}{(\alpha^l_0)^2(\alpha^l_0+1)}+2(\frac{\alpha^l_i}{\alpha^l_0})^3\nonumber \\
&=\frac{2(\alpha_0^l)^2\alpha_i^l-6\alpha_0^l(\alpha_i^l)^2+4(\alpha_i^l)^3}{(\alpha^l_0)^3(\alpha^l_0+1)(\alpha^l_0+2)}\nonumber \\
&=\frac{4(\alpha^l_i/\alpha^l_0)^3-6(\alpha^l_i/\alpha^l_0)^2+2(\alpha^l/\alpha^l_0)}{(\alpha^l_0+1)(\alpha^l_0+2)}\nonumber
\end{align}

Case 2: $i=j\neq k$.
\begin{align}
&E[(\boldsymbol\delta^l-E[\boldsymbol\delta^l])_i(\boldsymbol\delta^l-E[\boldsymbol\delta^l])_j(\boldsymbol\delta^l-E[\boldsymbol\delta^l])_k]\nonumber \\
&=E[(\boldsymbol\delta^l_i)^2(\boldsymbol\delta^l_k)]-2E[\boldsymbol\delta^l_i]E[(\boldsymbol\delta^l_i)(\boldsymbol\delta^l_k)]-E[\boldsymbol\delta^l_k]E[(\boldsymbol\delta^l_i)^2]+2E[\boldsymbol\delta^l_i]^2E[\boldsymbol\delta^l_k]\nonumber \\
&=\frac{\Gamma(\alpha^l_0)}{\Gamma(\alpha^l_0+3)}\cdot\frac{\Gamma(\alpha^l_i+2)}{\Gamma(\alpha^l_i)}\cdot\frac{\Gamma(\alpha^l_k+1)}{\Gamma(\alpha^l_k)}-\frac{\Gamma(\alpha^l_0)}{\Gamma(\alpha^l_0+2)}\cdot(2\frac{\alpha^l_i}{\alpha^l_0}\frac{\Gamma(\alpha^l_i+1)}{\Gamma(\alpha^l_i)}\frac{\Gamma(\alpha^l_k+1)}{\Gamma(\alpha^l_k)}+\frac{\alpha^l_k}{\alpha^l_0}\frac{\Gamma(\alpha^l_i+2)}{\Gamma(\alpha^l_i)}))+2(\frac{\alpha^l_k}{\alpha^l_0})(\frac{\alpha^l_i}{\alpha^l_0})^2\nonumber \\
&=\frac{\alpha^l_i(\alpha^l_i+1)\alpha^l_k}{\alpha^l_0(\alpha^l_0+1)(\alpha^l_0+2)}-\frac{2(\alpha^l_i)^2\alpha^l_k+\alpha^l_i(\alpha^l_i+1)\alpha^l_k}{(\alpha^l_0)^2(\alpha^l_0+1)}+2(\frac{\alpha^l_i}{\alpha^l_0})^2\frac{\alpha^l_k}{\alpha^l_0}\nonumber \\
&=\frac{-2\alpha_0^l\alpha_i^l\alpha_k^l+4(\alpha_i^l)^2\alpha_k^l}{(\alpha^l_0)^3(\alpha^l_0+1)(\alpha^l_0+2)}\nonumber \\
&=\frac{4(\alpha^l_i/\alpha^l_0)^2(\alpha^l_k/\alpha^l_0)-2(\alpha^l_i/\alpha^l_0)(\alpha^l_k/\alpha^l_0)}{(\alpha^l_0+1)(\alpha^l_0+2)}\nonumber
\end{align}

Case 3: $i\neq j\neq k$.
\begin{align}
&E[(\boldsymbol\delta^l-E[\boldsymbol\delta^l])_i(\boldsymbol\delta^l-E[\boldsymbol\delta^l])_j(\boldsymbol\delta^l-E[\boldsymbol\delta^l])_k]\nonumber \\
&=\frac{\Gamma(\alpha^l_0)}{\Gamma(\alpha^l_0+3)}(\alpha^l_i\alpha^l_j\alpha^l_k)-\frac{\Gamma(\alpha^l_0)}{\Gamma(\alpha^l_0+2)}\cdot(3\alpha^l_i\alpha^l_j\alpha^l_k/\alpha^l_0)+2\alpha^l_i\alpha^l_j\alpha^l_k/(\alpha^l_0)^3\nonumber \\
&=\frac{4\alpha^l_i\alpha^l_j\alpha^l_k}{(\alpha^l_0)^3(\alpha^l_0+1)(\alpha^l_0+2)}\nonumber \\
&=\frac{4(\alpha^l_i/\alpha^l_0)(\alpha^l_j/\alpha^l_0)(\alpha^l_k/\alpha^l_0)}{(\alpha^l_0+1)(\alpha^l_0+2)}\nonumber
\end{align}\hfill\Halmos
\endproof

\proof{Proof of Lemma \ref{lemma:theta_increasing}.}
For $i=2,\ldots,n-1$,
\begin{equation*}
\theta^i=\frac{2p_i-\sum_{j=1}^{i-1}\theta^j}{n-i}\geq \frac{2p_{i-1}-\sum_{j=1}^{i-1}\theta^j}{n-i} =\frac{(n-i+1)\theta^{i-1}-\theta^{i-1}}{n-i}=\theta^{i-1}
\end{equation*}\hfill\Halmos
\endproof

\proof{Proof of Theorem \ref{thm:delta1}.}
For the first moment condition \eqref{first order}, we notice that $E[\boldsymbol\delta^{l,i}]=(\mathbf e_i+\mathbf e_l)/2$, and $E[\boldsymbol\delta^n] = \mathbf e_n$. Thus
\begin{align}
E[\boldsymbol\delta]&=\sum_{l,i}\theta^{l,i}(\mathbf e_i+\mathbf e_l)/2+\theta^n\mathbf e_n\nonumber \\
&=\frac{1}{2}\sum_{i=1}^{n-1} \mathbf e_i(\sum_{j=1}^{i}\theta^j+(n-i)\theta^i)+(\frac{1}{2}\sum_{i=1}^{n-1}\theta^i+\theta^n)\mathbf e_n\nonumber \\
&=\sum_{i=1}^n p_i\mathbf e_i = \mathbf p\nonumber
\end{align}

For the second moment condition \eqref{second order}, we know that for each $l,i$
\begin{equation*}
E[(\boldsymbol\delta^{l,i}-E[\boldsymbol\delta^{l,i}])(\boldsymbol\delta^{l,i}-E[\boldsymbol\delta^{l,i}])'] = \frac{1}{C(\theta^{l,i})^2}(diag(\mathbf e_i+\mathbf e_l)/2-(\mathbf e_i+\mathbf e_l)(\mathbf e_i+\mathbf e_l)'/4)
\end{equation*}
and for $n$
\begin{equation*}
E[(\boldsymbol\delta^n-E[\boldsymbol\delta^n])(\boldsymbol\delta^n-E[\boldsymbol\delta^n])']=0
\end{equation*}

Thus,
\begin{align}
E[(\boldsymbol\delta - E[\boldsymbol\delta])(\boldsymbol\delta - E[\boldsymbol\delta])'] &= \sum_{l,i}(\theta^{l,i})^2E[(\boldsymbol\delta^{l,i}-E[\boldsymbol\delta^{l,i}])(\boldsymbol\delta^{l,i}-E[\boldsymbol\delta^{l,i}])']~~~\text{by independence of $\boldsymbol\delta^{l,i}$}\nonumber \\
&=\frac{1}{C}\sum_{l,i} (diag(\mathbf e_i+\mathbf e_l)/2-(\mathbf e_i+\mathbf e_l)(\mathbf e_i+\mathbf e_l)'/4)\nonumber \\
&= \frac{1}{4C}\sum_{l,i}(diag(\mathbf e_i+\mathbf e_l)-\mathbf e_i\mathbf e_l'-\mathbf e_l\mathbf e_i')\nonumber \\
&=\frac{1}{4C}(n\mathbf I-\mathbbm 1\mathbbm 1')\nonumber
\end{align}

For the third moment condition \eqref{third order}, we first write
\begin{align}
&E[(\boldsymbol\delta-E[\boldsymbol\delta])_i(\boldsymbol\delta-E[\boldsymbol\delta])_j(\boldsymbol\delta-E[\boldsymbol\delta])_k] \nonumber \\
&= \sum_{(l_1,i_1),(l_2,i_2),(l_3,i_3)}\theta^{l_1,i_1}\theta^{l_2,i_2}\theta^{l_3,i_3}E[(\boldsymbol\delta^{l_1,i_1}-E[\boldsymbol\delta^{l_1,i_1}])_i(\boldsymbol\delta^{l_2,i_2}-E[\boldsymbol\delta^{l_2,i_2}])_j(\boldsymbol\delta^{l_3,i_3}-E[\boldsymbol\delta^{l_3,i_3}])_k]\nonumber \\
&= \sum_{(l_0,i_0)}(\theta^{l_0,i_0})^3E[(\boldsymbol\delta^{l_0,i_0}-E[\boldsymbol\delta^{l_0,i_0}])_i(\boldsymbol\delta^{l_0,i_0}-E[\boldsymbol\delta^{l_0,i_0}])_j(\boldsymbol\delta^{l_0,i_0}-E[\boldsymbol\delta^{l_0,i_0}])_k]\nonumber \\
&=\sum_{(l_0,i_0)}(\theta^{l_0,i_0})^3 0=0\nonumber
\end{align}
Here, Lemma \ref{lemma_third_order} and our design of $\boldsymbol\delta^*$ imply that $E[(\boldsymbol\delta^{l_0,i_0}-E[\boldsymbol\delta^{l_0,i_0}])_i(\boldsymbol\delta^{l_0,i_0}-E[\boldsymbol\delta^{l_0,i_0}])_j(\boldsymbol\delta^{l_0,i_0}-E[\boldsymbol\delta^{l_0,i_0}])_k]=0$ for all $l_0,i_0$ and all indices $i,j,k$.\hfill\Halmos
\endproof

\proof{Proof of Theorem \ref{thm:delta2}.}
For the first moment condition \eqref{first order}, 
\begin{align}
E[\boldsymbol\delta] &= \sum_{l=1}^n \theta^l E[\boldsymbol\delta^l]\nonumber \\
&= np_1 \frac{\mathbbm 1}{n}+\sum_{l=2}^n (p_l-p_1)\mathbf e_l\nonumber \\
&=\sum_{l=1}^n p_l\mathbf e_l = \mathbf p\nonumber
\end{align}

For the second moment condition \eqref{second order},
\begin{align}
\gamma E[(\boldsymbol\delta-\mathbf p)(\boldsymbol\delta-\mathbf p)']&=\gamma \sum_{l=1}^n (\theta^l)^2E[(\boldsymbol\delta^l-E[\boldsymbol\delta^l])(\boldsymbol\delta^l-E[\boldsymbol\delta^l])']~~~\text{by independence of $\boldsymbol\delta^{l}$}\nonumber \\
& = \gamma (\theta^1)^2E[(\boldsymbol\delta^1-E[\boldsymbol\delta^1])(\boldsymbol\delta^1-E[\boldsymbol\delta^1])']\nonumber \\
&=n(n^{\eta+1}+1)\frac{diag(\mathbbm 1/n)-(\mathbbm 1/n)(\mathbbm 1/n)'}{n^{\eta+1}+1}\nonumber \\
&=I - \frac{1}{n}\mathbbm 1\mathbbm 1'\nonumber
\end{align}

For the third moment,
\begin{align}
\gamma E[(\boldsymbol\delta -\mathbf p)_i(\boldsymbol\delta -\mathbf p)_j(\boldsymbol\delta -\mathbf p)_k]&=\sum_{l_1,l_2,l_3}\theta^{l_1}\theta^{l_2}\theta^{l_3}E[(\boldsymbol\delta^{l_1} -E[\boldsymbol\delta^{l_1}])_i(\boldsymbol\delta^{l_2} -E[\boldsymbol\delta^{l_2}])_j(\boldsymbol\delta^{l_3} -E[\boldsymbol\delta^{l_3}])_k]\nonumber \\
&= \gamma \sum_{l}(\theta^{l})^3E[(\boldsymbol\delta^{l} -E[\boldsymbol\delta^{l}])_i(\boldsymbol\delta^{l} -E[\boldsymbol\delta^{l}])_j(\boldsymbol\delta^{l} -E[\boldsymbol\delta^{l}])_k]\nonumber \\
&=\gamma (\theta^1)^3E[(\boldsymbol\delta^{1} -E[\boldsymbol\delta^{1}])_i(\boldsymbol\delta^{1} -E[\boldsymbol\delta^{1}])_j(\boldsymbol\delta^{1} -E[\boldsymbol\delta^{1}])_k]\nonumber \\
&=\frac{n^2p_1}{n^{\eta+1}+2}\begin{cases}
4/n^3-6/n^2+2/n~~~i=j=k\\
4/n^3-2/n^2~~~i=j\neq k\\
4/n^3~~~i\neq j\neq k
\end{cases}\nonumber
\end{align}

Thus, the $k$-th component of the Hessian term in the bias becomes
\begin{align}
&\frac{\gamma c}{2} E[(\sum_{i,j}(\boldsymbol\delta -\mathbf p)_i(\boldsymbol\delta -\mathbf p)_j\nabla^2 Z)(\boldsymbol\delta - \mathbf p)_k]\nonumber\\
&=\frac{n^2p_1c}{2(n^{\eta+1}+2)}(\mathbbm 1 \nabla ^2 Z\mathbbm 1\cdot \frac{4}{n^3}-\mathbbm 1'diag(\nabla^2 Z)\cdot \frac{2}{n^2}+\nabla^2Z_{kk}\cdot (\frac{2}{n}-\frac{4}{n^2}))\nonumber \\
&=\frac{n^2p_1c}{2(n^{\eta+1}+2)}(\mathbbm 1 \nabla ^2 Z\mathbbm 1\cdot \frac{4}{n^3}-\mathbbm 1'diag(\nabla^2 Z)\cdot \frac{2}{n^2})+\frac{np_1c(1-2/n)}{n^{\eta+1}+2}\nabla^2Z_{kk}\nonumber
\end{align}
Notice that the first term of above appears in $\frac{\gamma c}{2} E[(\sum_{i,j}(\boldsymbol\delta -\mathbf p)_i(\boldsymbol\delta -\mathbf p)_j\nabla^2 Z)(\boldsymbol\delta - \mathbf p)_k]$ for all index $k$, and thus  we can see it as part of a translation of $\mathbbm 1$. The only difference between different $k$ is the second term $\lvert \frac{np_1c(1-2/n)}{n^{\eta+1}+2}\nabla^2Z_{kk} \rvert \leq c\lVert \nabla ^2 Z \rVert_{\infty}n^{-\eta -1}$. \hfill\Halmos
\endproof

\proof{Proof of Corollary \ref{cor:delta2_to_zero}.}
This is a direct consequence of Theorem \ref{thm:delta2}.\hfill\Halmos
\endproof

\proof{Proof of Corollary \ref{cor:delta2_comparable}.}
Theorem \ref{thm:delta2} says $ \lVert \frac{c}{2}E[(\boldsymbol\delta - \mathbf p)'\nabla ^2 Z(\mathbf p)  (\boldsymbol\delta - \mathbf p)S(\mathbf p,\boldsymbol\delta)] -\epsilon_{0} \mathbbm 1\rVert\leq \sqrt{n}\cdot O(c\lVert \nabla ^2 Z \rVert_{\infty} n^{-\eta-1}) = O(cM_2n^{-\eta-\frac{1}{2}})$ where the first inequality is because $\lVert \mathbf x \rVert \leq \sqrt{n}\lVert \mathbf x \rVert_{\infty}\forall \mathbf x\in \mathbbm R ^n$, and the second inequality is because $\|\nabla^2 Z\|_{\infty} = \max_{i,j}\lvert \mathbf e_i'\nabla^2Z\mathbf e_j\rvert \leq \max_{i,j} \lVert \mathbf e_i\rVert \cdot \lVert \nabla^2 Z\rVert \cdot \lVert \mathbf e_j\rVert=\lVert \nabla^2 Z\rVert\leq M_2$. And the corollary follows by taking $\eta = -1$.\hfill\Halmos
\endproof

\section{Proofs for Section \ref{sec:comparisons}}

\proof{Proof of Theorem \ref{thm:FD}.}
For the expectation, we have 
\begin{align}
(\psi_{FD,standard})_i &= E[(\hat \psi_{FD,standard})_i] = \frac{Z(\mathbf p +c(\mathbf e_i-\mathbf p))-Z(\mathbf p)}{c}\nonumber \\
& =\frac{c\nabla Z(\mathbf p+c\eta(\mathbf e_i-\mathbf p))'(\mathbf e_i-\mathbf p)}{c}~~~\text{for some $\eta \in [0,1]$ by Lemma \ref{lemma:taylor}} \nonumber \\
& =\nabla Z(\mathbf p) '(\mathbf e_i-\mathbf p) + c\eta (\mathbf e_i-\mathbf p)'\nabla^2 Z(\mathbf p+c\eta\xi(\mathbf e_i-\mathbf p))(\mathbf e_i-\mathbf p)~~~\text{for some $\xi \in [0,1]$} \nonumber \\
& = \nabla Z(\mathbf p)_i-\nabla Z(\mathbf p) '\mathbf p +\epsilon_i\nonumber
\end{align}
where $\epsilon_i = c\eta (\mathbf e_i-\mathbf p)'\nabla^2 Z(\mathbf p+c\eta\xi(\mathbf e_i-\mathbf p))(\mathbf e_i-\mathbf p)$, and $\lvert \epsilon_i \rvert \leq cM_2\lVert \mathbf e_i - \mathbf p\rVert^2\leq 4cM_2$ where $M_2 = \sup_{\mathbf p\in \mathcal P} \lVert \nabla^2 Z(\mathbf p) \rVert$. 

Thus,
\begin{equation*}
\psi_{FD,standard} = \nabla Z(\mathbf p) - \nabla Z(\mathbf p)'\mathbf p \mathbbm 1 +\epsilon
\end{equation*}
and we have $\lVert \epsilon \rVert \leq 4cM_2\sqrt{n}$.

For the variance,
\begin{equation*}
Var[(\hat \psi_{FD,standard})_i] \leq \frac{2\sigma^2}{R_pc^2}
\end{equation*}
so 
\begin{equation*}
tr[Var[\hat \psi_{FD,standard}]] \leq \frac{2n\sigma^2}{R_pc^2}= \frac{2n^2\sigma^2}{Rc^2}
\end{equation*}\hfill\Halmos
\endproof

\proof{Proof of Theorem \ref{thm:FDrd}.}
For the expectation, 
\begin{align}
\psi_{FD,random} &= E[\hat \psi_{FD,random}] =E[\hat \psi_{FD,random}^1] \nonumber\\
&= \frac{n}{c}E[(Z(\mathbf p +c(\mathbf e_{l_1} - \mathbf p))-Z(\mathbf p))\mathbf e_{l_1}] \nonumber\\
&=\psi_{FD}= \nabla Z(\mathbf p) - \nabla Z(\mathbf p)'\mathbf p \mathbbm 1+\epsilon\nonumber
\end{align}
so the bias of $\hat \psi_{FD,random}$ is the same as the bias of $\hat \psi_{FD}$.

For the variance, we note that
\begin{equation*}
Var[ \hat \psi_{FD,random}] = \frac{1}{R}Var[\hat \psi_{FD,random}^1]
\end{equation*}

So we use the decomposition again 
\begin{equation}
Var[ \hat \psi_{FD,random}^1]=E[Var[\hat \psi_{FD,random}^1|l_1]]+Var[E[\hat \psi_{FD,random}^1|l_1]]
\label{eq:fdrd_var0}
\end{equation}

For the first term,
\begin{align}
    tr[E[Var[\hat \psi_{FD,random}^1|l_1]]] &= E[tr[Var[\hat \psi_{FD,random}^1|l_1]]]\nonumber\\
    &= E[tr[\frac{n^2(\sigma^2(\mathbf p+c(\mathbf e_{l_1}-\mathbf p))+\sigma^2(\mathbf p))}{c^2}\mathbf e_{l_1}\mathbf e_{l_1}']]\nonumber\\
    &\leq  \frac{2n^2\sigma^2}{c^2}\label{eq:fdrd_var1}
\end{align}

For the second term,
\begin{align}
&tr[Var[E[\hat \psi_{FD,random}^1|l_1]]]\nonumber\\
&= tr[Var[n \frac{Z(\mathbf p + c(\mathbf e_{l_1}-\mathbf p))-Z(\mathbf p)}{c}\mathbf e_{l_1}]]\nonumber\\
&\leq n^2E[\lVert \frac{Z(\mathbf p + c(\mathbf e_{l_1}-\mathbf p))-Z(\mathbf p)}{c}\mathbf e_{l_1} \rVert^2]\nonumber \\
& = n^2E[(\nabla Z(\mathbf p)'(\mathbf e_{l_1}-\mathbf p)+ (\nabla Z(\mathbf p + c(\mathbf e_{l_1}-\mathbf p))-\nabla Z(\mathbf p))'(\mathbf e_{l_1}-\mathbf p))^2]\nonumber\\
&\leq 2n^2E[(\nabla Z(\mathbf p)'(\mathbf e_{l_1}-\mathbf p))^2+ ((\nabla Z(\mathbf p + c(\mathbf e_{l_1}-\mathbf p))-\nabla Z(\mathbf p))'(\mathbf e_{l_1}-\mathbf p))^2]~~~~~(a+b)^2\leq 2a^2+2b^2\nonumber\\
&\leq 2n^2(4M_1^2+16c^2L_1^2)
\label{eq:fdrd_var2}
\end{align}

Combining \eqref{eq:fdrd_var1} and \eqref{eq:fdrd_var2} with \eqref{eq:fdrd_var0}, we get 
\begin{equation*}
tr[Var[\hat \psi_{FD,random}]] \leq \frac{2n^2}{Rc^2}[\sigma^2+2c^2L_0^2]
\end{equation*}\hfill\Halmos
\endproof

\proof{Proof of Lemma \ref{lm:gamma}.}
Condition \eqref{second order} tells us $E[\gamma \lVert \boldsymbol \delta -\mathbf p\rVert^2] = n-1$, and with the simplex constraint $\lVert \boldsymbol \delta -\mathbf p \rVert \leq 2$, we get $\gamma\geq \frac{n-1}{4}$.\hfill\Halmos
\endproof

\section{Proofs for Section \ref{sec:SA}}
\proof{Proof of Lemma \ref{lm:bdd_FWSA}.}
First, by the update rule, we must have $\min_{i\in[m],j\in n^i} p^i_{k,j}\geq \prod_{l=1}^{k-1}(1-\epsilon_l)p_0$. Moreover, we choose $\epsilon_l = \frac{a}{l}$. So 
\begin{align}
\prod_{l=1}^{k-1}(1-\epsilon_l)&=\prod_{l=1}^{k-1}\frac{l-a}{l}\nonumber\\
&=\frac{\Gamma(k-a)}{\Gamma(1-a)(k-1)!}\nonumber\\
&\geq \frac{\Gamma(k-a)}{\Gamma(1/2)(k-1)!}=\frac{\Gamma(k-a)}{\sqrt{\pi}(k-1)!}
\label{eq:eps_bound}
\end{align}
where the inequality holds because the Gamma function is decreasing in $(0,1)$, and $a< \frac{1}{2}$.

From Lemma \ref{lm:gamma_bound} below, if we take $N=1$, then $\lvert R_N(z)\rvert \leq 1/(\pi ^2 z)$, and we get
\begin{align}
&\Gamma(z)\geq \sqrt{\frac{2\pi}{z}}(\frac{z}{e})^z(1-\frac{1}{\pi^2z})\nonumber\\
&\Rightarrow \Gamma(k-a)\geq \sqrt{\frac{2\pi}{e}}(\frac{k-a}{e})^{k-a-1/2}(1-\frac{1}{\pi^2(k-a)})
\label{eq:gamma_lower_bd}
\end{align}

Using Stirling's approximation, we get,
\begin{align}
(k-1)!\leq e(\frac{k-1}{e})^{k-1}(k-1)^{1/2}=e^{3/2}(\frac{k-1}{e})^{k-1/2}
\label{eq:stirling_upper_bd}
\end{align}

Combining \eqref{eq:gamma_lower_bd} and \eqref{eq:stirling_upper_bd} with \eqref{eq:eps_bound}, we get 
\begin{align}
&\frac{\Gamma(k-a)}{\sqrt{\pi}(k-1)!}\nonumber\\
&\geq (\frac{k-a}{k-1})^{k-1/2}\frac{\sqrt{2}e^{-(2-a)}(1-1/(\pi^2(k-a)))}{(k-a)^a}~~~\text{by the lower and upper bound in \eqref{eq:gamma_lower_bd} \eqref{eq:stirling_upper_bd}}\nonumber\\
&\geq \frac{\sqrt{2}e^{-(2-a)}(1-1/(\pi^2(k-a)))}{(k-a)^a}~~~\text{because $a<\frac{1}{2}<1$}\nonumber\\
&\geq \frac{\sqrt{2}e^{-2}(1-1/(\pi^2(k-a)))}{k^a}\nonumber\\
&= \frac{A}{k^a}~~~\text{where }A=\sqrt{2}e^{-2}(1-1/(\pi^2(1-a))>0\nonumber
\end{align}\hfill\Halmos
\endproof

\begin{lemma}[Adapted from \citealt{NIST:DLMF} Section 5.11(ii)]
Let $z\in \mathbbm C$ such that $\lvert \arg(z)\rvert \leq \frac{\pi}{2}$ and define for any $N\geq 1$ the remainder $R_N(z)$ by
\begin{equation*}
    \Gamma(z) = \sqrt{2\pi}z^{z-\frac{1}{2}}e^{-z}(\sum_{n=0}^{N-1}\frac{g_n}{z^n}+R_N(z))
\end{equation*}
where $g_0=1$, $g_1 = 1/12$, $g_2 = 1/288$. (The explicit formulas for $g_k$ is stated in \cite{Nemes:2013:GFSN}).

Then we have 
\begin{equation*}
    \lvert R_N(z)\rvert \leq \frac{(1+\zeta(N))\Gamma(N)}{2(2\pi)^{N+1}\lvert z\rvert^N}(1+\min(\sec(\arg(z)),2\sqrt{N}))
\end{equation*}
where $\zeta(\cdot)$ is the Riemann zeta function, and in the case $N=1$, $1+\zeta(N)$ is replaced with 4. 
\label{lm:gamma_bound}
\end{lemma}

\proof{Proof of Theorem \ref{main_FWSA}.}
The proof is generalized from \cite{ghosh2019robust} that uses only an unbiased gradient estimator to deal with the bias in our zeroth-order estimator. To begin with, by Assumption \ref{assumption:smooth} $Z(\cdot)$ is $L_0$-Lipschitz, and since $\prod_{i=1}^m \mathcal P^{n^i}$ is bounded, we have $Z(\cdot)$ bounded. For notational convenience, we write $\mathbf d_k = \mathbf q_k -\mathbf p_k$ and $\hat {\mathbf d}_k = \hat {\mathbf q}_k -\mathbf p_k$ where $\mathbf q_k$ is the solution to the subproblem with the exact gradient $\nabla Z(\mathbf p_k)$, and $\hat {\mathbf q}_k$ is the solution with the estimated gradient $\hat \psi(\mathbf p_k)$, and $\psi(\mathbf p_k):=E[\hat \psi(\mathbf p_k)] $. 

Thus, the update at each iteration is $\mathbf p_{k+1} = \mathbf p_k+\epsilon_k\hat{\mathbf d}_k$. By Lemma \ref{lemma:taylor}, we have
\begin{equation}
Z(\mathbf p_{k+1}) = Z(\mathbf p_k)+\epsilon_k \nabla Z(\mathbf p_k)'\hat{\mathbf d}_k+\frac{\epsilon_k^2}{2}\hat{\mathbf d}_k'\nabla ^2Z(\mathbf p_k+\xi_kc_k\hat{\mathbf d}_k)\hat{\mathbf d}_k
\label{ec25}
\end{equation}
for some $\xi_k\in[0,1]$. We can decompose the second term of the RHS of \eqref{ec25} as
\begin{align}
\nabla Z(\mathbf p_k)'\hat{\mathbf d}_k &= \hat \psi(\mathbf p_k)'\hat{\mathbf d}_k + (\nabla Z(\mathbf p_k) - \hat \psi(\mathbf p_k))'\hat{\mathbf d}_k \nonumber\\
&\leq \hat \psi(\mathbf p_k)'\mathbf d_k + (\nabla Z(\mathbf p_k) - \hat \psi(\mathbf p_k))'\hat{\mathbf d}_k ~~~~~\text{by optimality of $\hat{\mathbf d}_k$}\nonumber\\
&=\nabla Z(\mathbf p_k)'\mathbf d_k + (\nabla Z(\mathbf p_k) - \hat \psi(\mathbf p_k))'(\hat{\mathbf d}_k - \mathbf d_k)\nonumber\\
&=\nabla Z(\mathbf p_k)'\mathbf d_k+(\psi(\mathbf p_k)-\hat \psi(\mathbf p_k))'(\hat{\mathbf d}_k-\mathbf d_k)+(\nabla Z(\mathbf p_k)+\epsilon_0(\mathbf p_k)\mathbbm 1-\psi(\mathbf p_k))'(\hat{\mathbf d}_k-\mathbf d_k)
\label{ec26}
\end{align}

Combining \eqref{ec26} with \eqref{ec25}, we get
\begin{align}
Z(\mathbf p_{k+1}) &\leq  Z(\mathbf p_k)+\underbrace{\epsilon_k\nabla Z(\mathbf p_k)'\mathbf d_k}_{\Theta_{1,k}}+\underbrace{\epsilon_k(\psi(\mathbf p_k)-\hat \psi(\mathbf p_k))'(\hat{\mathbf d}_k-\mathbf d_k)}_{\Theta_{2,k}}\nonumber\\
&+\underbrace{\epsilon_k(\nabla Z(\mathbf p_k)+\epsilon_0(\mathbf p_k)\mathbbm 1-\psi(\mathbf p_k))'(\hat{\mathbf d}_k-\mathbf d_k)}_{\Theta_{3,k}}+\underbrace{\frac{\epsilon_k^2}{2}\hat{\mathbf d}_k'\nabla ^2Z(\mathbf p_k+\xi_kc_k\hat{\mathbf d}_k)\hat{\mathbf d}_k}_{\Theta_{4,k}}\nonumber
\end{align}

Let $\mathcal F_k$ be the filtration generated by $\mathbf p_1,\ldots,\mathbf p_k$, we then have
\begin{equation}
 E[Z(\mathbf p_{k+1})|\mathcal F_k] \leq  Z(\mathbf p_k)+E[\Theta_{1,k}|\mathcal F_k]+E[\Theta_{2,k}|\mathcal F_k]+E[\Theta_{3,k}|\mathcal F_k]+E[\Theta_{4,k}|\mathcal F_k]
\label{ec27}   
\end{equation}

We analyze the RHS of \eqref{ec27} term by term.

For $E[\Theta_{4,k}|\mathcal F_k]$, we get 
\begin{align}
\lvert E[\Theta_{4,k}|\mathcal F_k] \rvert&\leq \frac{\epsilon_k^2}{2}E[\lvert \hat{\mathbf d}_k'\nabla ^2Z(\mathbf p_k+\xi_kc_k\hat{\mathbf d}_k)\hat{\mathbf d}_k\rvert|\mathcal F_k]\nonumber\\
&\leq \frac{\epsilon_k^2}{2}E[\lVert \hat{\mathbf d}_k \rVert^2 \cdot \lVert \nabla ^2Z(\mathbf p_k+\xi_kc_k\hat{\mathbf d}_k) \rVert|\mathcal F_k]\nonumber\\
&\leq \frac{\epsilon_k^2}{2} E[4mM_2|\mathcal F_k]\nonumber\\
&= \frac{\epsilon_k^2}{2}\cdot 4mM_2=2mM_2\epsilon_k^2
\label{eq:bd_theta4}
\end{align}
where $M_2 =\sup_{\mathbf p \in \prod_{i=1}^m \mathcal P^{n^i}} \lVert \nabla^2 Z(\mathbf p ) \rVert$ and the sup is finite because $\nabla^2 Z$ is $L_2$-Lipschitz and $\prod_{i=1}^m \mathcal P^{n^i}$ is bounded. 

For $E[\Theta_{3,k}|\mathcal F_k]$, we have
\begin{align}
    \lvert E[\Theta_{3,k}|\mathcal F_k] \rvert&\leq \epsilon_kE[\lVert \nabla Z(\mathbf p_k)+\epsilon_0(\mathbf p_k)\mathbbm 1-\psi(\mathbf p_k)\rVert \cdot\lVert \hat{\mathbf d}_k-\mathbf d_k\rVert |\mathcal F_k]\nonumber\\
    &\leq B_1\epsilon_kc_k E[\lVert \hat{\mathbf d}_k-\mathbf d_k\rVert |\mathcal F_k]~~~~~\text{by assumption on the bias of $\hat{\psi}$}\nonumber\\
    &\leq 4\sqrt{m}B_1\epsilon_kc_k
\label{eq:bd_theta3}
\end{align}

For $E[\Theta_{2,k}|\mathcal F_k]$, we first get
\begin{align}
&E[(\psi(\mathbf p_k)-\hat \psi(\mathbf p_k))'(\hat{\mathbf d}_k-\mathbf d_k)|\mathcal F_k]\nonumber \\
&\leq \sqrt{E[\lVert \psi(\mathbf p_k)-\hat \psi(\mathbf p_k) \rVert^2|\mathcal F_k]E[\lVert \hat{\mathbf d}_k-\mathbf d_k \rVert^2|\mathcal F_k]}~~~\text{by the Cauchy-Schwarz inequality}\nonumber \\
&\leq 2\sqrt{m}\cdot \sqrt{E[\lVert \psi(\mathbf p_k)-\hat \psi(\mathbf p_k) \rVert^2|\mathcal F_k]}~~~\text{since $\lVert \hat{\mathbf d}_k-\mathbf d_k \rVert \leq \lVert \hat{\mathbf d}_k \rVert+\lVert \mathbf d_k \rVert\leq 2\sqrt{m} $}\nonumber\\
&\leq 2\sqrt{mB_2}\sqrt{\frac{\gamma_k}{R_kc_k^2}}~~~\text{by assumption on the variance of $\hat{\psi}$}\nonumber
\end{align}
Thus, we get 
\begin{equation}
 \lvert E[\Theta_{2,k}|\mathcal F_k] \rvert \leq 2\sqrt{mB_2}\epsilon_k\sqrt{\frac{\gamma_k}{R_kc_k^2}}
 \label{eq:bd_theta2}
\end{equation}

From Lemma \ref{lm:bdd_FWSA}, if we let $p_0 = \min_{i\in[m], j\in[n^i]} p^i_{1,j}$, we know that $\min_{i\in[m], j\in[n^i]} p^i_{k,j}\geq \frac{Ap_0}{k^a}$ for some $A>0$. Thus, $\gamma_k = \gamma_0 (\min_{i\in[m], j\in[n^i]} p^i_{1,j})^{-2}\leq \tilde A k^{2a}$ for some $\tilde A>0$.

Since $\nabla Z(\mathbf p_k)'\mathbf d_k \leq 0$ by the optimality of $\mathbf d_k$, we have $E[\Theta_{1,k}|\mathcal F_k]\leq 0$. Combining the above bounds \eqref{eq:bd_theta2}, \eqref{eq:bd_theta3}, and \eqref{eq:bd_theta4} with \eqref{ec27} we have 
\begin{align}
E[Z(\mathbf p_{k+1})-Z(\mathbf p_{k})|\mathcal F_k] &\leq C_1\epsilon_k\sqrt{\frac{\gamma_k}{R_kc_k^2}} + C_2\epsilon_kc_k + C_3\epsilon_k^2\nonumber\\
& =\frac{C_1'}{k^{1+\beta/2-a-\theta}}+ \frac{C_2'}{k^{\theta+1}}+\frac{C_3'}{k^{2}}\nonumber
\end{align}
for some constants $C_1,C_2,C_3,C_1',C_2',C_3'>0$.

So we get 
\begin{equation}
\sum_{k=1}^{\infty}E[E[Z(\mathbf p_{k+1})-Z(\mathbf p_k)|\mathcal F_k]^{+}]\leq \sum_{k=1}^{\infty} (\frac{C_1'}{k^{1+\beta/2-a-\theta}}+ \frac{C_2'}{k^{\theta+1}}+\frac{C_3'}{k^{2}})\label{interim appendix}
\end{equation}

Thus, for $\beta/2>a+\theta$, the RHS of \eqref{interim appendix} converges. By Lemma \ref{lm:ec21-cor}, we have $Z(\mathbf p_k)$ converge to an integrable random variable $Z_\infty$.

Now we can take expectation on both sides of \eqref{ec27},
\begin{align}
E[Z(\mathbf p_{k+1})] &\leq  E[Z(\mathbf p_k)]+\epsilon_kE[\nabla Z(\mathbf p_k)'\mathbf d_k]+\epsilon_kE[(\psi(\mathbf p_k)-\hat \psi(\mathbf p_k))'(\hat{\mathbf d}_k-\mathbf d_k)]\nonumber\\
&+\epsilon_k(\nabla Z(\mathbf p_k)+\epsilon_k(\mathbf p_k) \mathbbm 1-\psi(\mathbf p_k))'E[\hat{\mathbf d}_k-\mathbf d_k]+\frac{\epsilon_k^2}{2}E[\hat{\mathbf d}_k'\nabla ^2Z(\mathbf p_k+\xi_kc_k\hat{\mathbf d}_k)\hat{\mathbf d}_k]\nonumber\\
&=E[Z(\mathbf p_k)] + \epsilon_kE[\nabla Z(\mathbf p_k)'\mathbf d_k] + E[\Theta_{2,k}+\Theta_{3,k}+\Theta_{4,k}]
\label{eq:full_exp}
\end{align}
Applying the bounds in \eqref{eq:bd_theta2}, \eqref{eq:bd_theta3}, and \eqref{eq:bd_theta4} and telescoping \eqref{eq:full_exp} we get
\begin{align}
E[Z(\mathbf p_{k+1})]\leq E[Z(\mathbf p_1)]+\sum_{j=1}^{k} \epsilon_kE[\nabla Z(\mathbf p_k)'\mathbf d_k]+ \sum_{j=1}^{k}  (\frac{C_1'}{j^{1+\beta/2-a-\theta}}+ \frac{C_2'}{j^{\theta+1}}+\frac{C_3'}{j^{2}})
\label{ec33}
\end{align}

Take limit on both sides of \eqref{ec33}. First, we have $E[Z(\mathbf p_{k+1})]\to E[Z_{\infty}]$ by the dominated convergence theorem. Also, $Z(\mathbf p_1)<\infty$. Therefore, since $E[\nabla Z(\mathbf p_j)'\mathbf d_j]\leq 0$, we must have $\sum_{j=1}^k\epsilon_jE[\nabla Z(\mathbf p_j)'\mathbf d_j]$ converge a.s., which implies that $\lim\sup_{k\to \infty}E[\nabla Z(\mathbf p_j)'\mathbf d_j]=0$. So there exists a subsequence $k_i$ such that $\lim_{i\to \infty}E[\nabla Z(\mathbf p_{k_i})'\mathbf d_{k_i}]=0$. Thus, $\nabla Z(\mathbf p_{k_i})'\mathbf d_{k_i}\to_{p}0$. Then there exists a further subsequence $l_i$ such that $\nabla Z(\mathbf p_{l_i})'\mathbf d_{l_i}\to0$ a.s..

Now let $S^* = \{\mathbf p\in\prod_{i=1}^m\mathcal P^{n^i}:g(\mathbf p) =0\}$. Since $g(\cdot)$ is continuous, we have $D(\mathbf p_{l_i},S^*)\to 0$ a.s.. Since $Z(\cdot)$ is continuous, we have $D(Z(\mathbf p_{l_i}),\mathcal Z^*)\to 0$ a.s.. But since we have shown that $Z(\mathbf p_k)$ converges a.s., we have $D(Z(\mathbf p_k),\mathcal Z^*)\to 0$ a.s.. This gives part 1 of the theorem.

Now, under Assumptions \ref{assumption:FW gap} and \ref{assumption:uniqueness}, since $\mathbf p^*$ is the only $\mathbf p$ such that $g(\mathbf p)=0$, we must have $\mathbf p_{l_i}\to \mathbf p^*$ a.s.. Since $Z(\cdot)$ is continuous, we have $Z(\mathbf p_{l_i})\to Z(\mathbf p^*)$. But since $Z(\mathbf p_k)$ converges a.s. as shown above, we must have $Z(\mathbf p_k)\to Z(\mathbf p^*)$. Then by Assumption 3, since $\mathbf p^*$ is the unique optimizer, we have $\mathbf p_k\to \mathbf p^*$ a.s.. This gives part 2 of the theorem. \hfill\Halmos
\endproof
\begin{lemma}[Adapted from the corollary in Section 3 in \citealt{blum1954multidimensional}] Consider a sequence of integrable random variable $Y_k,k=1,2,\ldots$. Let $\mathcal F_k$ be the filtration generated by $Y_1,\ldots,Y_k$. Assume
\begin{equation*}
\sum_{k=1}^{\infty} E[E[Y_{k+1}-Y_{k}|\mathcal F_k]^{+}]<0
\end{equation*}
where $x^{+}$ denotes the positive part of $x$, i.e. $x^{+} = x $ if $x\geq 0$ and $0$ if $x<0$. Moreover, assume that $Y_k$ is bounded uniformly from below. Then $Y_k\to Y_{\infty}$ a.s., where $Y_{\infty}$ is a random variable.
\label{lm:ec21-cor}
\end{lemma}


To show a.s. convergence for MDSA, we first recall the following definition and  result:
\begin{definition}[\citealt{doi:10.1137/070704277}]
Let $X\subset \mathbbm R^n$ be a non-empty bounded closed convex set, $w(\cdot):X\to \mathbbm R$ is a \textit{distance generating function} modulus $\alpha>0$ if $w$ is convex and continuous on $X$, the set $X^o:= \{\mathbf x\in X~|~\exists \mathbf p \in \mathbbm R^n \text{ such that } \mathbf x \in argmin_{\mathbf u\in X} [\mathbf p'\mathbf u+w(\mathbf u)]\}$ is convex, and restricted to $X^o$, $w(\cdot)$ is continuously differentiable and strongly convex with parameter $\alpha$, i.e.
\begin{equation*}
    (\mathbf x'-\mathbf x)'(\nabla w(\mathbf x') - \nabla w(\mathbf x))\geq \alpha \lVert \mathbf x'-\mathbf x\rVert^2
\end{equation*}
\end{definition}

\begin{lemma}[Lemma 2.1 in \citealt{doi:10.1137/070704277}]
Let $w(\cdot):X\to \mathbbm R$ be a distance generating function  and $V:X\times X^o\to \mathbbm R$ be defined as 
\begin{equation*}
    V(\mathbf x,\mathbf z):= w(\mathbf z) -(w(\mathbf x)+\nabla w(\mathbf x)'(\mathbf z - \mathbf x))
\end{equation*}

Then for any $\mathbf u\in X$ and $\mathbf x \in X^o$ and $\mathbf y$ the following inequality holds
\begin{equation*}
V(prox_{\mathbf x}(\mathbf y),\mathbf u)\leq V(\mathbf x,\mathbf u)+\mathbf y'(\mathbf u - \mathbf x)+\frac{\lVert \mathbf y \rVert^2}{2\alpha}
\end{equation*}
\label{lemma_21}
\end{lemma}
We have the following corollary:
\begin{corollary}
\begin{equation*}
\lVert prox_{\mathbf x}(\mathbf y)-\mathbf x\rVert\leq \frac{\lVert \mathbf y\rVert}{\alpha}
\end{equation*}
\label{cor_MDSA}
\end{corollary}
\proof{Proof of Corollary \ref{cor_MDSA}.}
Taking $\mathbf u = \mathbf x$ in Lemma \ref{lemma_21}, we get
\begin{equation*}
V(prox_{\mathbf x}(\mathbf y),\mathbf x)\leq \frac{\lVert \mathbf y \rVert^2}{2\alpha}
\end{equation*}
Since $w(\cdot)$ is strongly convex	with parameter $\alpha$, we have
\begin{equation*}
V(prox_{\mathbf x}(\mathbf y),\mathbf x) \geq \frac{\alpha}{2}\lVert prox_{\mathbf x}(\mathbf y)-\mathbf x\rVert^2
\end{equation*}
Together, we get the corollary.\hfill\Halmos
\endproof

Notice that for our choice of the entropic distance-generating function, Lemma \ref{lemma_21} and Corollary \ref{cor_MDSA} hold, so we get the following: 
\proof{Proof of Theorem \ref{main_MDSA}.}
The proof is generalized from \cite{goeva2019optimization} that uses only an unbiased gradient estimator to deal with the bias in our zeroth-order estimator. We analyze the evolution of $V(\mathbf p_{k},\mathbf p^*)$. First, the update rule is $\mathbf p_{k+1} = prox_{\mathbf p_k}(\rho_k \hat{\psi}(\mathbf p_k))= prox_{\mathbf p_k}(\rho_k (\hat{\psi}(\mathbf p_k)-\epsilon_0(\mathbf p_k)\mathbbm 1))$. Thus Lemma \ref{lemma_21} implies that
\begin{align}
V(\mathbf p_{k+1},\mathbf p^*) &\leq V(\mathbf p_k,\mathbf p^*)+\rho_k (\hat{\psi}(\mathbf p_k)-\epsilon_0(\mathbf p_k)\mathbbm 1)'(\mathbf p^* - \mathbf p_k)+\frac{\rho_k^2 \lVert \hat{\psi}(\mathbf p_k)-\epsilon_0(\mathbf p_k)\mathbbm 1\rVert ^2}{2}\nonumber\\
&\leq V(\mathbf p_k,\mathbf p^*)+\underbrace{\rho_k\nabla Z(\mathbf p_k)'(\mathbf p^* - \mathbf p_k)}_{\Theta_{0,k}}+\underbrace{\rho_k (\hat{\psi}(\mathbf p_k)-\epsilon_0(\mathbf p_k)\mathbbm 1-\nabla Z(\mathbf p_k))'(\mathbf p^* - \mathbf p_k)}_{\Theta_{1,k}}\nonumber\\
&~+\underbrace{\rho_k^2 \lVert \hat{\psi}(\mathbf p_k)-\epsilon_0(\mathbf p_k)\mathbbm 1-\nabla Z(\mathbf p_k)\rVert ^2}_{\Theta_{2,k}}+\underbrace{\rho_k^2 \lVert \nabla Z(\mathbf p_k) \rVert^2}_{\Theta_{3,k}}~~~\text{by the triangle inequality}
\label{eq_vk}
\end{align}

By Assumption \ref{assumption:properties}, $\Theta_{0,k}\leq 0$, so \eqref{eq_vk} gives
\begin{equation*}
    V(\mathbf p_{k+1},\mathbf p^*)\leq V(\mathbf p_{k},\mathbf p^*)+ \Theta_{1,k}+\Theta_{2,k}+\Theta_{3,k}
\end{equation*}

Let $\mathcal F_k$ be the filtration generated by $\{\mathbf p_1,\mathbf p_2,\ldots,\mathbf p_k\}$. Then given $\mathcal F_k$, 
\begin{align}
E[V(\mathbf p_{k+1},\mathbf p^*)| \mathcal F_k] &\leq V(\mathbf p_k,\mathbf p^*)+E[\Theta_{1,k}|\mathcal F_k]+E[\Theta_{2,k}|\mathcal F_k]+E[\Theta_{3,k}|\mathcal F_k]
\label{eq:evol_V}
\end{align}

We analyze the RHS of \eqref{eq:evol_V} term by term.

For the second term,
\begin{align}
    \lvert E[\Theta_{1,k}|\mathcal F_k]\rvert &= \rho_k\lvert E[(\hat{\psi}(\mathbf p_k)-\epsilon_0(\mathbf p_k)\mathbbm 1-\nabla Z(\mathbf p_k))'(\mathbf p^* - \mathbf p_k)|\mathcal F_k]\rvert \nonumber\\
    &=\rho_k E[\lvert (\psi(\mathbf p_k)-\epsilon_0(\mathbf p_k)\mathbbm 1-\nabla Z(\mathbf p_k))'(\mathbf p^* - \mathbf p_k) \rvert|\mathcal F_k]\nonumber\\
    &\leq \rho_k E[\lVert  \psi(\mathbf p_k)-\epsilon_0(\mathbf p_k)\mathbbm 1-\nabla Z(\mathbf p_k)\rVert \cdot \lVert \mathbf p^* - \mathbf p_k\rVert|\mathcal F_k]\nonumber\\
    &\leq \rho_k E[ B_1c_k\cdot 2\sqrt{m} |\mathcal F_k]=2B_1\sqrt{m}c_k\rho_k
    \label{eq:bd_mdsa1}
\end{align}

For the third term, we first get
\begin{equation}
\lVert \hat{\psi}(\mathbf p_k)-\epsilon_0(\mathbf p_k)\mathbbm 1-\nabla Z(\mathbf p_k)\rVert ^2 \leq 2\lVert \psi(\mathbf p_k)-\epsilon_0(\mathbf p_k)\mathbbm 1-\nabla Z(\mathbf p_k)\rVert ^2 + 2\lVert \psi(\mathbf p_k) - \hat{\psi}(\mathbf p_k)\rVert^2
\label{eq:bd_mdsa3_0}
\end{equation}
and under the assumptions on the bias and variance of the gradient estimator, \eqref{eq:bd_mdsa3_0} gives 
\begin{equation}
    \lvert E[\Theta_{2,k}|\mathcal F_k]\rvert \leq 2\rho_k^2 (B_1^2c_k^2+B_2\frac{\gamma_k}{R_kc_k^2})
    \label{eq:bd_mdsa2}
\end{equation}

For the last term, since $\nabla^2 Z$ is $L_2$-Lipschitz and $\prod_{i=1}^m \mathcal P^{n^i}$ is bounded, $M_1 = \sup_{\mathbf p \in \prod_{i=1}^m \mathcal P^{n^i}} \lVert \nabla Z(\mathbf p ) \rVert<\infty$, and we have 
\begin{equation}
    \lvert E[\Theta_{3,k}|\mathcal F_k]\rvert \leq \rho_k^2M_1^2
    \label{eq:bd_mdsa3}
\end{equation}

Telescoping \eqref{eq:evol_V}, and together with the bounds in \eqref{eq:bd_mdsa1}, \eqref{eq:bd_mdsa2}, and \eqref{eq:bd_mdsa3}, we get
\begin{align}
\sum_{k=1}^{\infty} E[E[V(\mathbf p_{k+1},\mathbf p^*)-V(\mathbf p_{k},\mathbf p^*)|\mathcal F_k]^{+}] &=O(\rho_kc_k)+O(\rho_k^2)+O(\frac{\rho_k^2\gamma_k}{c_k^2R_k})+O(\rho_k^2c_k^2)\nonumber\\
&\leq \sum_{k=1}^{\infty}\frac{C_1}{k^{\alpha +\theta}}+\frac{C_2}{k^{2\alpha}}+\frac{C_3}{k^{2\alpha -2\theta-2d+\beta}}+\frac{C_4}{k^{2\alpha +2\theta}}< \infty\nonumber
\end{align}
for some $C_1, C_2,C_3,C_4\geq 0$.

Thus, by the martingale convergence theorem as stated in Lemma \ref{lm:ec21-cor}, $V(\mathbf p_k,\mathbf p^*)$ converges a.s. to a random variable $V_{\infty}$. Next, we will argue that $V_{\infty} = 0$ a.s.. To show this, we take the expectation of \eqref{eq_vk},
\begin{align}
E[V(\mathbf p_{k+1},\mathbf p^*)] &\leq E[V(\mathbf p_k,\mathbf p^*)]+E[\Theta_{0,k}+\Theta_{1,k}+\Theta_{2,k}+\Theta_{3,k}]
\label{eq:expectation}
\end{align}
Telescoping \eqref{eq:expectation} we get 
\begin{equation}
E[V(\mathbf p_{k+1},\mathbf p^*)] \leq E[V(\mathbf p_1,\mathbf p^*)]+\sum_{i=1}^{k}\rho_iE[\nabla Z(\mathbf p_i)'(\mathbf p^* - \mathbf p_i)]+\sum_{i=1}^{k} E[\Theta_{1,k}+\Theta_{2,k}+\Theta_{3,k}]
\label{eq_tele}
\end{equation}


Sending $k\to \infty$ in \eqref{eq_tele}, the bounds in \eqref{eq:bd_mdsa1}, \eqref{eq:bd_mdsa2} and \eqref{eq:bd_mdsa3} give that 
\begin{equation*}
    \sum_{i=1}^{\infty} E[\Theta_{1,k}+\Theta_{2,k}+\Theta_{3,k}] =\sum_{i=1}^{\infty} E[E[\Theta_{1,k}+\Theta_{2,k}+\Theta_{3,k}|\mathcal F_k] ]<\infty
\end{equation*}

Using the fact that $V(\mathbf p_k,\mathbf p^*)\to V_{\infty}$ a.s. and $V(\cdot,\cdot)\geq 0$, we get
\begin{equation*}
\sum_{i=1}^{\infty} \rho_iE[\nabla Z(\mathbf p_i)'(\mathbf p_i-\mathbf p^* )] <\infty~a.s.
\end{equation*}

By Assumption \ref{assumption:properties}, $\nabla Z(\mathbf p_i)'(\mathbf p_i-\mathbf p^* )\geq 0$, and since $\sum_{i=1}^{\infty} \rho_i = \infty$ we get 
\begin{equation*}
\limsup_{k\to\infty} E[\nabla Z(\mathbf p_k)'(\mathbf p_k-\mathbf p^* )]=0
\end{equation*}
By an argument similar to the proof of Theorem \ref{main_FWSA}, we can find a subsequence $k_i$ such that $\nabla Z(\mathbf p_{k_i})'(\mathbf p_{k_i}-\mathbf p^* ) \to 0~a.s.$. By Assumption \ref{assumption:properties}, this implies that $\mathbf p_{k_i} \to \mathbf p^*~a.s.$ and so $V(\mathbf p_{k_i},\mathbf p^*)\to 0~a.s.$. Since we have already shown the a.s. convergence of $V(\mathbf p_k,\mathbf p^*)$, the limit must be identically 0. Thus, by Pinsker’s inequality we have $\mathbf p_k \to \mathbf p^*$ in total variation a.s.. This concludes the theorem.\hfill\Halmos
\endproof

\end{document}